\documentclass[lettersize,journal]{IEEEtran}
\usepackage{amsmath,amsfonts}
\usepackage{algorithmic}
\usepackage{algorithm}
\usepackage{array}
\usepackage[caption=false, font=footnotesize]{subfig}
\usepackage{textcomp}
\usepackage{stfloats}
\usepackage{url}
\usepackage{verbatim}
\usepackage{graphicx}
\usepackage{cite}
\usepackage{subcaption}

\usepackage{xspace}
\usepackage{hyperref}
\usepackage{cleveref}
\usepackage{caption}
\usepackage{siunitx}
\usepackage{tabularray}
\UseTblrLibrary{varwidth}

\usepackage[export]{adjustbox}
\usepackage{pgffor}
\usepackage{caption}
\usepackage{stmaryrd}
\usepackage{float}
\usepackage{enumerate}
\usepackage{xcolor}
\usepackage{scrhack} % necessary for listings package
\usepackage{listings}
\usepackage{dsfont} % blackboard-bold 1 \mathds{1}
\definecolor{brickred}{rgb}{0.8, 0.25, 0.33}
\lstdefinestyle{codeStyle}{
    showspaces=false,
    keywordstyle=\ttfamily\bfseries\color{purple},
    basicstyle=\fontsize{7}{7}\ttfamily,
    numbers=left,
    numberstyle=\tiny,
    commentstyle=\color{gray},
    breaklines=true,
    postbreak=\mbox{\textcolor{lightgray}{$\hookrightarrow$}\space},
    showstringspaces=false,
    stringstyle=\color{brickred},
    morekeywords={},
    escapeinside={(*}{*)}
}
\lstdefinestyle{PseudoPython}{
    language=python,
    morekeywords={function,do,new,mod,True,False}
}

\DeclareMathOperator{\pivot}{pivot}
\DeclareMathOperator{\im}{im}
\DeclareMathOperator{\Dgm}{Dgm}

\newcommand{\code}[1]{\texttt{#1}}

\begin{document}

\title{Efficient Betti Matching Enables Topology-Aware 3D Segmentation via Persistent Homology}

%\author{IEEE Publication Technology,~\IEEEmembership{Staff,~IEEE,}
        % <-this % stops a space
%\thanks{This paper was produced by the IEEE Publication Technology Group. They are in Piscataway, NJ.}% <-this % stops a space
%\thanks{Manuscript received April 19, 2021; revised August 16, 2021.}
%}

\author{Nico Stucki, Vincent Bürgin, Johannes C. Paetzold, Ulrich Bauer

\thanks{Nico Stucki, Vincent Bürgin, Ulrich Bauer and are with the Technical University of Munich, Germany. Johannes C. Paetzold is with Imperial College London, UK. Vincent Bürgin is supported by the DAAD programme Konrad Zuse Schools of Excellence in Artificial Intelligence, sponsored by the Federal Ministry of Education and Research.
N. Stucki and U. Bauer are supported by the Munich Data Science Institute (MDSI – ATPL4IS).}
}

%\markboth{IEEE Transactions on Pattern Analysis and Machine Intelligence}{}

\maketitle


\begin{abstract}
In this work, we propose an efficient algorithm for the calculation of the \textit{Betti matching}, which can be used as a loss function to train topology aware segmentation networks. 
Betti matching loss builds on techniques from topological data analysis, specifically persistent homology. 
A major challenge is the computational cost of computing persistence barcodes. 
In response to this challenge, we propose a new, highly optimized implementation of Betti matching, implemented in C++ together with a python interface, which achieves significant speedups compared to the state-of-the-art implementation Cubical Ripser. % \cite{kaji_cubical_2020}.
We use \textit{Betti matching 3D} to train segmentation networks with the Betti matching loss and demonstrate improved topological correctness of predicted segmentations across several datasets.
\thanks{ The source code is available at \href{https://github.com/nstucki/Betti-Matching-3D}{https://github.com/nstucki/Betti-Matching-3D}.}
\end{abstract}


\section{Introduction}

\IEEEPARstart{T}{opology}-aware segmentation is a growing field of research in Computer Vision. 
In recent years, numerous optimization functions have been developed to enhance the preservation of topology. 
This is a relevant task since topologically accurate segmentation is sometimes more important than pixel/voxel accurate segmentation, i.e., some voxels are more important than others.
For instance, a network aiming to detect cells may produce noisy predictions with false positives scattered in empty areas, skewing the cell count. 
The loss contribution from these false positives is marginal, and hence the signal to eliminate these predictions is too small. 
Here, it might be more beneficial to refine the edges of correctly detected cells rather than eliminating falsely created connected components. 
Similarly, in vessel segmentation, vessels appear as thin curvilinear structures that can be disconnected or wrongly connected by a few mispredicted voxels. 
A small error in the Dice score can significantly alter the topology. 

Some of the proposed approaches utilize persistent homology, while others are based on overlap techniques. 
However, persistent homology-based losses have not yet been effectively applied to three-dimensional (3D) data. 
The Betti matching framework \cite{stucki_topologically_2023} has shown promise in two-dimensional (2D) applications, and its theoretical principles generalize to 3D. 
Despite this, the current algorithms face significant runtime challenges, even when used in 2D. 

% \begin{itemize}
%     \item multiple optimization functions have been proposed which optimize topological structure preservation 
%     \item some approaches employ persistent homology others are overlap based 
%     \item so far no persisten homology based losses have been efficiently applied to 3D 
%     \item the recent Betti matching framework is useful in 2D and its theoretical aspects extend to 3D however the runtime of existing algorithms makes it already challengig to be used in 2D 
%     \item In this work we resolve this issue and propose an efficient algorithm implemented in C++ which allows the rapid calculation of the Betti matching loss in 2D and 3D 
%\end{itemize}

\paragraph*{Contributions}

In this paper, we address this issue by introducing an efficient algorithm implemented in C++. 
Our proposed solution enables the calculation of the Betti matching for inputs of arbitrary dimensions, with an highly optimized version in 1D, 2D and 3D contexts.
%Specifically, we ... 1) 2) ,

\subsection{Related works}
Our work builds on two main fields:
1) \textit{topological deep learning}, in particular topology-aware segmentation, and 2) \textit{computational topology}, in particular the theory of persistent homology and the literature on efficient algorithms for computing it.

\paragraph{Topology-Aware Segmentation and Topological Losses}
Our work extends \textbf{Betti matching} introduced in \cite{stucki_topologically_2023}. 
The work defines the Betti matching loss, which uses persistent homology and induced matchings between persistence barcodes to define a segmentation loss based on the lifetimes of matched topological features. 
The matching implicitly takes the spatial relationship of features into account. 
A precursor to the Betti matching loss is the \textbf{TopoNet loss} \cite{hu_topology-preserving_2019}, which was the first paper to propose using persistent homology for a segmentation loss.
The crucial difference to the Betti matching loss is that the TopoNet loss is based on the Wasserstein matching between persistence diagrams and hence matches features only based on their lifetimes (i.e.,  their local contrast differences between the feature and its background). 
Since features on one side of the matching come from a binary label map, they all have the same lifetime, and the loss essentially only measures if there is an equal number of topological features between prediction and ground truth, ignoring their spatial relation.

Another prominent topological loss function is the centerline Dice loss \textbf{clDice} \cite{shit_cldice_2021}, which is motivated by the tubular structure found in vessel networks. 
It computes approximate skeletonizations (centerlines) of the label and prediction and, in analogy to the precision and sensitivity metrics that the Dice loss is based on, defines the topological precision and topological sensitivity metrics based on the overlap of the prediction's skeleton with the non-skeletonized label and vice versa. 
The clDice score is defined as the harmonic mean between topological precision and topological sensitivity.

Other recent works have proposed several quite different approaches towards topologically correct predictions. 
In \cite{hu_structure-aware_2022}, prediction and ground truth masks are warped in a homotopy-type preserving way to resemble each other, and the volumetric loss is amplified on those critical pixels which cannot be included without changing the homotopy type.
The uncertainty quantification method of \cite{gupta_topology-aware_2023} uses discrete morse theory to capture areas of high topological uncertainty such that they can be brought to the attention of human experts in an active learning scenario.
In \cite{li_robust_2023}, Euler numbers of local patches are computed and compared between the prediction and label, which is used to form a topology violation map that serves as the input signal to a second prediction network.
A multiclass segmentation method based on domain knowledge, notably working on 3D data, is proposed in \cite{gupta_learning_2022}. 
It enforces problem-specific topological properties by forbidding certain classes to be adjacent in the prediction:
For example, the aortic lumen being enclosed by the aortic wall is expressed by a rule that the aortic lumen class should never be adjacent to the background class.

The recent survey and position paper \cite{papamarkou_position_2024} gives an overview of the field of topological deep learning, including many other directions apart from segmentation losses, and identifies open challenges.

\paragraph{Persistent Homology and Computational Topology}
The notion of persistent homology was introduced in \cite{edelsbrunner_topological_2002}, even though similar ideas had been independently introduced previously several times  (see \cite{edelsbrunner_persistent_2017}).
% The survey \cite{edelsbrunner_24_2017} attributes the attention the topic got following  \cite{edelsbrunner_topological_2002} to the introduced algorithm which allowed to turn the theoretical notion into a practical computational tool.  
Persistent homology is often used to analyse filtrations of simplicial complexes, for example arising from geometric constructions on point cloud inputs.
However, it can be equally well applied to cubical complexes, a natural candidate to model the topology of digital images. Cubical persistent homology is extensively described in \cite{kaczynski_computational_2004}.

An important line of work in the persistent homology literature is concerned with algorithmic improvements, which over time have sped up the computation times for persistence barcodes by several orders of magnitude (cf. \cite[Table 1]{bauer_ripser_2021}). 
The basic algorithm used to compute barcodes is a restricted form of Gaussian elimination of a matrix, the boundary matrix, derived from the topological complex.
The boundary matrix has a lot of structure that can be exploited.
The first proposal of this algorithm in \cite{edelsbrunner_topological_2002} already used a sparse representation of the matrix (see \cite{edelsbrunner_short_2014}). 
Later optimizations include the clearing optimizations \cite{chen_persistent_2011}, the usage of cohomology instead of homology \cite{silva_dualities_2011}, and the (non-trivial) parallelization \cite{bauer_distributed_2013}.
The software \textit{Ripser} \cite{bauer_ripser_2021} introduces the technique of implicitly representing the boundary matrix and the apparent and emergent pair optimizations (combined with the earlier clearing and cohomology techniques), which together allows major speedups over previous software packages. 
\textit{Ripser++} \cite{zhang_gpu-accelerated_2020} realizes massive parallelization of parts of the Ripser algorithm on the GPU (namely the filtration construction, clearing, and identification of apparent pairs). 
\textit{Oineus} \cite{nigmetov_oineus_2021} introduces a lock-free CPU-based parallelization of the full Ripser reduction algorithm (not only of some parts as Ripser++). 
Building on this parallelized algorithm, \textit{giotto-ph} \cite{perez_giotto-ph_2021} is another recent implementation and incorporates further optimizations.
While these implementations works on simplicial complexes, \textit{Cubical Ripser} \cite{kaji_cubical_2020} adapts the algorithm for three-dimensional cubical complexes which correspond to 3D image inputs.

The Betti matching relies on computing so-called image barcodes, barcodes of the image of a morphism between two persistence modules. 
An algorithm for computing image persistence was first presented in \cite{cohen-steiner_persistent_2009}, and generalized in \cite{bauer_efficient_2022}.

The idea of optimizing objective functions based on persistent homology is not restricted to the computer vision or medical imaging literature, see e.g. \cite{carriere_optimizing_2021}, \cite{nigmetov_topological_2024}. 
The latter is a recent work that describes the general setting of a quadratic loss computed between matched points in a persistence diagram, a setting that encompasses both the TopoNet loss and the Betti matching loss as special cases.
While the standard way of optimizing such loss functions involves prescribing gradients to two simplices (or cubes) in a topological complex, the paper proposes a more efficient optimization procedure that propagates gradients to a larger region of the input.

\section{Persistent Homology and Betti matching}

\begin{figure}[!t]
    \centering
    \subfloat
    {%
        \centering
        \includegraphics[height=0.6\linewidth]{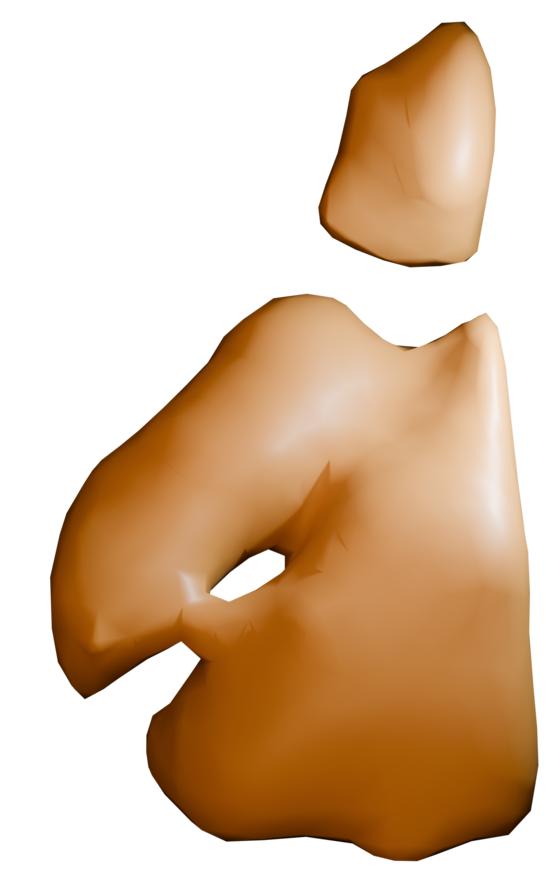}
    }\hfill
    \subfloat
    {%
        \centering
        \includegraphics[height=0.6\linewidth]{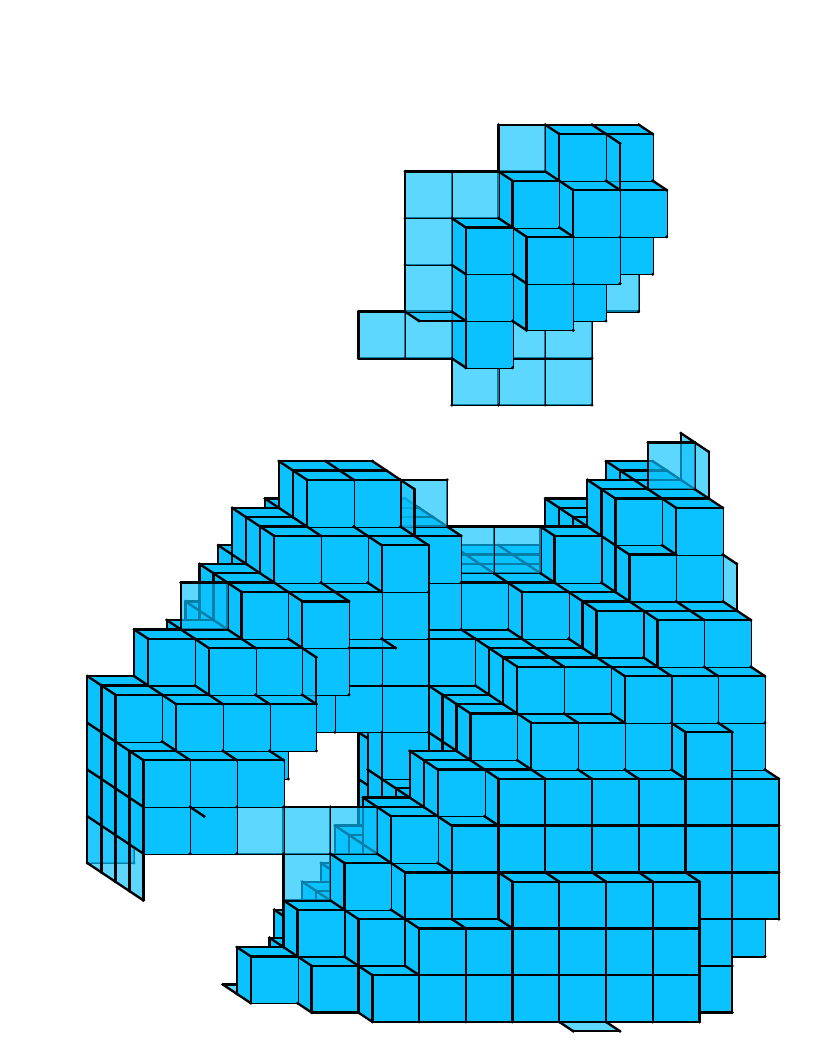}
    }
    \caption{A predicted segmentation 97\% certainty superlevel set on a patch from the BBBC027 cell dataset that has two connected components and a loop, in smoothed surface and cubical grid complex representations.
    }
    \label{fig:persistent-homology-bbbc027-intro-example}
\end{figure}

We first introduces the necessary mathematical background from topological data analysis, namely the homology of cubical complexes, persistent homology, and induced matchings, and use these foundations to recap the definition of Betti matching between grayscale digital images \cite{stucki_topologically_2023}. 
We will then go into the algorithms to compute barcodes and Betti matching. 
We keep the description of the theory as specific as possible to our setting, that is, the persistent homology of \textit{digital images} in \textit{three dimensions} with \textit{$\mathbb F_2$ coefficients}. 
This restriction simplifies the presentation and allows us to make specific algorithmic optimizations.

As an introductory example, consider Figure \ref{fig:persistent-homology-bbbc027-intro-example}. 
The left picture shows a segmentation network prediction on the BBBC027 dataset, with the voxel-based prediction displayed as a smoothed surface.
On the right, we see the same data represented as a cubical grid complex. 
This representation is finer than the voxel representation and is made up of vertices, edges, faces, and cubes.
The topological information we can get from the picture is that there are two connected components, one of which contains a hole. 
This example shows some of the information we can intuitively capture with homology: connected components correspond to classes in the homology of dimension zero and holes correspond to classes in homology of dimension one. 
We can furthermore capture cavities, i.e., empty space enclosed by faces, which correspond to classes in homology of dimension two.

The complex in Figure \ref{fig:persistent-homology-bbbc027-intro-example} is not the whole picture: 
It only represents the topology of the input image at a certain threshold $\alpha$. 
As we vary $\alpha$ and gradually take more voxels into consideration, new structures might come into existence and existing ones might vanish.
We can track the lifespan of each structure by \textit{persistent homology}. 
The idea is to deem the long-living features important, as they have high contrast, and the short-lived ones unimportant or likely to stem from noise.

\subsection{Cubical Complexes}
\label{sec:persistent-homology-cubical-sets}

\begin{figure}[!t]
    \centering
    \includegraphics[width=0.8\linewidth]{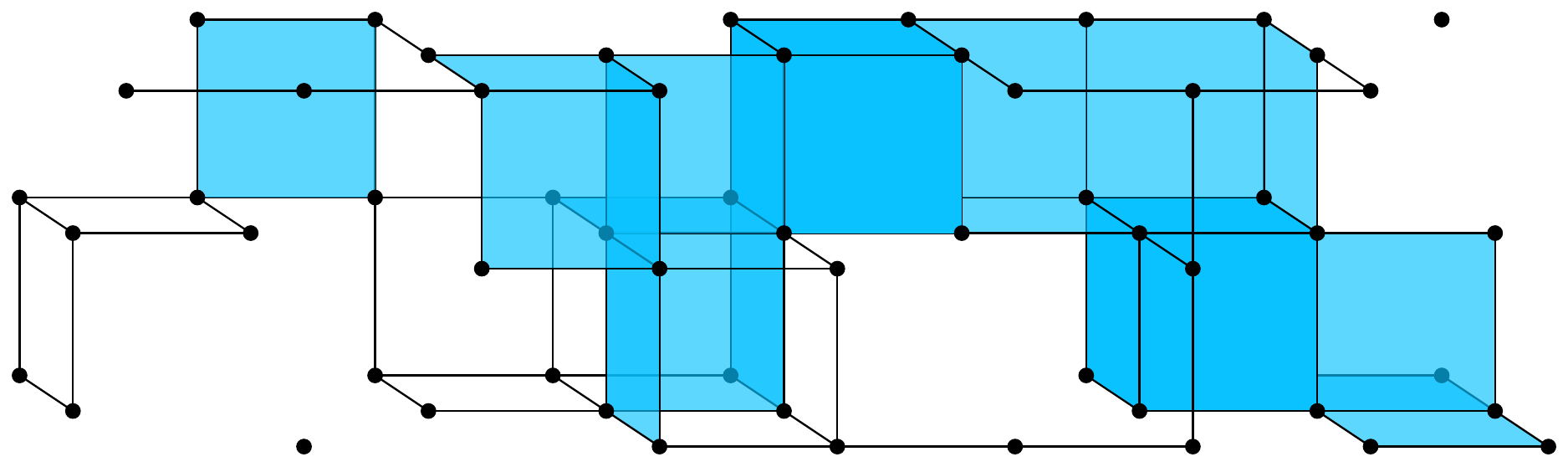}
    \caption{A cubical grid complex consisting of vertices, edges, 2-cubes and 3-cubes}
    \label{fig:cubicalGridComplex}
\end{figure}

We represent 3D greyscale images as \textit{cubical complexes}, which are composed of \textit{cubical cells}. 
A \textit{cubical cell}, or \textit{cube}, is a set $c$ formed by a cartesian product of the form
\begin{equation*}
    c = [k_1, k_1 + \delta_1] \times [k_2, k_2 + \delta_2] \times [k_3, k_3 + \delta_3],
\end{equation*}
for $k_i \in \mathbb Z, ~\delta_i \in \{0, 1\}.$
The cube extends in $d = |\{i \mid \delta_i \neq 0\}|$ dimensions and we say that $c$ is a $d$-dimensional cube (\textit{$d$-cube}). 
We also call $0$-cubes \textit{vertices} and $1$-cubes \textit{edges} (see Figure \ref{fig:cubicalGridComplex} for an example). 
If $c_1, c_2$ are cubes of dimension $d_1, d_2$ and $c_1 \subseteq c_2$, we call $c_1$ a \textit{face} of $c_2$ (denoted by $c_1 < c_2$) of \textit{codimension} $d_2 - d_1$, and $c_2$ a \textit{co-face} of $c_1$. 
If the codimension is one, we also call $c_1$ a \textit{facet} of $c_2$.
A \textit{cubical complex} is a set of cubes that is closed under the face relation: 
If a cube is included, all its faces must be included (\cite{stucki_topologically_2023}, cf. \cite[Definition 2.9]{kaczynski_computational_2004}).
They arise naturally from digital image data, which is the application we build toward: 
There are two common constructions, representing voxels either as 0-cubes or as 3-cubes, and we will postpone the details to \cref{sec:persistent-homology-digital-images}.

An arbitrary collection of $d$-cubes constitutes a $d$-dimensional \textit{chain}. Despite the name, chains do not need to be connected.
We can define chains as arbitrary sets of cubes, and two chains can be combined by taking the symmetric difference (exclusive or) between the sets. 
A different, but equivalent way of viewing chains is not as sets, but as vectors: 
The chains form a vector space $C_d$ over the field $\mathbb F_2$ and the space is generated by the set of $d$-cubes. 
That is, the cubes are seen as linearly independent vectors and can be combined via addition to form chains. When we add chains duplicate cubes cancel out (equivalent to the exclusive-or operation). 
We will not use the set-based formulation but only use the vector space view from here, since all the concepts that follow can be expressed conveniently in the linear algebra framework. (\cite{edelsbrunner_short_2014}, \cite{kaczynski_computational_2004})

The \textit{boundary} of a $d$-cube $c$ is a $(d-1)$-chain $\partial_d \,c$ consisting of the facets of $c$.
This definition extends linearly to a \textit{boundary map} $\partial_d: C_d \to C_{d-1}$.
A \textit{cycle} is a chain $c$ with empty boundary: $\partial_d\,c = 0$. 
Since $\partial$ is linear, the sum of two cycles forms a cycle.
The space of $d$-dimensional cycles $Z_d$ again forms a vector space, a subspace of $C_d$. (\cite{edelsbrunner_short_2014}, \cite{kaczynski_computational_2004})

Another distinguished subset of chains are the $d$-dimensional \textit{boundaries} $B_d$, which are those $d$-chains which are the boundary of some $(d+1)$-chain in $C_{d+1}$. Every boundary is a cycle, and since $B_d$ is closed under addition, $B_d$ is a subspace of $Z_d$. (\cite{edelsbrunner_short_2014}, \cite{kaczynski_computational_2004})

\subsection{Homology}
\label{sec:homology}

\textit{Homology} is an equivalence class construction that captures, intuitively speaking, whether cycles go around the same hole. In a topological complex, there may be a large number of cycles, but a much smaller number of holes, and grouping together cycles by the hole they encompass lets us identify the holes. \cite{edelsbrunner_short_2014}

In homology, we consider two cycles $z_1, z_2 \in Z_d$ equivalent, or \textit{homologous}, if they only differ by a boundary $b \in B_d$: That is, $z_1 - z_2 = b$. 
In this case, they are both represented by the same \textit{homology class} $[z_1] = [z_2]$. 
Since we use coefficients in $\mathbb F_2$, this condition is equivalent to $z_1 + z_2 = b$. 
In linear algebra terms, the homology $H_d$ in dimension $d$ of a cubical complex is the quotient space $H_d = Z_d / B_d$. 
The dimension of $H_d$ is denoted $\beta_d$ and called the \textit{d-th Betti number} of the complex 
(\cite{edelsbrunner_short_2014}, \cite{kaczynski_computational_2004})

The information that homology can capture is for example connectivity in dimension zero, loops in dimension one, and cavities in dimension two. 
In our $3$D image setting, the homology $H_3$ is trivial and contains no information, as there are no $3$-cycles.

\subsection{Persistent Homology}
\label{sec:persistence}

A \emph{filtration} of a cubical complex $K$ is a function $f \colon K \to \mathbb{R}$ that satisfies $f(c_1) < f(c_2)$ whenever $c_1 < c_2$.
It enables us to filter $K$ by \emph{sublevel sets} $D(f)_t = f^{-1}((-\infty,t])$, which form subcomplexes of $K$.

For $s \leq t$, the inclusion $D(f)_s \hookrightarrow D(f)_t$ extends linearly to a map $Z_d(D(f)_s) \to Z_d(D(f)_t)$ between the cycles, which descends to a map $h_d^{s,t}: H_d(D(f)_s) \to H_d(D(f)_t)$ in homology. 
The family $\{H_d(D(f)_t)\}_t$ together with the maps $\{h^d_{s,t}\}_{s,t}$ forms a \textit{persistence module} $H_d(f)$, the \emph{persistent homology} of the filtration. (\cite{edelsbrunner_short_2014}, \cite{stucki_topologically_2023})

Since any \emph{pointwise finite dimensional} persistence module admits a barcode decomposition \cite{crawleyboevey2014decompositionpointwisefinitedimensionalpersistence}, assuming that $K$ is finite, there exists a multiset $\mathcal{B}_d(f)$ of intervals, such that $H_d(f)$ decomposes into a direct sum of $\emph{interval modules}$. \cite{stucki_topologically_2023}
The barcode encodes \emph{birth} and \emph{death} of homological features in the filtration and how they relate to each other.

\begin{figure*}[!h]
    \centering

    \resizebox{\textwidth}{!}{
    \subfloat[$\scriptstyle r = 0.33$]{
        \centering
        \includegraphics[width=0.19\linewidth]{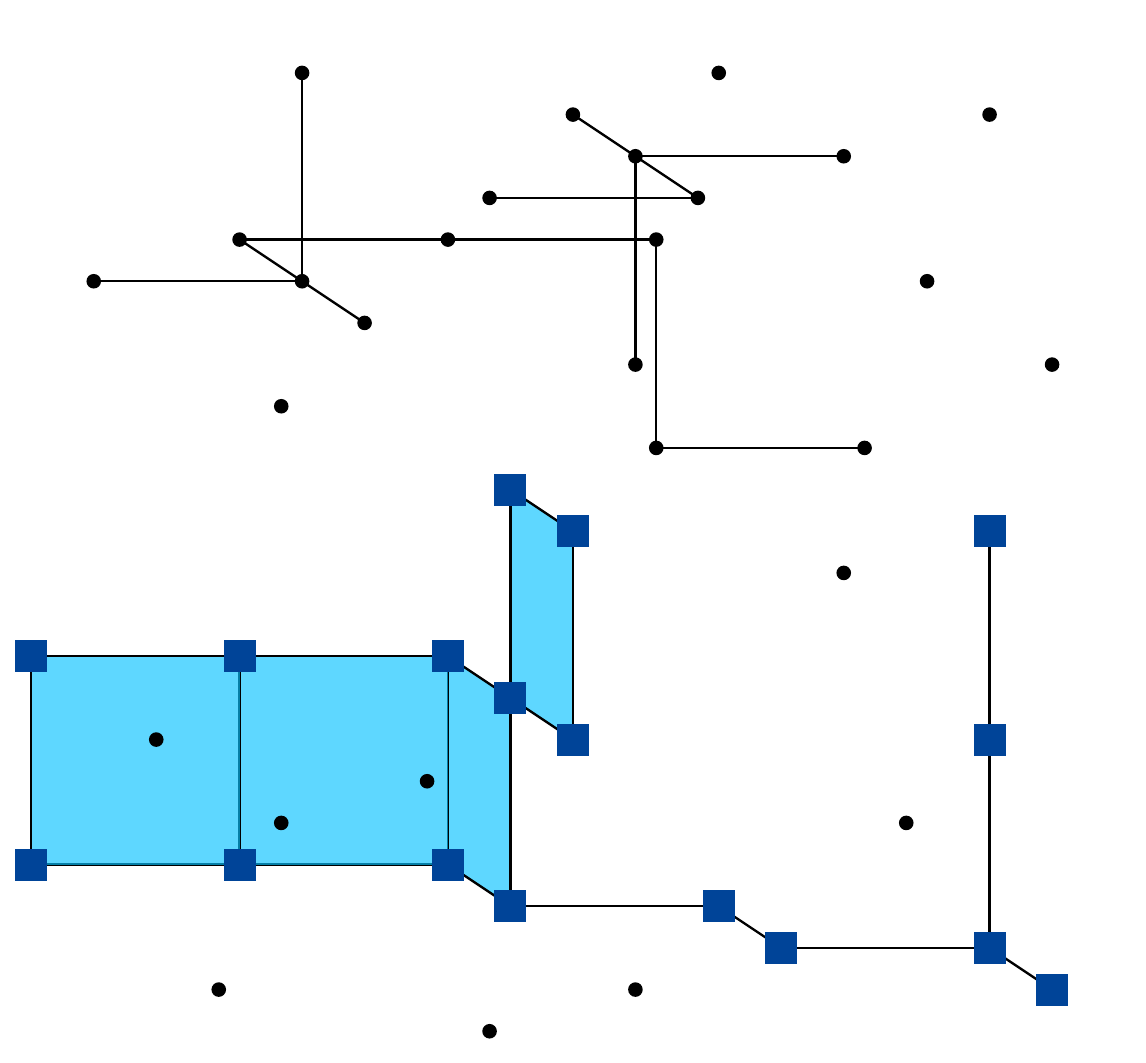}
    }
    \subfloat[$\scriptstyle r = 0.4$]{
        \centering
        \includegraphics[width=0.19\linewidth]{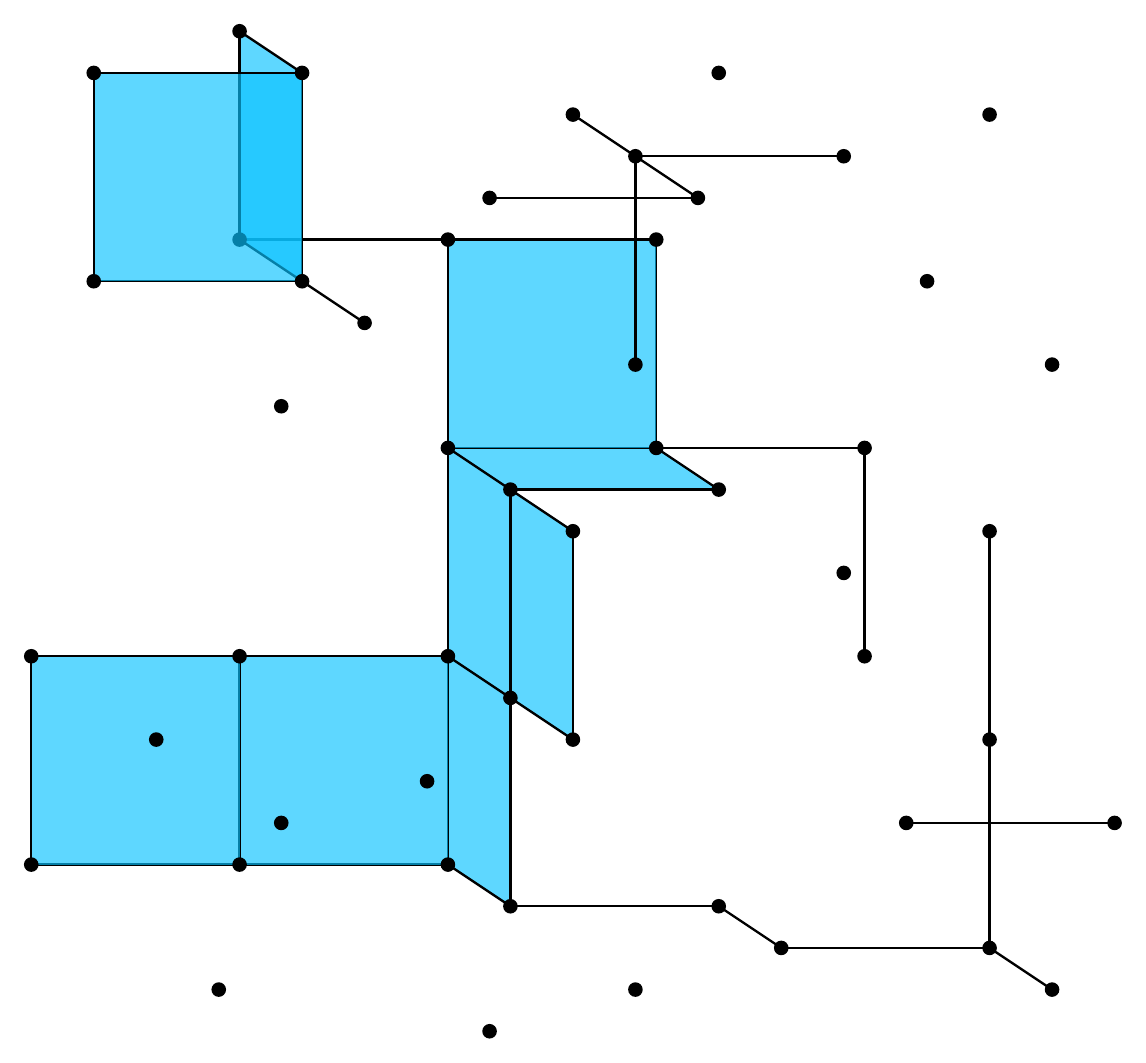}
    }
    \subfloat[$\scriptstyle r = 0.55$]{
        \centering
        \includegraphics[width=0.19\linewidth]{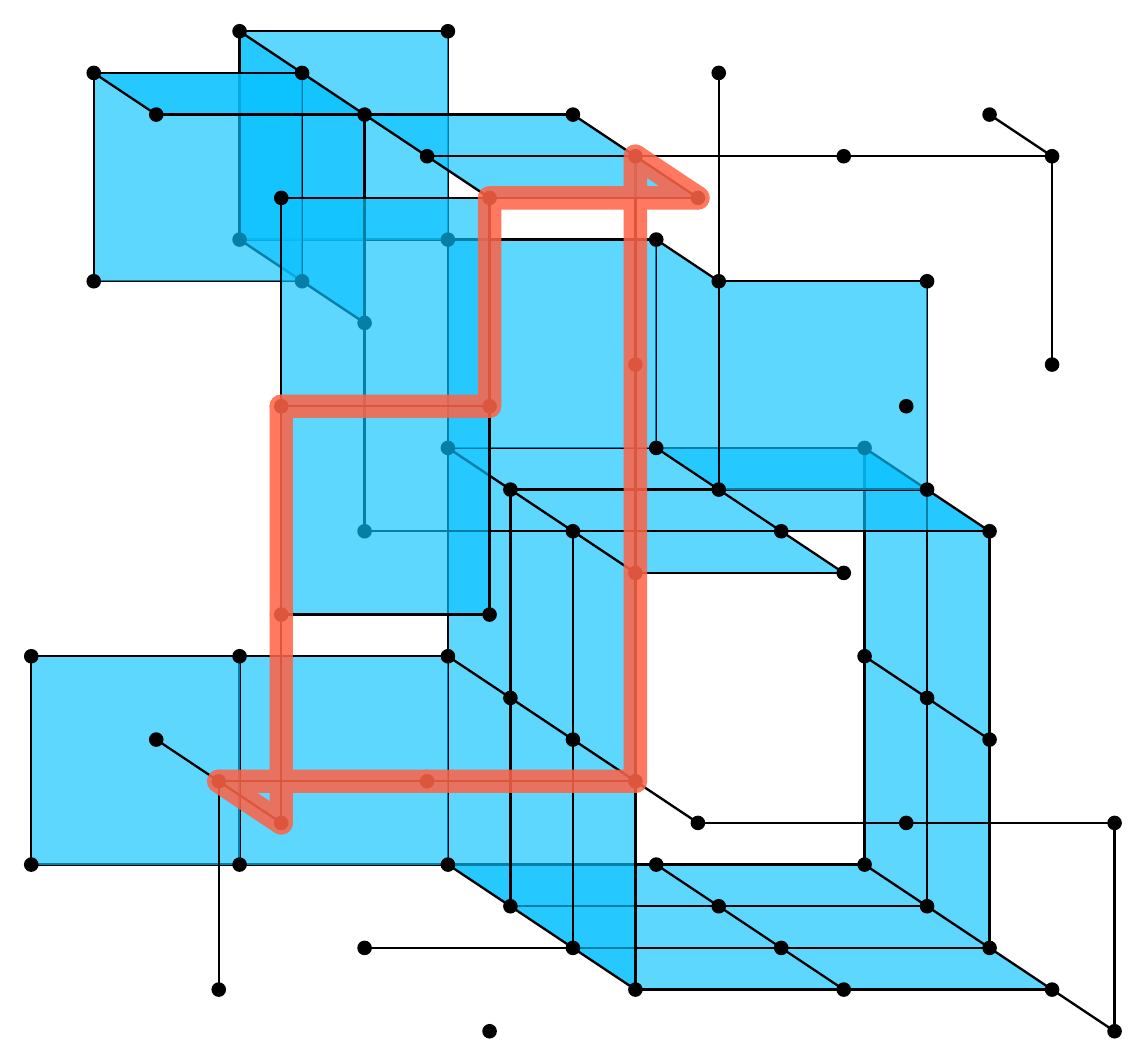}
    }
    \subfloat[$\scriptstyle r = 0.66$]{
        \centering
        \includegraphics[width=0.19\linewidth]{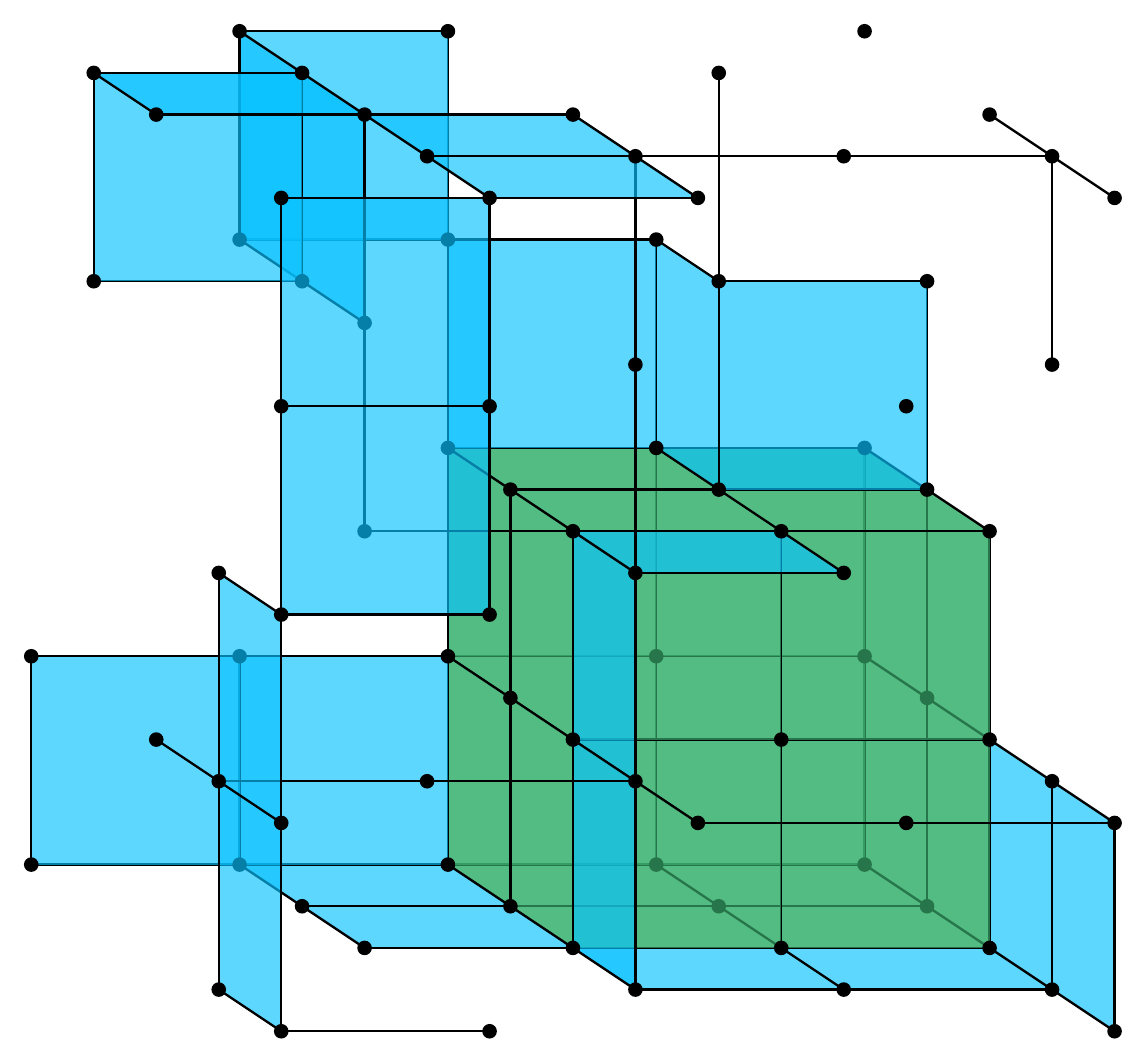}
    }
    \subfloat[$\scriptstyle r = 0.9$]{
        \centering
        \includegraphics[width=0.19\linewidth]{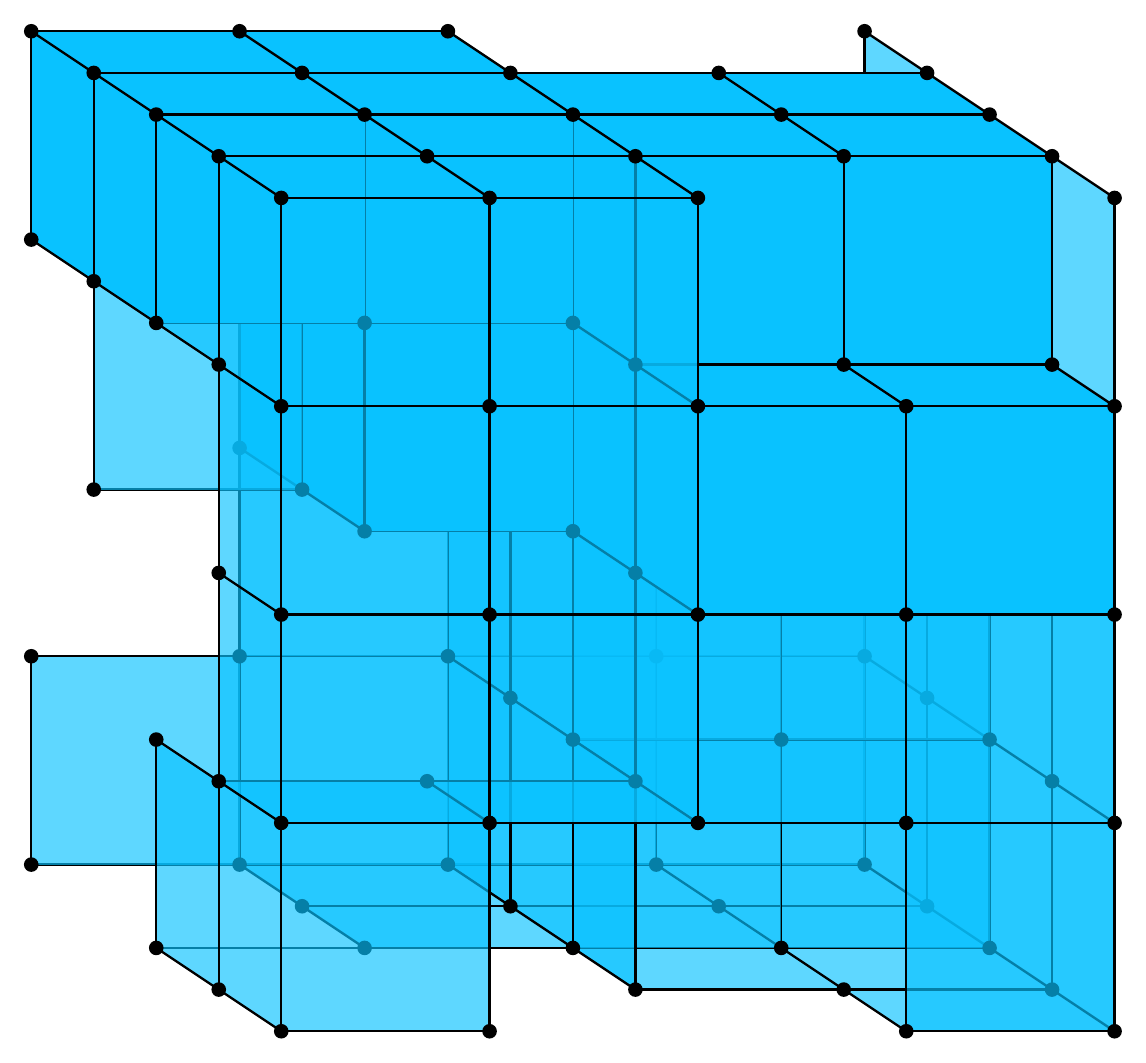}
    }%
    }
    % \\
    
    \hfill
    \resizebox{0.75\textwidth}{!}{
    \subfloat[Barcode]{
        \includegraphics[width=.3\textwidth]{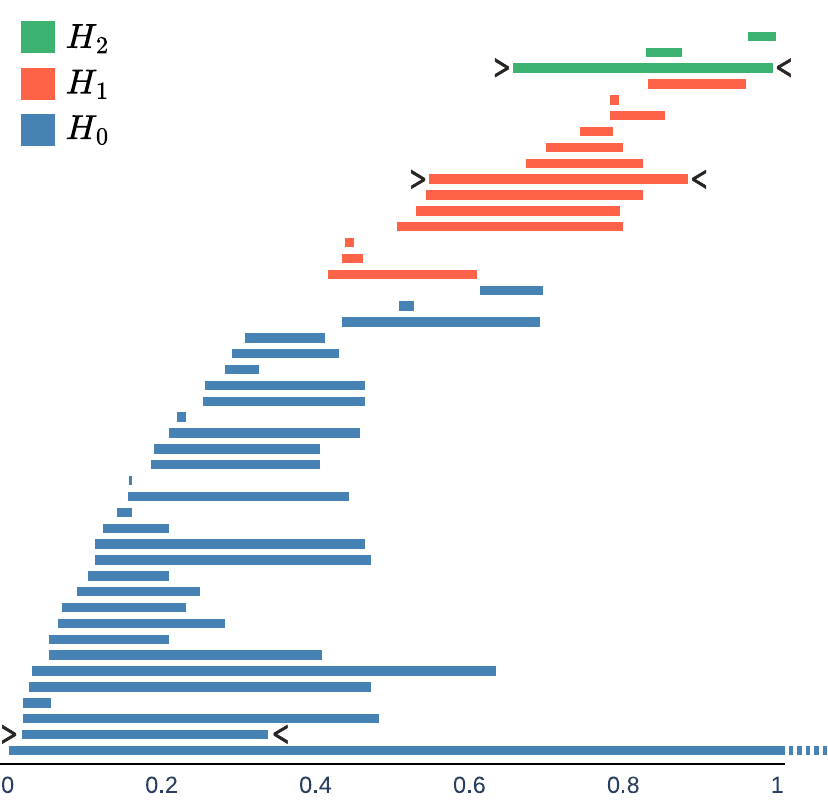}    
    }
    \hspace{0.1\textwidth}
    \subfloat[Persistence diagram]{
        \includegraphics[width=.3\textwidth]{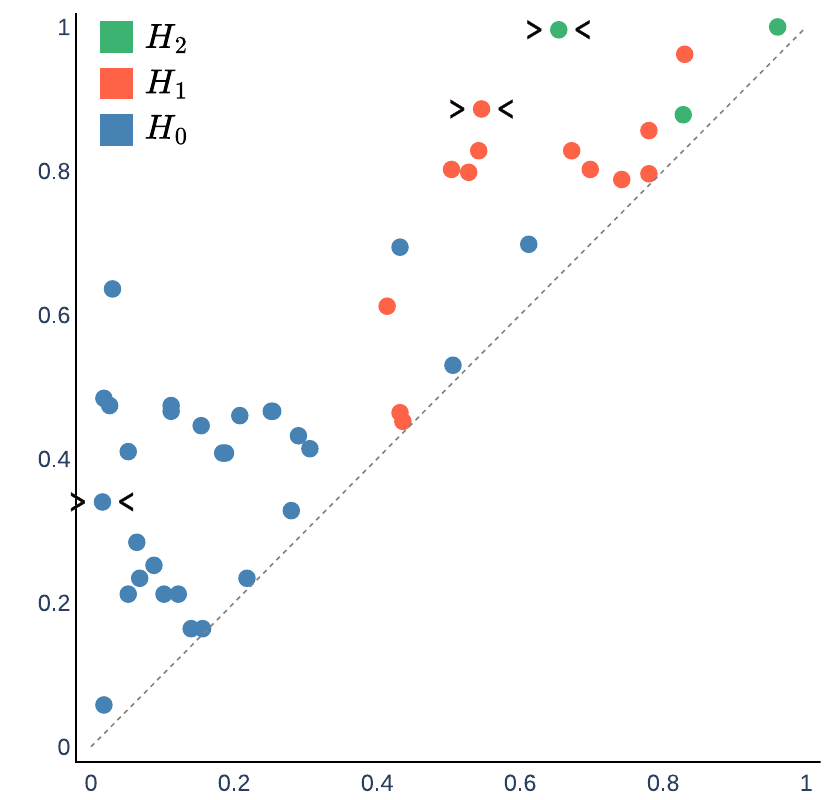}    
    }
    % \vspace{0.4cm}
    % \hfill
    % \begin{subfigure}{0.49\textwidth}
    %     \centering
    %     \includegraphics[width=\textwidth]{figures/exampleVI-barcode/barcode.pdf}
        
    %     % \vspace{0.045cm}
    %     \caption{Barcode}
    % \end{subfigure}
    % \hfill
    % \begin{subfigure}{0.49\textwidth}
    %     \centering
    %     \includegraphics[width=\textwidth]{figures/exampleVI-barcode/persistence_diagram.pdf}
    %     \caption{Persistence diagram}
    % \end{subfigure}
    % \hfill
    % {}
    % }
    }
    \hfill
    {}
    
    \caption{
    Barcode and persistence diagram of an example complex. We highlight one feature each from $H_0$, $H_1$, $H_2$ at different steps in the filtration, and the corresponding bars and points in the barcode and persistence diagram.
    }
    \label{fig:barcode}
\end{figure*}

Another representation for the persistent homology is its \textit{persistence diagram} $\Dgm_d(f) = \{(b, d) \mid [b,d) \in \mathcal{B}_d(f)\} \cup \{(x,x) \mid x \in \mathbb R\}$.
It records the persistence intervals as points in $\mathbb R^2$ above the diagonal $x=y$, augmented by infinitely many points on the diagonal. 
Figure \ref{fig:barcode} shows the barcode and persistence diagram of a 3D cubical complex.

\subsection{Persistent Homology of Digital Images}
\label{sec:persistent-homology-digital-images}

We represent a 3D grayscale image $\mathbf I \in \mathbb R^{N_1 \times N_2 \times N_3}$ by the cubical complex $K^{\mathbf I}$ that consists of all cubical cells contained in $[1, N_1] \times [1, N_2] \times [1, N_3]$. 
Voxels of the image correspond to vertices in $K^{\mathbf I}$  (known as the \textit{V-construction}).
The voxel intensity values are represented via a filtration $f_{\mathbf I}: K^{\mathbf I} \to \mathbb R$. 
At vertices $f_{\mathbf I}$ takes the value of the corresponding voxel and at higher-dimensional cubes $f_{\mathbf I}$ takes the \textit{maximum} of the values that it takes at its faces. \cite{stucki_topologically_2023}

We can now define the persistent homology of image $\mathbf I$
to be $H_d(\mathbf I) = H_d(D(f_{\mathbf I}))$ and its barcode $\mathcal{B}(\mathbf I) = \mathcal{B}(f_{\mathbf I})$.
An interval, or \textit{bar}, represents the span of threshold values within which the topological feature exists; we can interpret the length of the bar as the amount of contrast of the feature.

Following the standard algorithm for the barcode-computation of persistent homology (see \cite{edelsbrunner2008persistent})
we refine the filtration to a \emph{cube-wise} filtration by choosing a \emph{compatible total ordering} $c_1,\ldots,c_N$ of the cubes.
In this work, we will tie-break lexicographically first on the cube's dimension $d$, then the coordinates $(x, y, z)$ of the lexicographically smallest vertex among its faces, and then on its \textit{type}. 
We define the type as $0$, $1$ or $2$ for edges depending on if the edge extends in $x$, $y$ or $z$ direction, and $0, 1$ or $2$ for 2-cubes depending on if the 2-cube extends in $y$-$z$, $x$-$z$, or $x$-$y$ direction. 
Furthermore, cubes of dimension $0$ and $3$ are all of type $0$.

The algorithm will output pairs $(c_i,c_j)$ of cells, so called \emph{persistence pairs}, such that the barcode $\mathcal B(\mathbf I)$ consists of \emph{finite} intervals of the form
$[f_{\mathbf I}(c_i),f_{\mathbf I}(c_j))$ for persistence pairs $(c_i, c_j)$ that satisfy $f_{\mathbf I}(c_i) < f_{\mathbf I}(c_j)$ and \emph{essential} intervals of the form
$[f_{\mathbf I}(c_i),\infty)$ for cubes $c_i$ that are unpaired.

For the purpose of matching intervals by the means of induced matchings, we will make use of the tie-breaking on the cells, meaning that we record the obtained persistence pairs in the \emph{refined barcode} \[\mathcal B^{\text{fine}}(\mathbf I) = \{[i,j) \mid (c_i,c_j) \text{ is persistence pair}\}.\]

In the refined barcode and its dimension-specific versions, no number occurs as endpoint of a bar more than once. 
This means the respective persistence modules are so-called \textit{staggered persistence modules}. 
We will switch between these perspectives as needed: 
In the context of induced matchings, we use refined barcodes, whereas for visualization and defining the Betti matching loss, the non-refined barcodes/persistence diagrams are more useful.

\subsection{Induced Matchings and Betti Matching}
\label{sec:induced-matchings-betti-matching}

We are now ready to define induced matchings and image barcodes, which will allow us to define the Betti matching. 
The high-level idea is to match intervals in the barcodes of two grayscale images $\mathbf I, \mathbf J$ by defining a \textit{comparison image} $\mathbf C$ as their point-wise minimum and matching features between $\mathbf I$, $\mathbf J$ via matchings to the barcode of $\mathbf C$ induced by inclusions.

Consider two grayscale images $\mathbf I$, $\mathbf C$ of the same size such that $\mathbf I \geq \mathbf C$ (point-wise), and denote the corresponding cubical complexes by $K^{\mathbf I}$, $K^{\mathbf C}$. 
The inclusions $D(f_{\mathbf I})_s \hookrightarrow D(f_{\mathbf C})_s$ induce maps $\Phi_s \colon H_d(D(f_{\mathbf I})_s) \to H_d(D(f_{\mathbf C})_s)$ in homology that ensemble to a \emph{morphism} $\Phi \colon H_d(\mathbf I) \to H_d(\mathbf C)$ of persistence modules.
The \textit{image} of $\Phi$ forms another persistence module, denoted by $\im \Phi$ and defined as $(\im \Phi)_s = \im (\Phi_s)$. 
Intuitively, $\im(\Phi_s)$ captures topological features of $\mathbf C$ which also already have an equivalent and spatially corresponding feature in $\mathbf I$ at parameter $t$. 
The spatial correspondence comes from the fact that $\Phi_s$ is induced by the inclusion map.
The persistence module $\im \Phi$ comes with \emph{image persistence pairs} $(c_i, c_j)$, where $i$ refers to the cube-wise refinement with respect to $\mathbf I$ and $j$ refers to the cube-wise refinement with respect to $\mathbf J$.
They give rise to the \textit{image barcode} $\mathcal B(\im \Phi)$ consisting of finite intervals $[f_{\mathbf I}(c_i),f_{\mathbf C}(c_j))$ for each image persistence pair that satisfies $f_{\mathbf I}(c_i)<f_{\mathbf C}(c_j)$ and essential intervals $[f_{\mathbf I}(c_i),\infty)$ for each unpaired cell $c_i \in K^{\mathbf I}$. \cite{stucki_topologically_2023}

\let\exampleviifigurewidth\linewidth
\begin{figure}[!t]
    \centering
    \subfloat[Matched and unmatched cycles between $\mathbf{I}$, $\mathbf{C}$ and $\mathbf{J}$ 
    % \label{subfig:betti-matching-example-3d}
    ]{
        \begin{tblr}{colspec={X[1,c,m]X[30,c,h]X[30,c,h]X[30,c,h]X[30,c,h]X[30,c,h]},width=\exampleviifigurewidth,
        measure=vbox}
            $\mathbf{I}$  &
            \includegraphics[width=\linewidth]{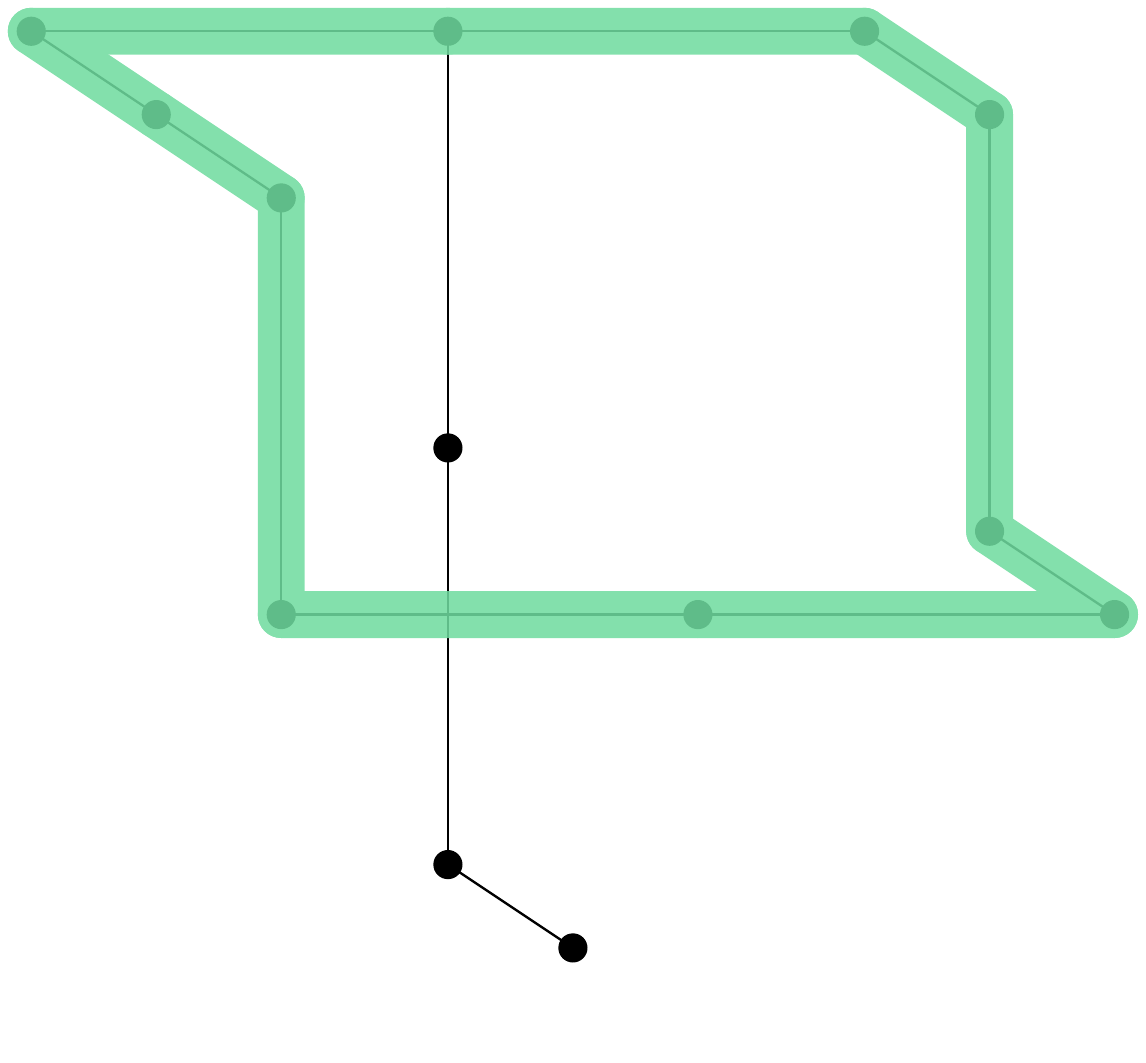}
            &
            \includegraphics[width=\linewidth]{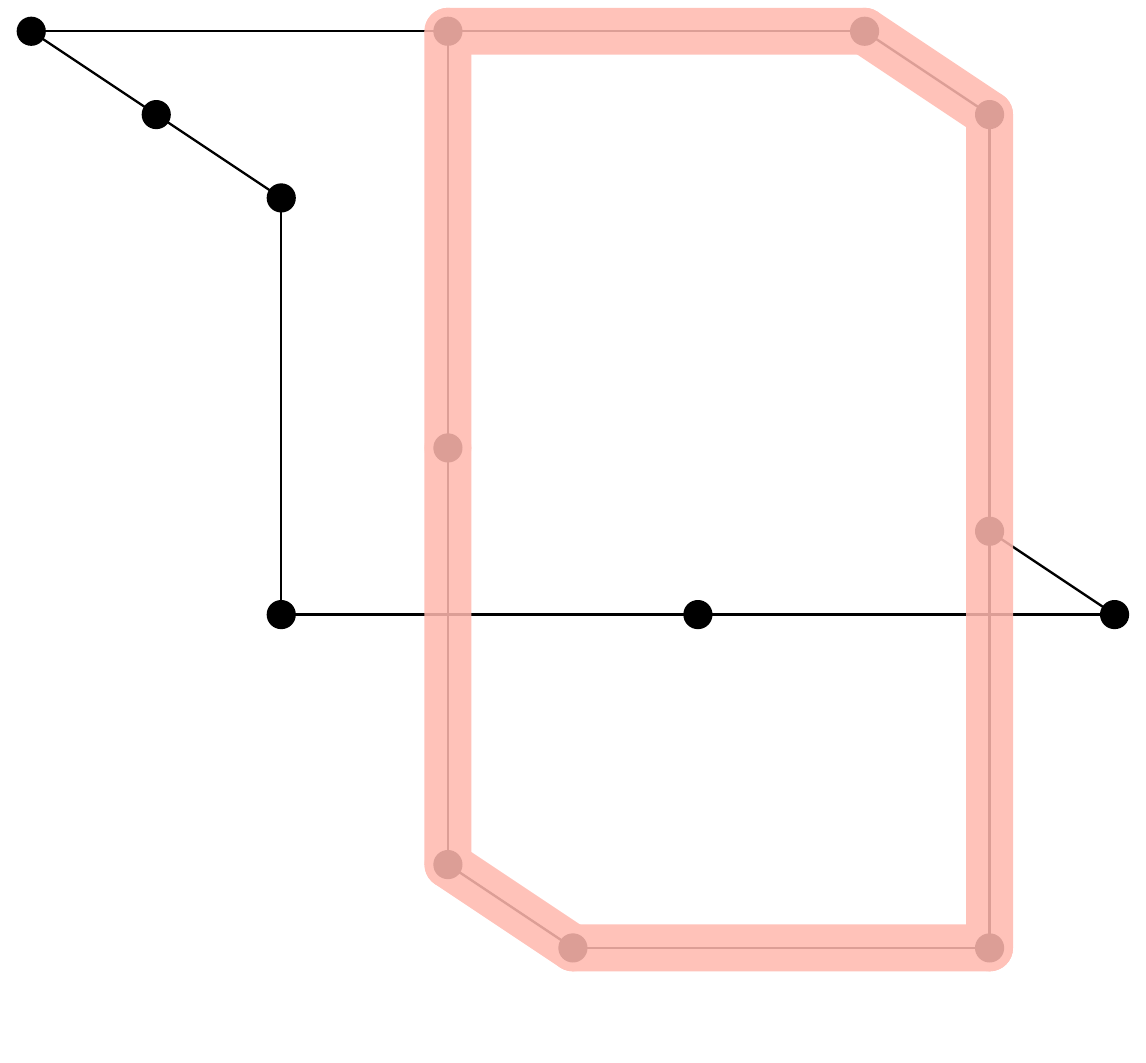}
            &
            \includegraphics[width=\linewidth]{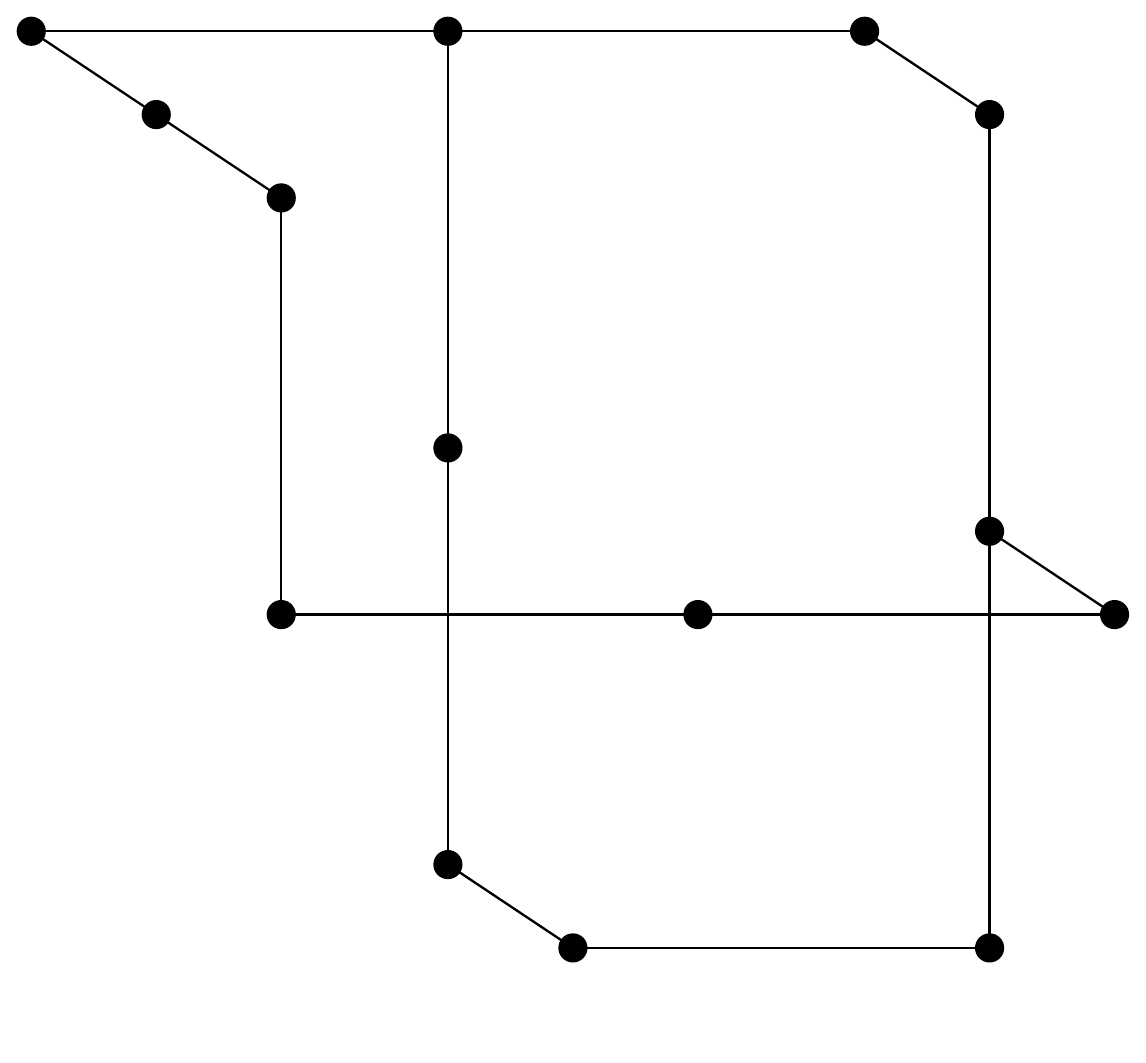}
            &
            \includegraphics[width=\linewidth]{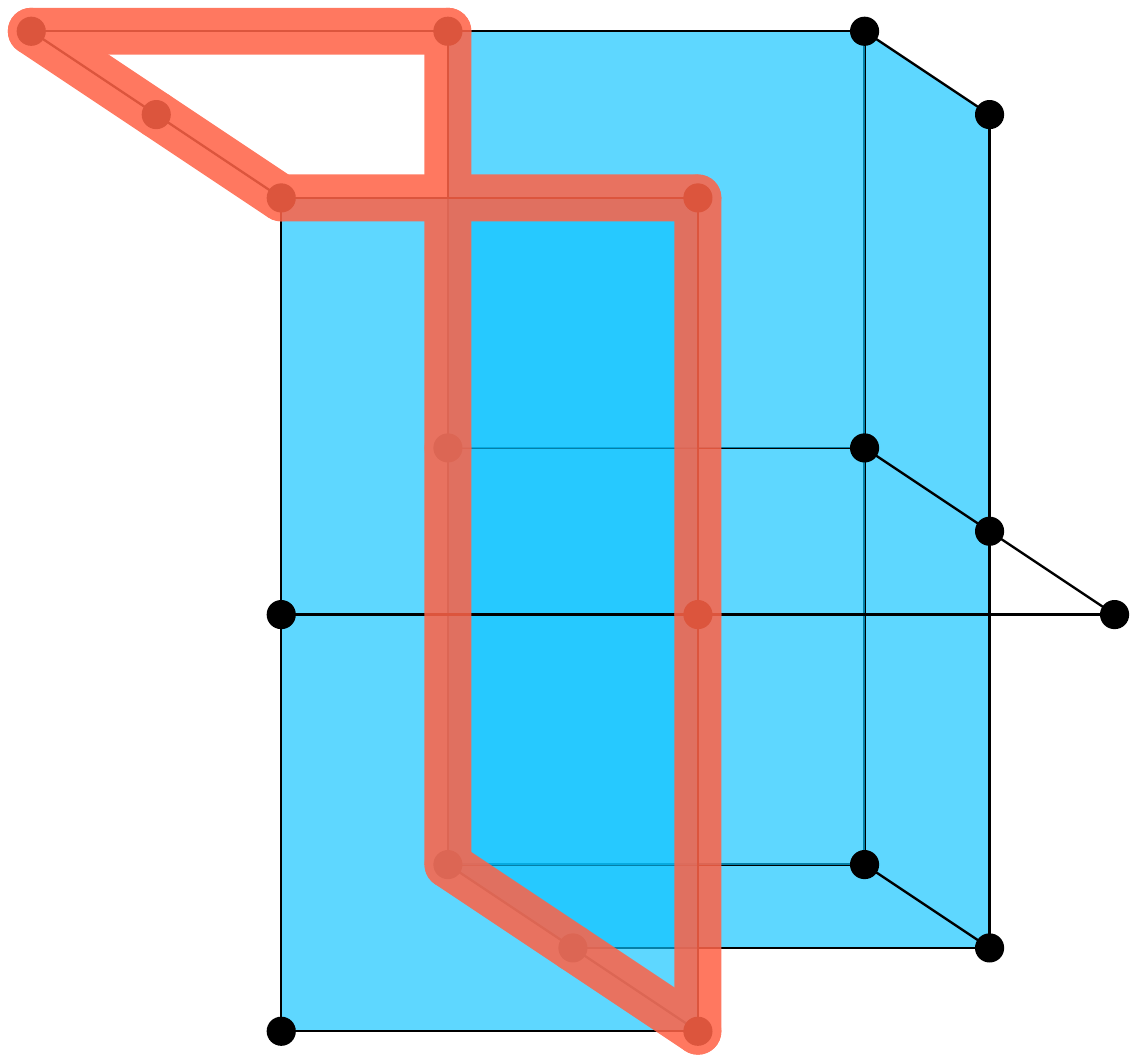}
            &
            \includegraphics[width=\linewidth]{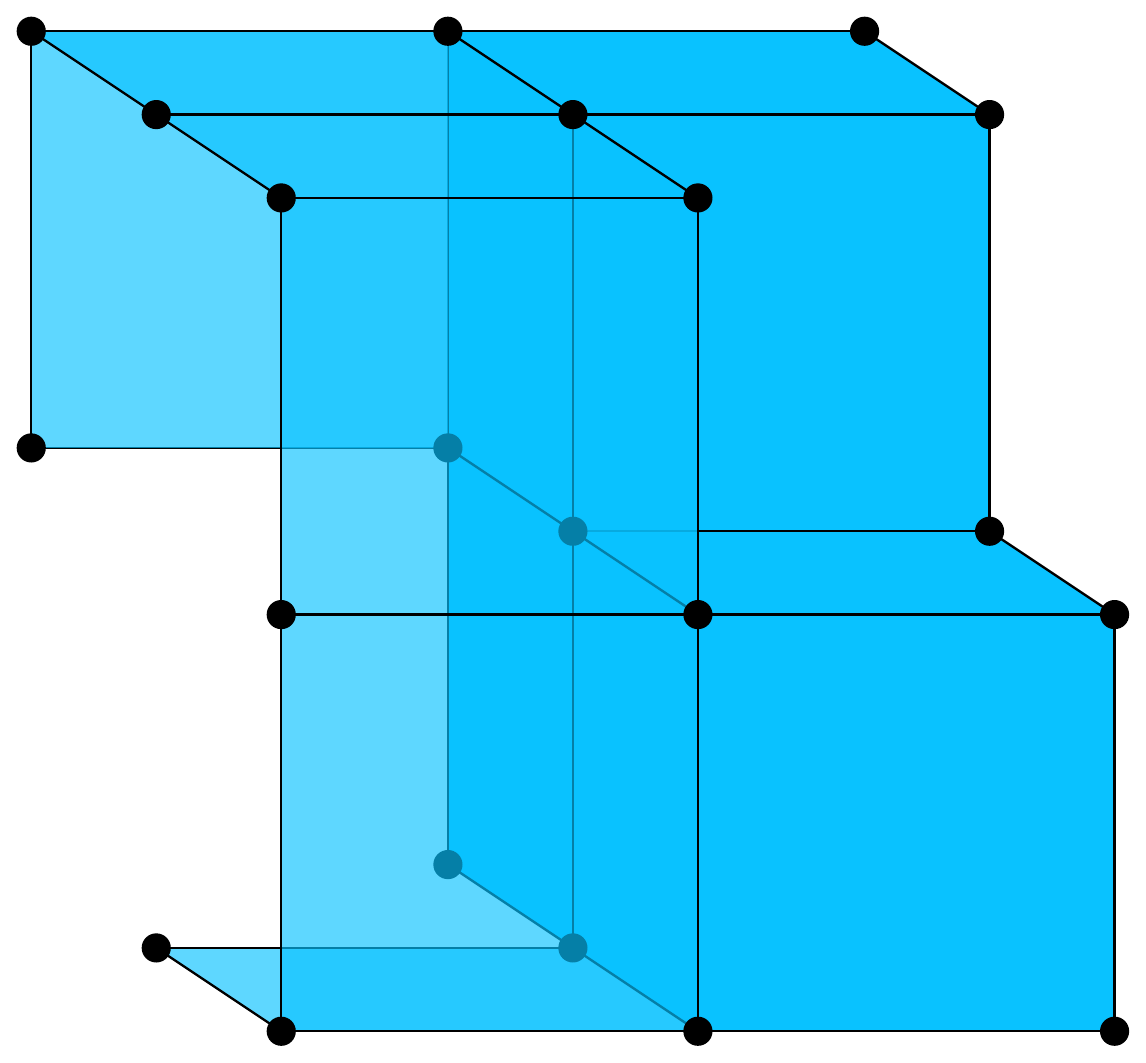}
            \\    
            
            $\mathbf{C}$ &
            \includegraphics[width=\linewidth]{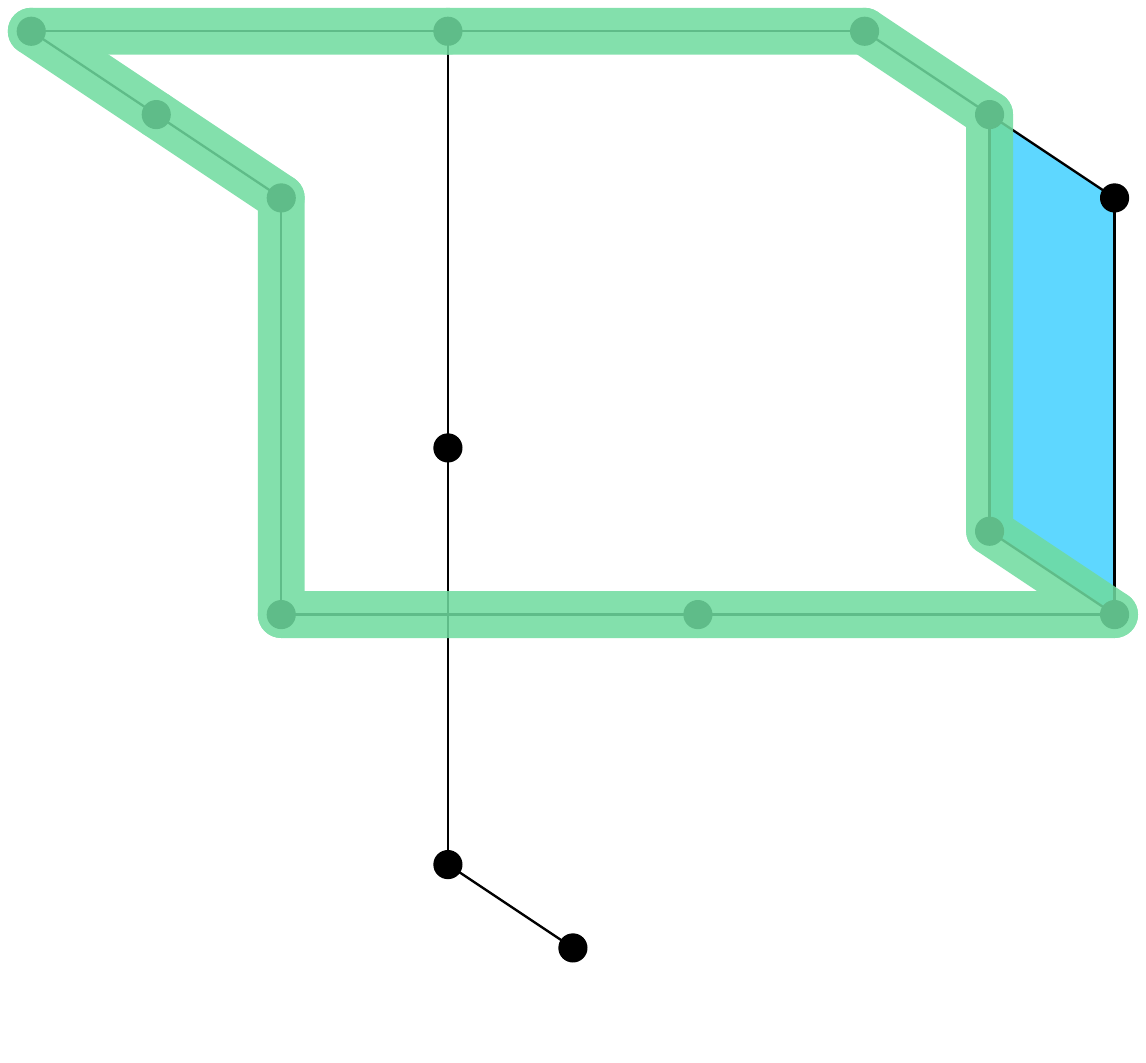}
            &
            \includegraphics[width=\linewidth]{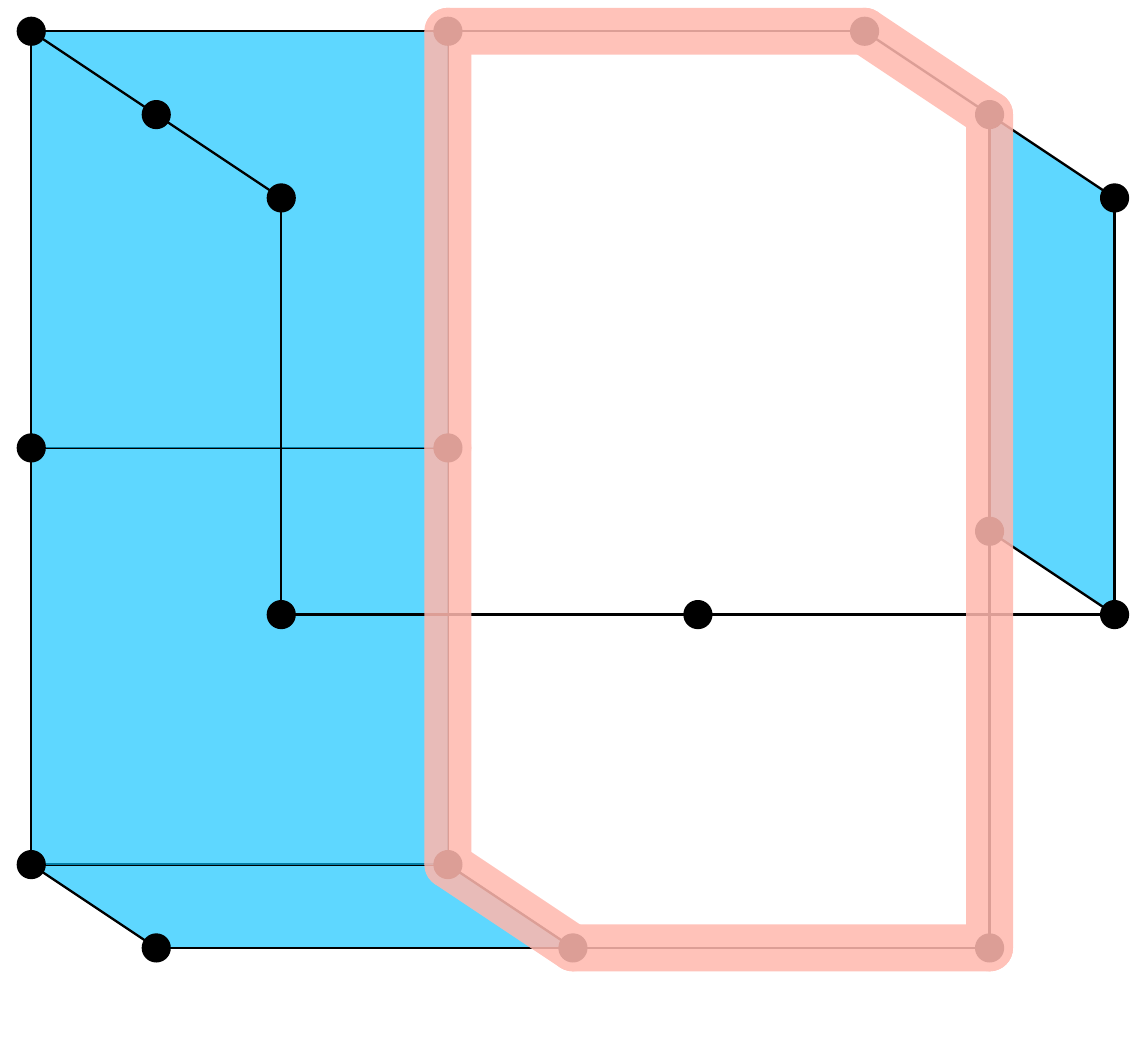}
            &
            \includegraphics[width=\linewidth]{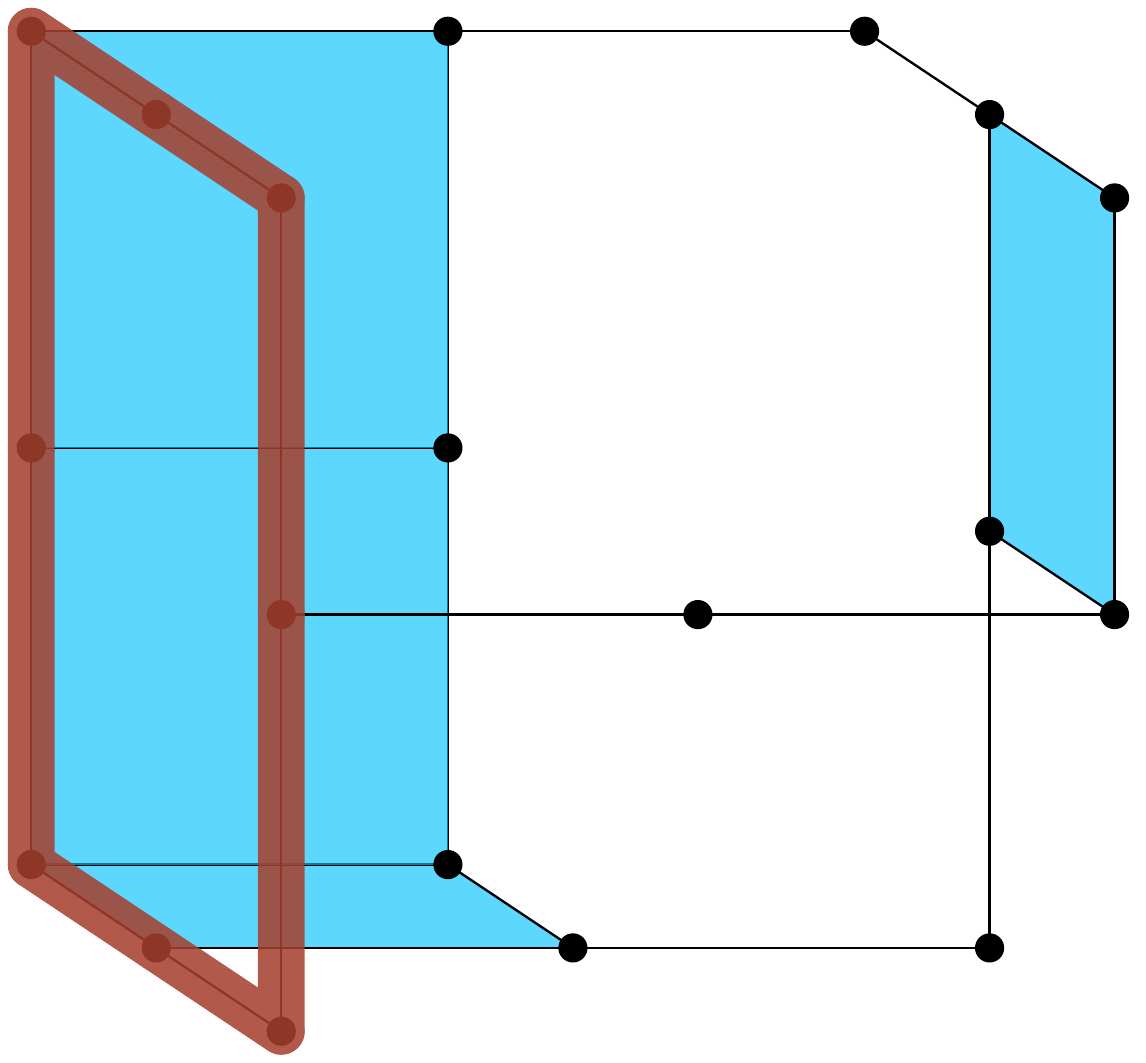}
            &
            \includegraphics[width=\linewidth]{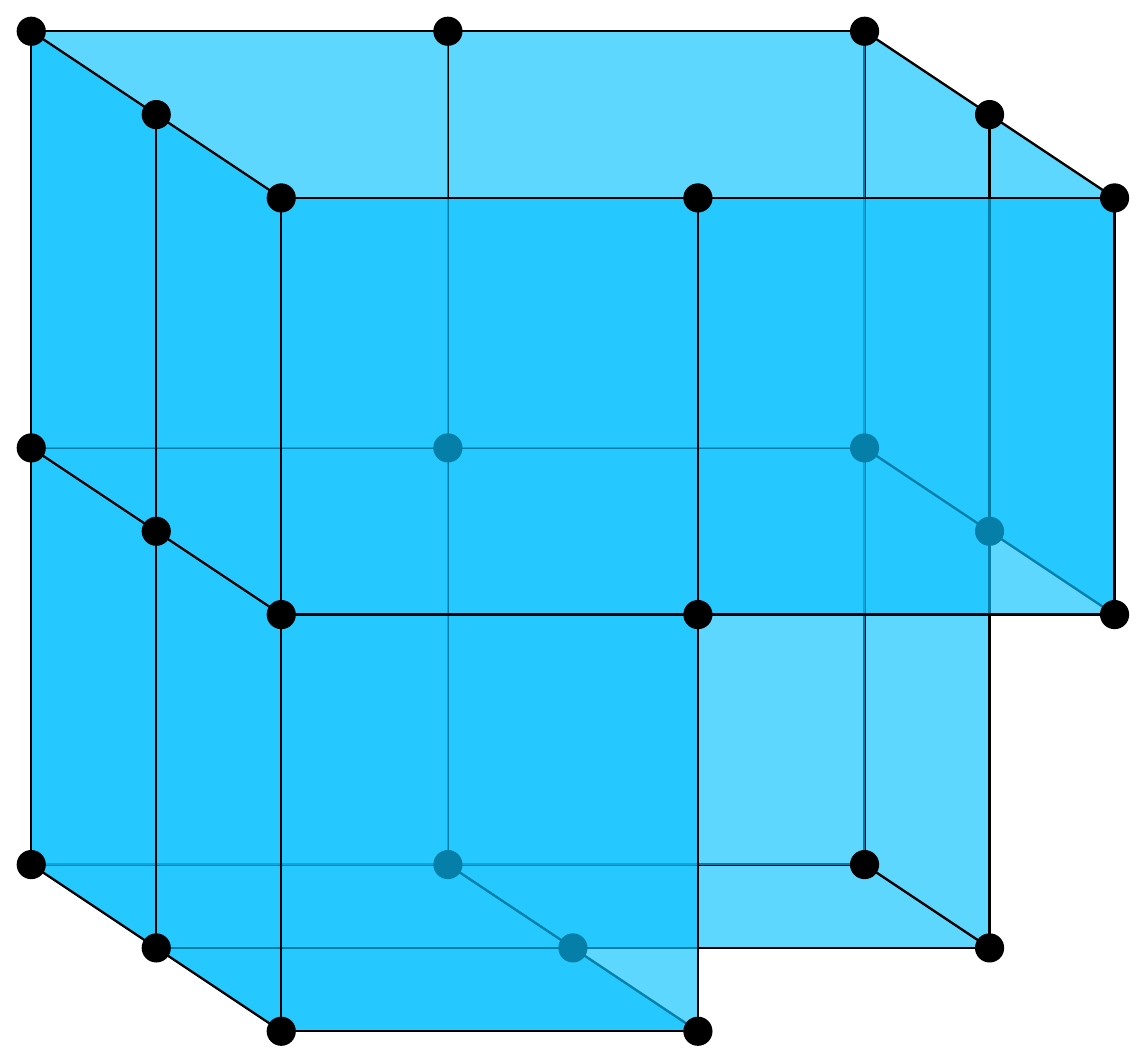}
            &
            \includegraphics[width=\linewidth]{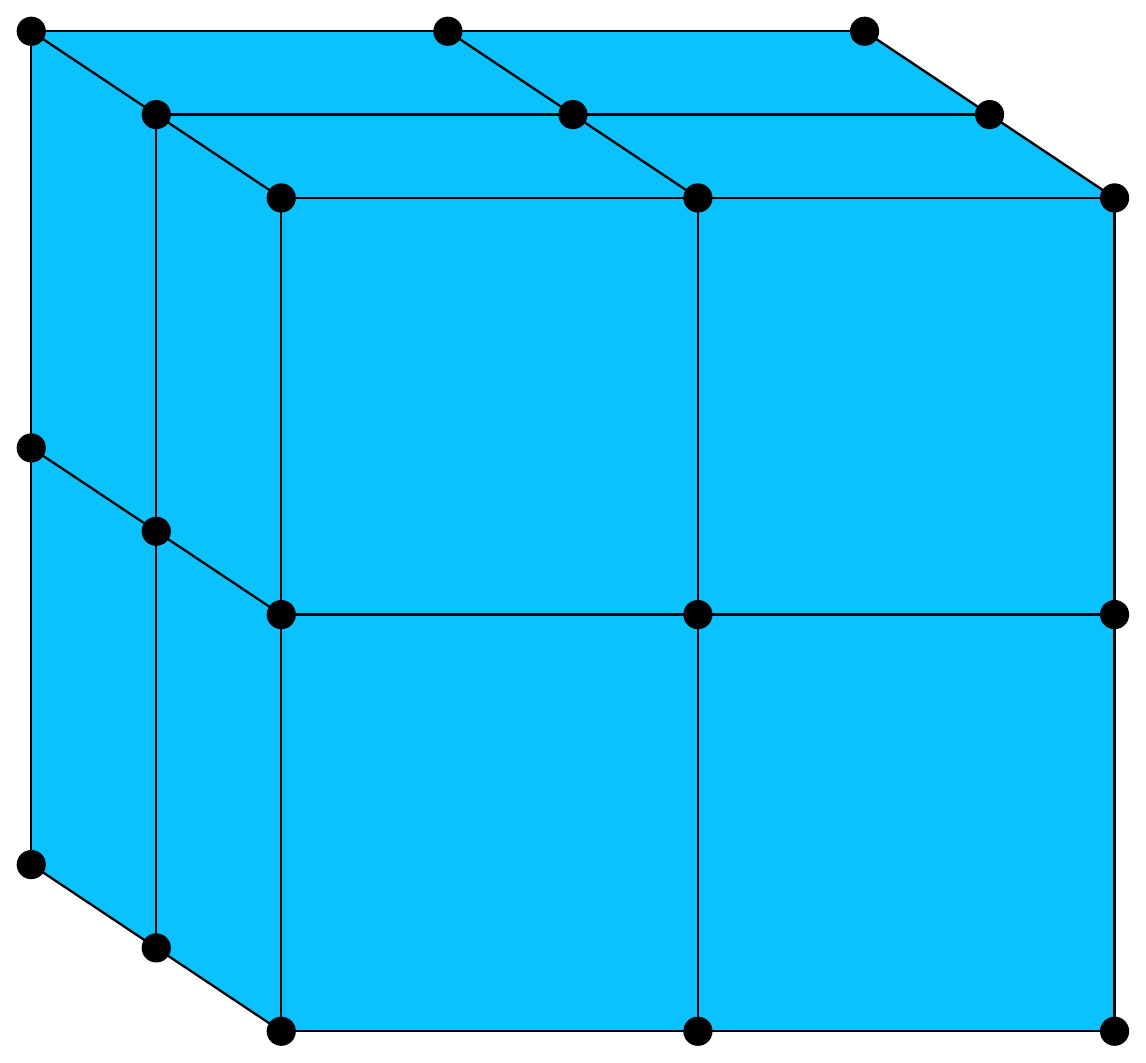}
            &
            \\
            
            $\mathbf{J}$ &
            \includegraphics[width=\linewidth]{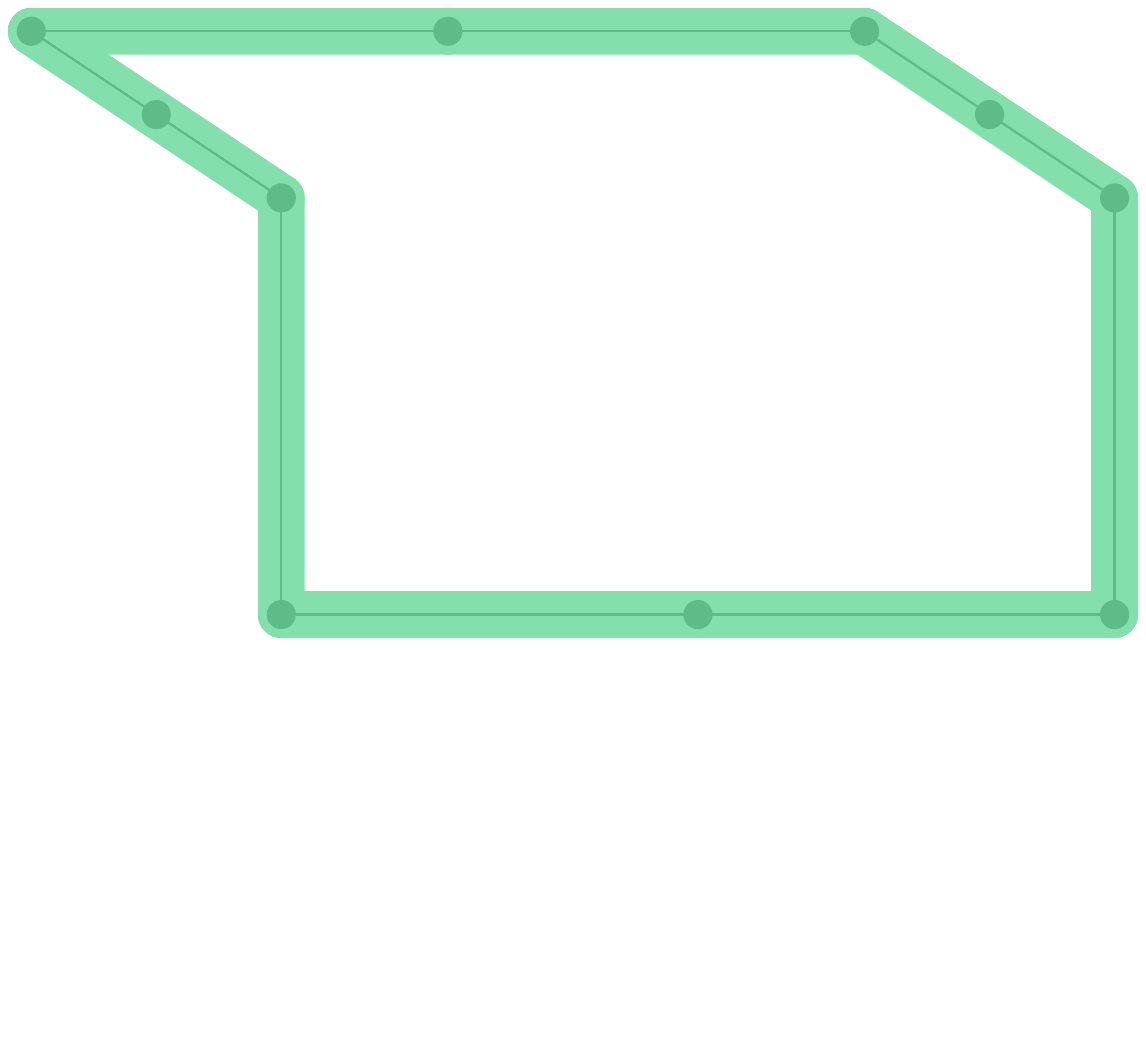}
            &
            \includegraphics[width=\linewidth]{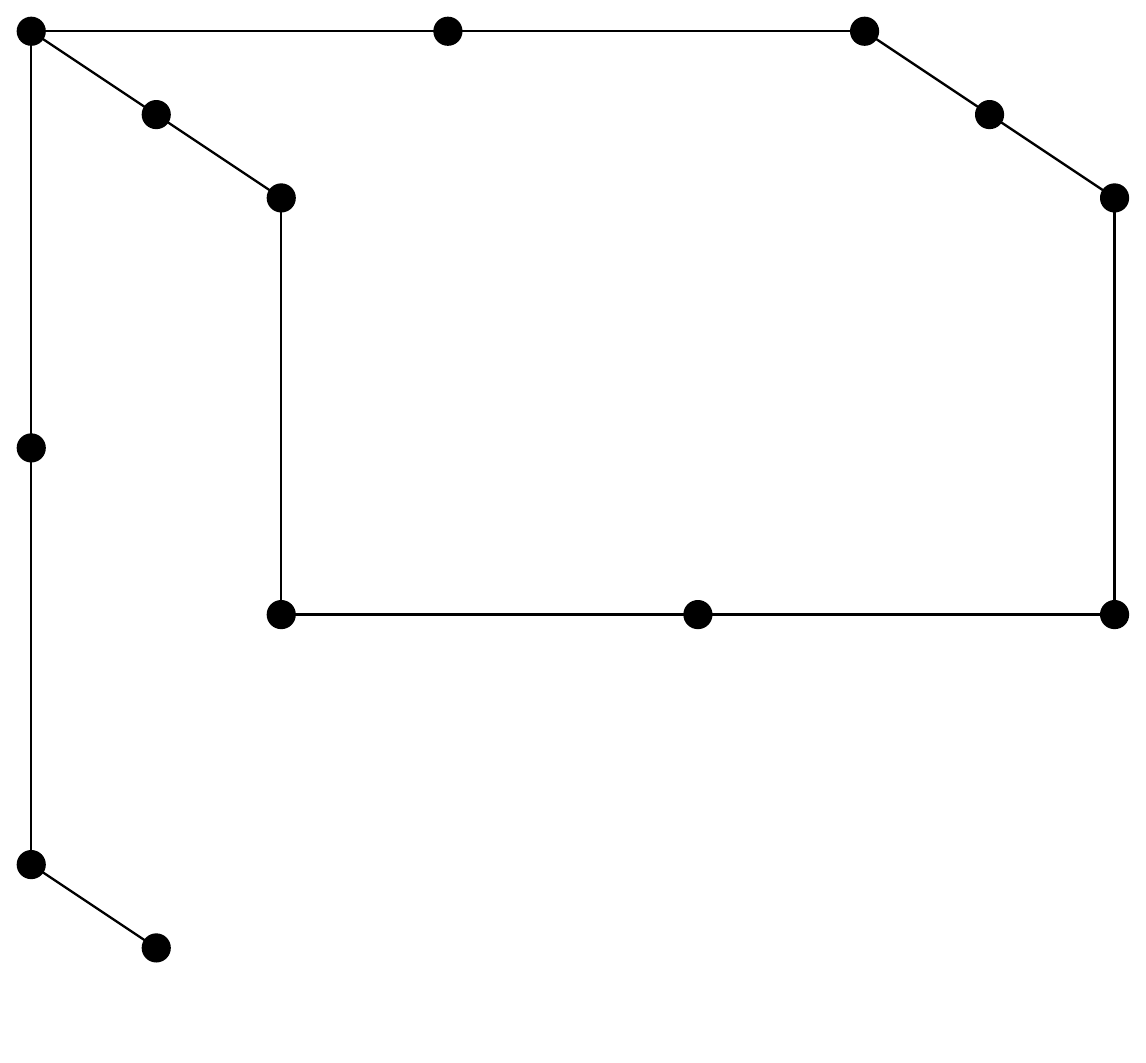}
            &
            \includegraphics[width=\linewidth]{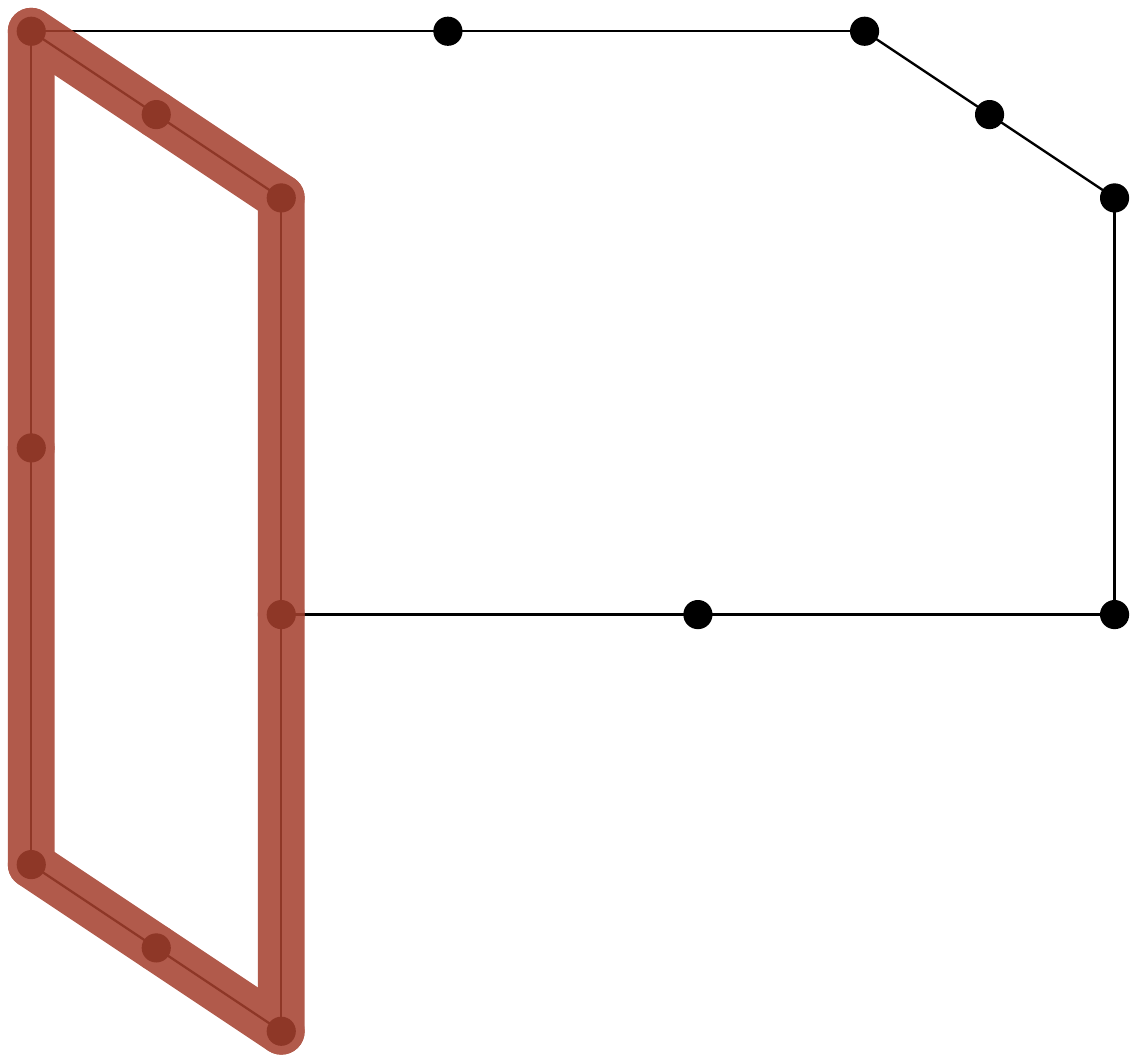}
            &
            \includegraphics[width=\linewidth]{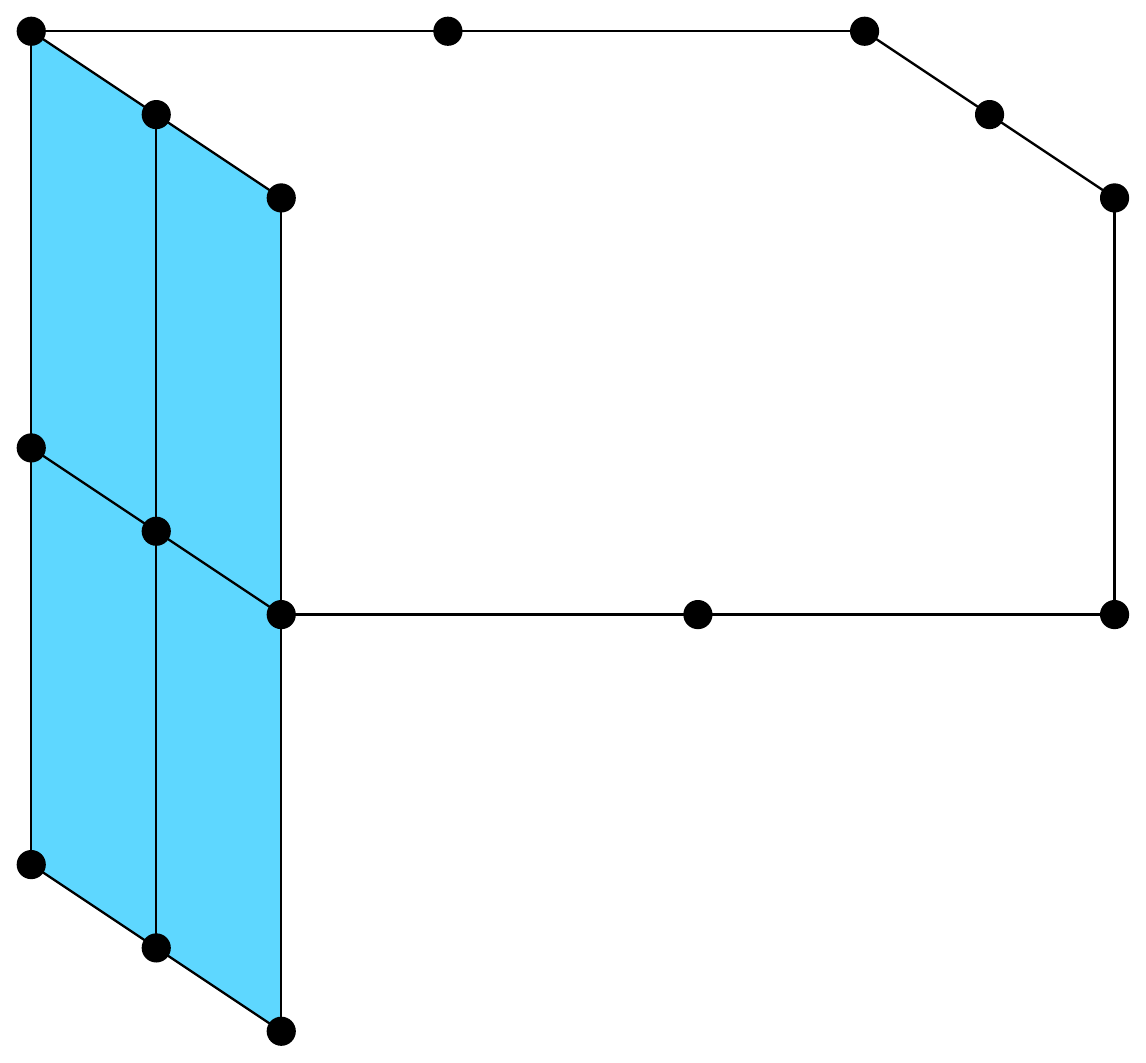}
            &
            \includegraphics[width=\linewidth]{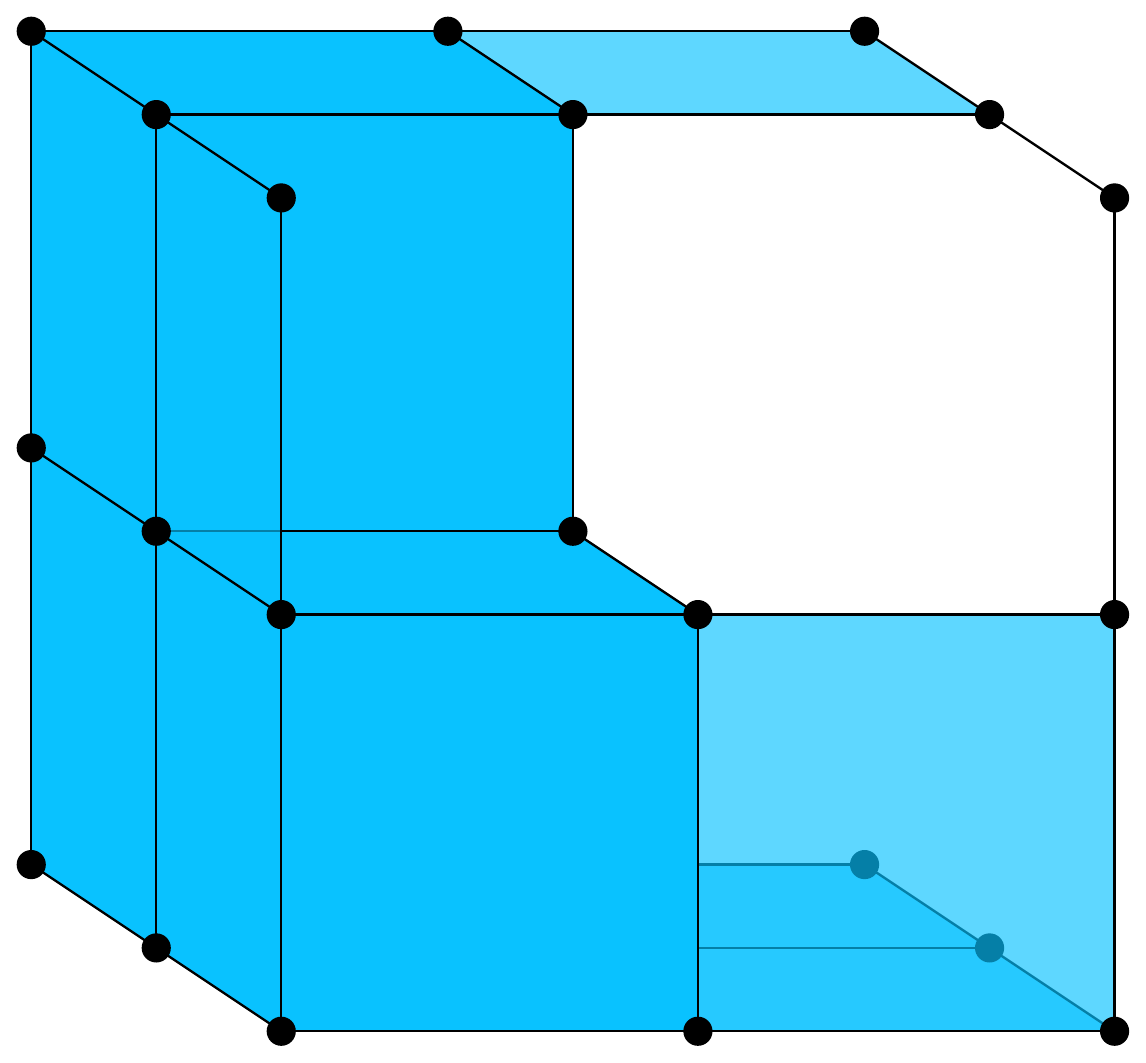}
            &
            \\

            &
            $\scriptstyle r=0.45$ & 
            $\scriptstyle r=0.55$ & 
            $\scriptstyle r=0.6$ & 
            $\scriptstyle r=0.8$ & 
            $\scriptstyle r=0.95$
        \end{tblr}
    }

    \subfloat[Input/image/comparison barcodes matched by the Betti matching
        % \label{subfig:image-barcode-betti-matching}
        ]{
        \centering
        \includegraphics[width=\exampleviifigurewidth]{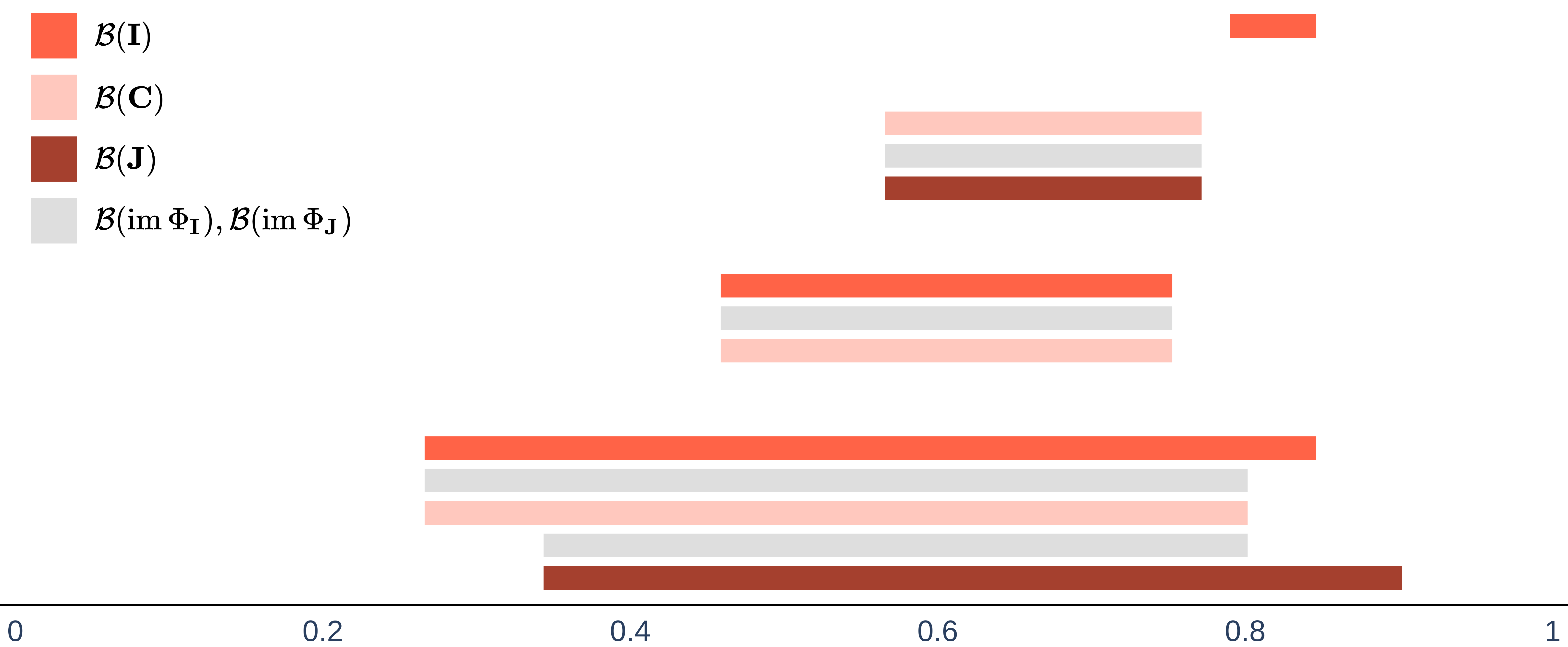}
    }

    \subfloat[Matching structure: a comparison barcode interval $[a,c)$ matches the input barcode intervals $[b_1, d_1)$, $[b_2, d_2)$ via image barcode intervals $[b_1, c)$, $[b_2, c)$.
        % \label{subfig:betti-matching-barcodes-structure}
    ]{
        \centering
        \includegraphics[width=\exampleviifigurewidth]{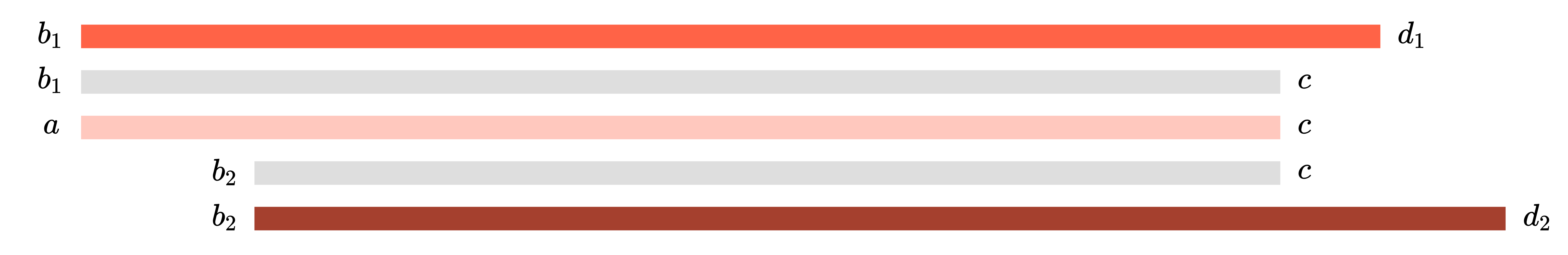}
    }

    \caption{Example of induced matchings and Betti matching in $H_1$. There is one matched pair of input bars, and there are three unmatched input bars, as shown in the barcode in (b) and visualized via representative cycles in the grid complex in (a). Two of the cycles that are unmatched by the Betti matching are still matched between their input barcode and the comparison barcode, but the comparison bar is not matched to the respective other input barcode, which is why the input bar is unmatched according to the Betti matching.
    The general relationship between the five bars involved in the Betti matching is shown in (c).}
    \label{fig:image-barcode-betti-matching}
\end{figure}

With this setup, we can invoke the induced matching theorem from \cite[Theorem 4.2]{bauer_induced_2015}, paraphrased in \cite[Theorem 2.1]{stucki_topologically_2023} for our situation.
The theorem states that there are unique injective maps $m_{\mathbf I}: \mathcal B(\im \Phi) \hookrightarrow \mathcal B({\mathbf I})$, $m_{\mathbf C}: \mathcal B(\im \Phi) \hookrightarrow \mathcal B({\mathbf C})$ which map each interval $[b, c)$ in the image barcode to an interval $[b, d) \in \mathcal B({\mathbf I})$, and to an interval $[a, c) \in \mathcal B({\mathbf C})$, such that $a \leq b < c \leq d$.
In other words: Every topological feature in $\im \Phi$ is uniquely mapped to one feature in $\mathbf I$ which is born at the same time, and one feature in $\mathbf C$ which dies at the same time.
The maps $m_{\mathbf I}, m_{\mathbf J}$ can be used to define a bijection between subsets of $\mathcal B({\mathbf I})$ and $\mathcal B({\mathbf C})$, the \textit{induced matching} $\sigma(\Phi) = m_{\mathbf C} \circ \big(m_{\mathbf I}{\big|}_{\im m_{\mathbf I}} \big)^{-1}$.

The notion of induced matchings allows us to define the Betti matching between images.
Given grayscale images $\mathbf I, \mathbf J \in \mathbb R^{N_1 \times N_2 \times N_3}$, define the comparison image $\mathbf C$ as their point-wise minimum.
The comparison image can be interpreted as a union, by the fact that the sublevel-sets satisfy $D(f_{\mathbf I})_s \cup D(f_{\mathbf J})_s \subseteq D(f_{\mathbf C})_s$.
Now we can construct induced matchings $\sigma_{\mathbf I}, \sigma_{\mathbf J}$ between $\mathbf{C}$ and $\mathbf I$ or $\mathbf{J}$, respectively.
The two matchings taken together induce a matching $\mu_{\mathbf I, \mathbf J}$, called the \textit{Betti matching} between $\mathbf{I}$ and $\mathbf J$ \cite{stucki_topologically_2023}. For each barcode entry $[a, c) \in \mathcal B(\mathbf C)$ that is matched according to both $\sigma_{\mathbf I}$ and $\sigma_{\mathbf J}$, the Betti matching $\mu_{\mathbf I, \mathbf J}$ matches $\sigma_{\mathbf I}\big([a, c)\big) = [b_1, d_1) \in \mathcal B(\mathbf{I})$ with $\sigma_{\mathbf J}\big([a, c)\big) = [b_2, d_2)\in\mathcal B(\mathbf{J})$.

In \cite{stucki_topologically_2023}, an extension to the induced matchings is proposed. 
The image barcode only captures topological features that are alive "at the same time", i.e., during an overlapping range of intensity values, in both images. 
Depending on the application, the goal can be to also match other topological features that do spatially correspond, but not live at the same time.
The algorithm for computing image barcodes (which we will describe in \cref{sec:computing-image-barcode-betti-matching}) computes, as a by-product, additional image persistence pairs $(c_i, c_j)$ with $f_{\mathbf I}(c_i) \geq f_{\mathbf C}(c_j)$, which according to \cite[Theorem 3.6]{bauer_efficient_2022}, are filtered out since they are not a part of the image barcode.
\cite{stucki_topologically_2023} propose to use those \textit{reverse pairs} for the matching, which allows to match some spatially corresponding features that would otherwise be unmatched.
To distinguish between the two versions, we shall call the versions including reverse pairs the \textit{extended image barcode}, the \textit{extended induced matchings} and the \textit{extended Betti matching}, as opposed to the \textit{Betti matching without reverse pairs}. 
From here on, when we talk about the Betti matching, we mean the extended Betti matching, and we also use this version in our implementation and our experiments.

\section{Algorithm and Implementation}

Being able to efficiently compute the barcodes and the Betti matching is essential, as we will rely on running these computations thousands of times when training a segmentation model. 
While \cite{stucki_topologically_2023} were able to train 2D segmentation models using a Python implementation for the persistent homology computations, the same turns out to be infeasible for 3D images.
In this work, a new highly optimized C++ implementation was used to compute the Betti matching on 3D images.
It was developed primarily by Nico Stucki, and extended, further optimized and equipped with a Python interface by Vincent Bürgin.
\emph{Betti-matching-3D} builds on the implementation of \emph{Cubical Ripser} \cite{kaji_cubical_2020}, with which it shares the basic data structures. 
Cubical Ripser in turn builds on the algorithm used by \emph{Ripser} \cite{bauer_ripser_2021}, adjusted for cubical complexes. 
In comparison to Cubical Ripser, Betti-matching-3D additionally implements the image barcode computation, the Betti matching, and several additional optimizations.

\subsection{Computing the Barcode}
\label{sec:computing-the-barcode}
The basic algorithm for computing persistent homology consists of reducing the \textit{filtered boundary matrix} using a variant of Gaussian elimination. 
We will first introduce the basic algorithm following \cite{edelsbrunner_short_2014} before describing the optimized algorithm used in Betti-matching-3D.

The \textit{filtered boundary matrix} $D$ is an $N \times N$ matrix over $\mathbb F_2$, where $N$ is the total number of cubes in the complex. 
For each cube $c_i$, indexed by the index of appearance $i$ of the cube according to the cube-wise filtration, the $i$-th column encodes the boundary of $c_i$: 
If $c_j$ is a facet of $c_i$, the entry $D_{j,i}$ will be 1, otherwise 0. 
The basic reduction algorithm reduces the matrix $D$ using only left-to-right column operations to obtain a \emph{reduced matrix} $R$. 
Define the \textit{pivot} of a column $j$ in a matrix $M$ as the largest row index $k$ where column $j$ has a non-zero entry, or $\pivot_M(j) = \bot$ if the column is zero.

\begin{lstlisting}[language=python,style=PseudoPython,keepspaces=True,escapeinside={(*}{*)}]
R = D
# iterate over all columns
for i in 1..N:
  j = 0
  # iterate over all columns left of i:
  while j < i and pivot(R, i) != (*$\bot$*):  
    if pivot(R, j) == pivot(R, i):
      # reduce i-th column by j-th column
      R[:, i] = (R[:, i] + R[:, j]) mod 2
      # restart at left-most column
      j = -1
    j += 1
\end{lstlisting}%

%
% 1_basicReductionAlgorithm.tex listing
%

The pivots of the reduced matrix $R$ encode the information needed to reconstruct the barcode: 
If $\pivot_R(i)=j$, then the row's cube $c_i$ creates a class that the column's cube $c_j$ destroys. 
In other words, $(c_i, c_j)$ is a persistence pair and we can find all persistence pairs from the pivots of the reduced matrix \cite{edelsbrunner_short_2014}. 

The basic reduction algorithm is not particularly efficient if implemented directly, but can be improved by exploiting the additional structure of the matrix. 
There has been a significant amount of research on optimizations since the algorithm was first introduced, reducing the run time by several orders of magnitude (see \cite{bauer_ripser_2021}). 
Betti-matching-3D uses several of these optimizations, and the remainder of this section describes them in detail.
There are four major algorithmic changes:
\begin{enumerate}
    \item \textbf{Separation of reductions by dimension:} Instead of working on the full boundary matrix $D$, we separate the subproblems of computing $H_d$ persistence for $d \in \{0, 1, 2\}$.
    We first compute the $H_2$ persistence, then the $H_1$ and finally the $H_0$ persistence. 
    Going from top to bottom enables the \textit{clearing optimization} \cite{chen_persistent_2011}, which allows to skip certain columns in dimension $d-1$ using information gained during the computation of dimension $d$.

    \item \textbf{Implicit matrix reduction:} The boundary matrix $D$ is not explicitly realized in memory, and instead its columns are computed on the fly and in a sparse representation. 
    Reduced columns are stored in a cache for the use in future reductions.

    \item \textbf{Union-find for $\mathbf{H_2}$ and $\mathbf{H_0}$:}
    The structure of $H_2$ and $H_0$ persistence lends itself to a more efficient algorithm than the matrix reduction: 
    $H_0$ captures connected components, and the persistence pairs can be computed by a union-find algorithm. 
    Furthermore, $H_2$ persistence pairs can be computed in a similar way by exploiting \textit{Alexander duality}.
    Therefore the expensive matrix reduction algorithm needs to be performed only for computing the $H_1$ persistence pairs.\footnote{A consequence of this is that computing Betti matching on 2-dimensional images, also available in Betti-matching-3D, can fully rely on union-find and is hence highly efficient.}

    \item \textbf{Emergent pairs:} The $H_1$ reduction algorithm makes use of the \textit{emergent pair} optimization. 
    This refers to certain zero persistence pairs which can be detected from the boundary matrix without reduction, which means the reduction can be skipped in this case \cite{bauer_ripser_2021} (see also \cite{zhang_hypha_2019}).
\end{enumerate}

To make matters concrete, let us build up our algorithm step by step. We use a Python-inspired pseudocode syntax to describe the core algorithm concisely and will mention C++-specific implementation details at the end.

The input image is converted to a list of cubes, which are represented by a data structure that holds the cube's birth intensity, its coordinates, and its \textit{type}.
We assume for all algorithms that the refined filtration order uses lexicographic tie-breaking on the $x$, then $y$, then $z$ coordinate, then the type.

The persistence pair algorithm computes in sequence the $H_2$ persistence pairs, the $H_1$ persistence pairs and the $H_0$ persistence pairs. The $H_2$ and $H_1$ computations furthermore pass down the list of 2-cubes and 1-cubes (i.e. columns) to be reduced during the $H_1$ and $H_0$ computation, respectively.
This is necessary to enable the clearing optimization, which means that the higher dimension algorithms can detect that certain columns will be unnecessary to reduce in the next-lower dimension. 
The higher dimension $d$ algorithm enumerates the $(d-1)$ cubes to be reduced, filters out the unnecessary ones, and passes it down to the dimension $(d-1)$ algorithm. 
The high-level structure is as follows:

%
% 2_computePersistencePairs_highLevel.tex listing
%
\begin{lstlisting}[language=python,style=PseudoPython,keepspaces=True]
function computePersistencePairs(gridComplex):
  (persistencePairs2, columnsToReduce2): ([(Cube2, Cube3)], [Cube2]) = computePersistencePairsDim2(gridComplex)
  (persistencePairs1, columnsToReduce1): ([(Cube1, Cube2)], [Cube1]) = computePersistencePairsDim1(gridComplex, columnsToReduce2)
  persistencePairs0: [(Cube0, Cube1)] =                computePersistencePairsDim0(gridComplex, columnsToReduce1)
  return (persistencePairs2, persistencePairs1, persistencePairs0)
\end{lstlisting}

\subsubsection{Dimension $1$}
\label{sec:barcode-dimension-1}

Let us first focus on the modified matrix reduction algorithm to compute $H_1$ persistence.
The first part to implicit matrix reduction is a \textit{boundary enumerator} which takes a $2$-cube and enumerates its four boundary 1-cubes on the fly. 
We will assume to we have a function \code{enumerateBoundary(cube, gridComplex)} that, given a \code{Cube}, enumerates its faces (using a hard-coded case distinction based on the input cube's \code{type}), and has access to the image array to compute birth values. 
The enumeration should be lazily evaluated (i.e. coroutine-/generator-style), since in some cases we want to be able to stop the enumeration early.\footnote{In the C++ implementation, we instead repeatedly call \code{hasNextFace()} on a \code{BoundaryEnumerator} object, but the generator-style formulation makes the pseudocode more concise.}

The second part is an efficient representation of the column that we are currently reducing, the \textit{working boundary}. 
It needs to support adding cubes and obtaining the pivot, i.e., the largest contained cube with respect to filtration order.
Crucially, the computation must be modulo 2: adding the same cube twice cancels it out. 
Furthermore, the working boundary can be converted to a list, with modulo 2 deduplication performed, which we use for caching.
We realize the working boundary with a data structure \code{CubeXorQueue}, based on a priority queue, that lazily cancels out duplicate cubes. 
We will describe the implementation after describing the main algorithm.

Using these ingredients, the $H_1$ reduction algorithm proceeds as follows: 
An outer loop iterates over the list of columns to reduce, which are sorted in filtration order, and in each iteration processes a $2$-cube $c_i$. 
The working boundary represents the partially reduced column and is initialized with the boundary of $c_i$. 
Then an inner loop repeatedly reduces the working boundary: 
\begin{enumerate}
    \item Compute the working boundary's pivot.
    \item Use a lookup table \code{columnIndexByPivot} to look up the index $j < i$ such that the reduced column corresponding to $c_j$ has the same pivot as the working boundary.
    \item Reduce the working boundary with the cached reduced column of $c_j$ or, if not cached, with the boundary of $c_j$ computed on the fly.
\end{enumerate}
The inner loop terminates if no such previous column exists or the working boundary is empty:
in both cases, the working boundary is fully reduced.
Finally, if not empty, the working boundary is saved to cache and the index $i$ is added to the lookup map of columns by pivot.
Then a non-zero persistence pair is recorded if the pivot's birth does not equal the birth of $c_i$.
The algorithm can be slightly shortened by not handling the initial working boundary as special case: instead initialize the working boundary as empty, initialize $j$ as $j = i$ and make the reduction step the first (not last) step in the inner loop.
This yields the following algorithm:

%
% 3_computePersistencePairsDim1.tex listing
%
\begin{lstlisting}[language=python,keepspaces=true,style=PseudoPython]
function computePersistencePairsDim1(gridComplex, columnsToReduce: [Cube2])    -> [(Cube1, Cube2)]:
  columnIndexByPivot: Dict[Cube1, int] = {}
  cache: Dict[Cube2, [Cube1]] = {}
  persistencePairs: [(Cube1, Cube2)] = []

  for i in 1..columnsToReduce.length:
    workingBoundary: CubeXorQueue = new CubeXorQueue()
    j = i
    pivot: Cube1 = (*$\bot$*)
    do:
      reduceWorkingBoundaryBy(workingBoundary, columnsToReduce[j], cache, gridComplex)
      pivot = workingBoundary.getPivot()
      if pivot != (*$\bot$*):
        j = columnIndexByPivot[pivot] # returns (*$\bot$*) if not present
    while (pivot != (*$\bot$*) and j != (*$\bot$*))

    if pivot != (*$\bot$*): # i.e. the columns was not reduced to zero:
      columnIndexByPivot[pivot] = i
      cache[columnsToReduce[i]] = workingBoundary.toList()
      # record persistence pair if it is a non-zero pair
      if pivot.birth != columnsToReduce[i].birth:
        persistencePairs.push((pivot, columnsToReduce[i]))

  columnsToReduce1: [Cube1] = enumerateDim0Columns(columnIndexByPivot, gridComplex)
  return (persistencePairs, columnsToReduce1)
\end{lstlisting}%

The sub-procedure that reduces the working boundary by a column is given in the listing \ref{alg:reduceWorkingBoundaryBy} in the appendix.

%
% 4_reduceWorkingBoundaryBy.tex listing
%

The \code{CubeXorQueue} is implemented internally with a priority queue that compares based on the filtration order and puts highest priority on the youngest cube, see listing \ref{alg:cubeXorQueue} in the appendix.
It deduplicates lazily when a pivot is requested in the \code{getPivot()} function: 
It retrieves the highest-priority elements $(c_1, c_2)$ and removes both if $c_1 = c_2$, or keeps both and returns $c_1$ if $c_1 \neq c_2$. 
It repeats the process until two unequal elements are found or the underlying queue is empty. 
The \code{toList()} operation employs the same technique and returns a modulo-2 reduced list.
The lazy deduplication technique was already employed in Ripser as well as in earlier implementations (cf. \cite{bauer_ripser_2021}).

%
% 5_CubeXorQueue.tex listing
%

\paragraph{Emergent Pairs}
\label{paragraph:emergent-pairs}
The dimension $1$ matrix reduction can be further optimized by taking \textit{emergent pairs} into account. 
A $1$-cube $\sigma$ and a $2$-cube $\tau$ form an emergent facet pair $(\sigma, \tau$), which is a zero persistence pair, iff \cite[Proposition 3.12]{bauer_ripser_2021}:
\begin{itemize}
    \item
    $\sigma$ is the lexicographically maximal facet of $\tau$ among those with the same birth as $\tau$: In other words, the youngest facet of $\tau$.
    
    \item
    $\sigma$ does not form a (non-zero or zero) persistence pair $(\sigma, \rho$) with any 2-cube $\rho$ that is older than $\tau$.
\end{itemize}%
This proposition allows our algorithm to skip the reduction loop for some columns. 
In essence, this optimization exploits that a pivot can be found faster for the initial single-cube boundary than for a general working boundary. 
Assume that in addition to \code{enumerateBoundary}, we have a function \code{enumerateBoundaryReverse(cube, gridComplex)} that lazily enumerates the boundary of a cube in \textit{reverse} lexicographic order. The necessary changes for the emergent pair optimization are shown in listing \ref{alg:emergentPairsOptimization} in the appendix.

%
% 6_emergentPairs.tex listing
%

\paragraph{Clearing} \label{paragraph:clearing-dim0-for-dim1}
Before returning the persistence pairs, \code{computePersistencePairsDim1()} calls \code{enumerateDim0Columns()} to enumerate the columns for dimension $0$ and perform the clearing optimization. 
The enumeration function uses a nested loop over the three spatial axes and three possible types, creates the non-cleared cubes and sorts them.

The clearing optimization \cite{chen_persistent_2011} is based on the following observation: 
If the reduction of the $i$-th column in dimension $d$ yields a pivot in the $j$-th row, this reveals not only that the cube $d$-cube $c_i$ kills a homology class, but also that the $(d-1)$-cube $c_j$ creates it.
Furthermore, it is known that the column corresponding to a cube that creates a homology class reduces to zero in the reduction algorithm. 
Hence, the column corresponding to $c_j$ can be skipped in the dimension 0 reduction algorithm.
In our case, this means that during the nested loop that enumerates the $1$-cubes for the dimension $0$ computation, we filter out those $1$-cubes which have been identified as pivots by consulting the \code{columnIndexByPivot} map (see listing \ref{alg:enumerateDim0Columns} in the appendix).

%
% enumerateDim0Columns.tex listing
%

\subsubsection{Dimension $0$}
\label{sec:barcode-dim0}
Computing the $H_0$ persistence can be reduced to a union-find algorithm since $H_0$ homology classes describe connected components. 
We set up a union-find data structure with a node for each vertex in the grid complex, and add the edges in filtration order. 
At each step we merge the connected components that the edge's endpoint vertices lie in. 
We are not interested in the end result, but rather in the process of getting there:
whenever two connected components are merged, this corresponds to the death of the younger one. 
For each component, the union-find data structure tracks a representative vertex, which on merge we set as the older one of the two representatives of the merged components. 
Assume adding an edge merges two previously unconnected components: 
In this case the younger component dies. That is, if two components represented by the cubes $c_h \neq c_i$ are merged by an edge $c_j$, and wlog. $c_h$ is older than $c_i$ in the refined order, we found a persistence pair $(c_i, c_j)$. 
If furthermore $c_i$ is strictly older than $c_j$ in the non-refined order, the persistence pair is a non-zero pair which we record.
Otherwise we found a zero persistence pair which we ignore. \cite{kaji_cubical_2020}
We implement the algorithm using a union-find forest data structure with path compression (cf. \cite[Section 21.3]{cormen_introduction_2009}).

\begin{lstlisting}[language=python,style=PseudoPython]
function computePersistencePairsDim0(gridComplex, columnsToReduce: [Cube1])    -> [(Cube0, Cube1)]:
  # initialize a UnionFind structure with all vertices in gridComplex
  UnionFind uf = new UnionFind(gridComplex)
  persistencePairs: [(Cube0, Cube1)] = []
  
  for edge in columnsToReduce:
    endpoint1, endpoint2 = findEndpoints(edge)
    representative1, representative2 = (uf.find(endpoint1), uf.find(endpoint2))
    if representative1 != representative2:
      # uf.merge() sets the older as new representative and returns the younger
      youngerRepresentative: Cube0 = uf.merge(representative1, representative2)
      if youngerRepresentative.birth < edge.birth:
        persistencePairs.append((youngerRepresentative, edge))
  return persistencePairs
\end{lstlisting}%
%
%
% 8_computePersistencePairsDim0.tex listing
%

\subsubsection{Dimension $2$}
\label{sec:barcode-dim2}
Perhaps surprisingly, the computation in dimension $2$ can also be reduced to a union-find procedure, by exploiting Alexander duality. 
This idea is discussed in detail in \cite{bleile_persistent_2022,garin_duality_2020}.
The intuition is that $2$-dimensional homology classes represent cavities, and cavities are nothing else but connected components in the inverted input image. 
More precisely, we construct a dual complex where $3$-cubes become vertices, $2$-cubes become edges, and the filtration order is reversed.

Consider a dual grid complex where the $3$-cubes become \emph{dual vertices}, and the $2$-cubes become \emph{dual edges}. 
The two facets of a dual edge are the dual vertices corresponding to the two co-facets of the corresponding $2$-cube. 
The $2$-cubes at the border of the complex only have one adjacent $3$-cube. 
Their dual edges are assigned a token $3$-cube $c_*$ as second endpoint that has birth intensity $\infty$ and represents the outside of the complex.
With this setup, we perform the union-find algorithm on the dual complex in an analogous way to the computation of dimension $0$. 
The required changes are:
\begin{itemize}
    \item Traverse the dual edges in \textit{reverse} filtration order of their corresponding 2-cubes.

    \item On merging two components, set the new representative vertex to the \textit{younger} of the two representatives (younger in forward filtration order, i.e., older in reverse filtration order).
    
    \item Change the way that persistence pairs are recorded: Assume a dual edge $c_i$ links two components represented by dual vertices $c_j, c_k$, where $c_j$ is older than $c_k$ in the refined forward filtration order (and hence \textit{younger than $c_k$} in reverse filtration order). If $c_i$ is strictly older than $c_j$ in the non-refined forward order, we record a persistence pair $(c_i, c_j)$.
\end{itemize}

The described algorithm works on the list of $2$-cubes, which it interprets as dual edges.
Since the second dimension is the top dimension, it must enumerate those 2-cubes itself before performing the union-find algorithm.
When it is finished, it performs clearing on this list before it passes it to dimension $1$:
the 2-cubes to be cleared are found by marking every dual edge that merges two distinct components during the union-find algorithm.

\subsection{Computing the Image Barcode}
\label{sec:computing-image-barcode-betti-matching}

As discussed in \cref{sec:induced-matchings-betti-matching}, the Betti matching relies on computing \textit{image barcodes}.
An algorithm to compute image barcodes was first described in \cite{cohen-steiner_persistent_2009} in a specific setting, and more generally in \cite{bauer_efficient_2022}.
The central result for our purposes is described in \cite[Theorem 3.6]{bauer_efficient_2022} and can be paraphrased as follows.

To compute the persistence pairs of $\im \Phi$, where $\Phi$ is the morphism induced by inclusions $D(f_{\mathbf I}) \hookrightarrow D(f_{\mathbf C})$, we can apply the matrix reduction algorithm in a modified form. 
We replace the boundary matrix $D$ with a matrix $D^\Phi$ that is obtained from the boundary matrix of the comparison image $\mathbf{C}$, with the columns ordered by the cube-wise refinement of $\mathbf C$, but the rows ordered by the cube-wise refinement of $\mathbf I$.
We reduce $D^\Phi$ using the standard reduction algorithm to obtain the reduced matrix $R^\Phi$.
If $\pivot_{R^\Phi}(j) = i$, we obtain an image persistence pair $(c_i, c_j)$.
If additionally $f_{\mathbf I}(c_i) < f_{\mathbf C}(c_j)$ it contributes an interval to the image barcode (see \cite[Theorem 3.6]{bauer_efficient_2022}).

As mentioned in \cref{sec:induced-matchings-betti-matching}, the Betti matching uses an extended notion of the image barcode that includes reverse pairs, i.e., image persistence pairs $(c_i, c_j)$ that satisfy $f_{\mathbf I}(c_i) \geq f_{\mathbf C}(c_j)$. 
Hence, we can ignore this condition and not discard such pairs in the algorithm if we want reverse pairs included.

Adjusting the optimized reduction algorithm in dimension 1 of \cref{sec:barcode-dimension-1} for image persistence is straightforward. The required changes are:
\begin{itemize}
    \item The \code{columnsToReduce} come from the comparison image $\mathbf C$, while the boundary cubes (i.e. rows) come from the input image $\mathbf I$. 
    In practice, this is achieved by passing $\mathbf C$'s array to the column enumeration, but passing $\mathbf I$'s array to  \code{enumerateBoundary()}, which, given a $2$-cube $c$, outputs boundary $1$-cubes with locations based on $c$, but birth values read from $\mathbf I$.

    \item In case reverse pairs should \textit{not} be included, we must check the relevant condition: 
    When a column is fully reduced and non-zero, yielding a potential image persistence pair $(c_i, c_j)$, one needs to check that $c_i$ was born before $c_j$ before recording the pair.

    \item In the emergent pair optimization, the condition \\ \code{facet.birth == columnsToReduce[i].birth} needs to be changed: 
    \code{columnsToReduce[i].birth} refers to $\mathbf C$, but needs to refer to $\mathbf I$ to correctly express the emergent pair condition (see \cite[section 3.4]{bauer_efficient_2022}). 
    Instead we look up the birth associated with the location of \code{columnsToReduce[i]} in $\mathbf I$'s array. 
    If reverse pairs should be included, emergent pairs are recorded instead of being skipped.
\end{itemize}%
To compute the image persistence pairs in dimension $0$ and dimension $2$, we can again exploit the extra structure and rely on a union-find algorithm.
In dimension $0$ we initialize a union-find data structure on the vertices of $\mathbf{I}$ and use the refined filtration order of $\mathbf C$ on the edges to link components.
In contrast, in dimension $2$ we initialize a union-find data structure on the dual vertices of $\mathbf{C}$ and use the reverse refined filtration order of $\mathbf I$ on the dual edges to link components.

\subsection{Computing the Betti matching}

In order to compute the Betti matching, we need to compute five different barcodes: 
One for each of the input images $\mathbf I$, $\mathbf J$, the \textit{input barcodes}; one for the comparison image $\mathbf C$, the \textit{comparison barcode}; and two \text{image barcodes} for the inclusion morphism from $\mathbf I$ to $\mathbf C$ and from $\mathbf J$ to $\mathbf C$, respectively.

Recall the structure of the matching (see figure \ref{fig:image-barcode-betti-matching}): 
If a comparison persistence pair $(a^{\mathbf C}, c^{\mathbf C})$ is matched via induced matchings with an input persistence pair $(b^{\mathbf I}, d^{\mathbf I})$ via an image persistence pair $(b^{\mathbf I}, c^{\mathbf C})$, and with an input persistence pair $(b^{\mathbf J}, d^{\mathbf J})$ via an image persistence pair $(b^{\mathbf J}, c^{\mathbf C})$. 
Then the input persistence pairs $(b^{\mathbf I}, d^{\mathbf I})$ and $(b^{\mathbf J}, d^{\mathbf J})$ are matched by the Betti matching.

Given the input persistence pairs of $\mathbf I$ and $\mathbf J$, the comparison persistence pairs and the image persistence pairs, it is straightforward to compute the Betti matching. 
One can create a lookup table of image persistence pairs for each input, mapping the death cube to the birth cube, and a lookup table for the persistence pairs of each input, mapping the birth cube to the pair. 
Then one can traverse the comparison pairs, follow the lookup maps to pairs $(b^{\mathbf I}, d^{\mathbf I})$ and $(b^{\mathbf J}, d^{\mathbf J})$ and match them. 
If any of the lookups fail, there are no matched input pairs induced by this comparison pair.

We show pseudocode for the general approach that we just described in listing \ref{alg:computeBettiMatching} in the appendix, using dimension $1$ as an example. 
The same structure can be used for dimension $0$ and dimension $2$. 
Keep in mind that for the clearing optimization, lists of columns to reduce are passed between the dimensions, a detail we omit here for conciseness.

%
% 9_computeBettiMatching.tex listing
%

After computing the matches, it is also straightforward to compute the unmatched input pairs, which we will also use in the Betti matching loss.
It should be noted that if one is only interested in the matching, not the persistence pairs, the Betti matching can be integrated more tightly with the persistence pair algorithms to gain a bit of efficiency. 
For example, we can refrain from computing lists of persistence pairs and instead place them directly into the lookup maps. 
We can also exploit the structure of the dimension $0$ and dimension $2$ union-find algorithms. 

In dimension 0, the comparison pairs and the image pairs can be computed jointly since both union-find algorithms traverse the edges of the comparison complex.
For the computation of comparison pairs we intialize a union-find data structure $\text{UF}_{\mathbf C}$ and for the computation of the image pairs we initialize a union-find data structure $\text{UF}_{\mathbf I}$.
Note that an edge $c_j$ links two distinct components in $\text{UF}_{\mathbf C}$ if and only if it links two distinct components in $\text{UF}_{\mathbf I}$.
When this happens, a persistence pair $(c_i,c_j)$ of $\mathbf C$ and an image persistence pair $(c_{\tilde{i}},c_j)$ is found and a match can directly be established with the previously computed input pairs if $f_{\mathbf C}(c_i) < f_{\mathbf C}(c_j)$.
At this point we can also check the condition that $f_{\mathbf I}(c_{\tilde i}) < f_{\mathbf C}(c_j)$ if we want to exclude reverse pairs.
In dimension $2$, the input and image pairs can be computed jointly, and afterwards the matching can be performed during the comparison pair computation.

\subsection{Implementation-Level Optimizations}
\label{sec:implementation-optimizations}

The optimized reduction algorithms described in pseudocode in the previous sections are close to our implementation, however a few details specific to the C++ implementation have been left out or simplified for presentation purposes. 
This section describes those details and measures the impact some of the optimizations have on performance. 
It is worth a considerable amount of effort to optimize the run time of the Betti matching computation since it will run thousands of times in a neural network training and in our experience dominates the run time of a training iteration.
Furthermore, it cannot trivially be parallelized except to a limited degree and in current implementations for cubical complexes (including ours) is CPU-bound (whereas neural network trainings can effectively utilize GPUs).

At this point, it should also be noted that our C++ codebase contains a number of deactivated algorithmic optimizations, such as \textit{apparent pairs} and \textit{image persistence clearing}. 
These have not proven beneficial for the performance on our benchmarking data and hence have been deactivated in the default setting (guarded by \code{\#ifdef} directives).
The presence of these optimizations (even if deactivated) does at some points require a more complicated code structure, and this is a cautious note that the structure described here slightly differs from the structure of the C++ code.

One of the implementation-level optimizations we will describe here is adopted from Cubical Ripser \cite{kaji_cubical_2020} and concerns the cube representation. 
The other implementation-level optimizations are improvements over Cubical Ripser or independent from Cubical Ripser as they concern the Betti matching. 
New improvements are the implementation of the \code{cache} and \code{columnIndexByPivot} dictionaries as flat arrays, the caching of working boundaries as lists instead of queues, sorting optimizations when enumerating edges, and parallelization of some of the five barcode computations that do not depend on each other. 
We will measure their impact on performance in \cref{sec:implementation-performance-experiments}. 
An implementation-level optimization independent of Cubical Ripser is the caching of collected faces while searching for emergent pairs (Cubical Ripser does not implement emergent pairs).
Besides runtime optimizations, we describe how the three types of barcode computations (input, comparison and image barcodes) can be unified to avoid code duplication.

\paragraph{Cube Representation}
We represent cubes of any dimension by the same data structure \code{Cube}. 
The dimension of the cube does not need to be stored, as in all our algorithms the dimension of each cube is clear from context.
We only store the cube's birth intensity, its coordinates $(x, y, z)$ and its \textit{type} (recall from \cref{sec:persistent-homology-digital-images} that the coordinates refer to the smallest-coordinate vertex that is a face, and the type indicates the direction of 1-cubes and 2-cubes and is set to 0 for 0-cubes and 3-cubes). 
To save memory, the coordinates and the type are compressed into a single 64-bit integer \code{index} (20 bits per coordinate and 4 bits for the type) and recovered via bit-shifting operations when required, a functionality which the \code{Cube} class offers via instance methods \code{x()}, \code{y()}, \code{z()}, \code{type()}. 
This increased memory efficiency is reflected in a noticeable runtime speedup. 
Furthermore, \code{index} serves as a sorting key for the lexicographic tiebreaking.
This representation is adopted from Cubical Ripser \cite{kaji_cubical_2020}, however Cubical Ripser uses a different lexicographic tie-breaking order (type, $z$, $y$, $x$).

\paragraph{Implementing Dictionaries as CubeMaps}
The dimension $1$ reduction algorithm in \ref{sec:barcode-dimension-1} uses dictionaries indexed by cubes for the \code{cache} and \code{columnIndexByPivot}. 
In Cubical Ripser, this is implemented via \code{std::unordered\_map} using the \code{Cube::index} as keys.
We replace this with a flat array, exploiting the fact that we can injectively map cubes by their coordinates to an index in $\{0, ..., N_1 \cdot N_2 \cdot N_3 \cdot T-1\}$, where $N_1, N_2, N_3$ are the input image dimensions and $T \in \{1, 3\}$ is the number of possible cube types, depending on the dimension.
It turns out that, although \code{std::unordered\_map} has on-average $O(1)$ search and insertion operations, using a flat array is considerably faster for our algorithm. 
We wrap the flat array in a template class \code{CubeMap<\_dim>} which abstracts the indexing and chooses $T$ based on the template argument \code{\_dim}.

\paragraph{Caching faces collected during emergent pairs search}
In the emergent pair optimization we traverse the facets of the initial $2$-cube, add them to the working boundary and continue to the next $2$-cube if we find an emergent pair. 
Adding a facet to the working boundary is moderately expensive, and if an emergent pair is found, the working boundary will not be needed.
To avoid this cost, one could either traverse the enumerator a second time after no emergent pair has been found and only then add the facets to the working boundary, or one could cache the facets into a \code{std::vector} (which has a less expensive append operation) and transfer it to the working boundary if no emergent pair is found. 
We experimentally find the last option to be most efficient and use it in our implementation.

\paragraph{Caching working boundaries as lists}
Cubical Ripser caches working boundaries in priority queue format after cleaning them (de-duplicating modulo 2). 
This is wasteful, as it involves not only push and pop operations during cleaning, but also a copy as well as pop operations when using a cached column. 
It is less expensive to convert the priority queue to a \code{std::vector} when caching and traversing its contents when reading from cache. 
We already allude to this technique with the function \code{toList} in the pseudocode for \code{CubeXorQueue}.

\paragraph{Stable sort and binary sort optimization}
The barcode computations in each dimension are preceded by cube enumeration and lexicographic sorting by \code{Cube::birth} and \code{Cube::index}.
This sorting makes up a small, but non-negligible part of the run time because of the high number of cubes. 
Cubical Ripser sorts the list of cubes with \code{std::sort}, performing the lexicographic comparison in a custom \code{CubeComparator}.
However, by the nature of the nested loop, the cubes are already sorted lexicographically by coordinates and type, and it suffices to perform a \code{std::stable\_sort} on the \code{Cube::birth}, which is a bit faster.

A second sorting optimization is specific to our supervised machine learning setting: 
As we will see in the next chapter, a typical pair of inputs will consist of a binary label volume and a non-binary prediction volume. 
A list of cubes with birth values restricted to $\{0, 1\}$ can be sorted highly efficiently using \code{std::stable\_partition}. 
We make use of this trick by checking during cube enumeration if all birth values are binary, using the partition sort in this case, and falling back to the standard stable sort otherwise. 
The sorting hence becomes significantly faster for the one barcode computation (out of five) that involves only the label input.

\paragraph{Parallelization of barcode computations}
Recall that the dimension $1$ Betti matching algorithm involves five barcode computations (two input barcodes, one comparison barcode, two image barcodes) which are largely independent.
This means that the barcode computations in dimension $1$ can be performed in parallel, which we implement by spawning a \code{std::async} task for each computation. 
The only dependence between these is that the comparison barcode computation can perform clearing to speed up the image barcode computations (see \cite[Proposition 3.10]{bauer_efficient_2022}). 
If this clearing is enabled, not all five computations can be started at once, but the input and comparison barcode computations can be started, and once the comparison barcode computation finishes, the image barcode computations can run.
We find it beneficial for performance to have the clearing enabled even if it prevents some parallelism. 
This can be explained by the observation that the image barcode computations are more expensive than the other computations and dominate the run time, so the clearing saving them some work has a greater impact than the full parallelization would have.

\paragraph{Unifying barcode computations}
The algorithms for input, comparison and image barcode computations are closely related. 
The input and comparison computations only differ by an additional clearing performed in the comparison case (and by some of the deactivated optimizations that we do not discuss), and the image barcode computation differs by recording persistence in a different way, detecting emergent pairs differently, and also an additional clearing step.
To avoid code duplication, all three can be implemented in one method, parametrized by an enum that specifies the mode. 
To avoid the overhead of checking for the mode in every iteration of a tight loop, we make the mode a template argument which gives an optimizing compiler the chance to optimize out any such checks.

\subsection{Performance experiments}
\label{sec:implementation-performance-experiments}

We benchmark the performance of the Betti-matching-3D implementation over several datasets. 
The first series of experiments concerns the implementation optimizations to the Betti matching described in the previous section. 
The optimizations were added in sequence, corresponding to different commits in our codebase, and we measure their influence to the performance by compiling and running the respective commits on a number of benchmarking datasets. 
Our code is compiled with clang, version \code{12.0.0-3ubuntu1-20.04.5}. 
All experiments are performed on a server with a AMD EPYC 7452 32-Core CPU. 
To reduce the run time variance that is due to system load, we run the full set of commit and dataset combinations ten times, each repetition in a randomized order, and report the mean and standard deviation of the ten runs.
For the benchmarking datasets, we focus on input volumes relevant to our deep-learning-based segmentation application that we will introduce in the next chapter. 
These input volumes come in pairs, consisting of a binary label volume and a prediction volume taking values in $[0, 1]$. 
We measure both the performance on the full volumes as well as on random patches of a size similar to what we use in our application setting, and compute the Betti matching on each input pair.

The results are shown in \cref{tab:betti-matching-3d-optimizations}. 
We find that the implementation optimizations give a significant speedup, in the most drastic case by a factor of more than 3.
Furthermore, the speedup is highly dataset-dependent and is bigger for large volumes, but still significant for smaller volumes. 
The single most effective optimization seems to be replacing \code{unordered\_map} by \code{CubeMap}. 
The parallelization also gives significant speed gains, especially for large volumes - for small volumes, we suspect that a single image barcode computation dominates the run time and makes the parallelization less effective. 
The influence of the sorting and caching as list optimizations is smaller and more volume-dependent (although in all cases a net reduction in runtime).

The second series of experiments concerns the performance difference between Cubical Ripser and the most optimized version of Betti-matching-3D. 
We use Betti-matching-3D in barcode-only mode, i.e., we only perform one of the five barcode computations required for the Betti matching. 
We benchmark the barcode computation on single label and prediction volumes as well as two volumes used in the Cubical Ripser publication.

The results, displayed in \cref{tab:betti-matching-3d-vs-cripser}, show that Betti-matching-3D typically outperforms Cubical Ripser significantly (it should be noted that the comparison is fair, as no parallelization is used in Betti-matching-3D for single barcodes). 
The exception is a label volume from the VesSAP dataset, where Cubical Ripser is significantly faster. 
To investigate this difference, we create a version of Betti-matching-3D that uses the same tiebreaking order as Cubical Ripser (type, $z$, $y$, $x$ instead of $x$, $y$, $z$, type). 
This version performs better on the volume in question, although still worse than Cubical Ripser, and performs worse than our main implementation on all other inputs. 
Again, the relative difference is bigger on larger volumes.

We show the Betti numbers of our benchmarking datasets in \cref{tab:dataset-betti-numbers}. 
Furthermore, we benchmark the run time of subprocedures of Betti-matching-3D and present the results in \cref{fig:subprocedure-timings}.

\section{Topologically Faithful 3D Segmentation}
\label{section:topological_3d_segmentation}

We are now ready to apply the persistent homology concepts introduced in the previous chapter to 3D image segmentation.
The task in (binary) image segmentation is to predict for each voxel in an input image, such as a medical scan, whether it belongs to a certain defined structure or not.
Neural networks for image segmentation are commonly trained with a volumetric loss such as the Dice loss (the voxel-wise $F_1$ score). 
A network trained with the Dice loss may not respect topological properties particularly well, as the Dice loss does not take the spatial relationship of voxels into account.
Major topological errors may be caused by a small number of mispredicted voxels, such as a thin bridge falsely connecting two distinct objects. 
Several approaches have  emerged recently that directly address topological correctness, often by explicitly modeling topological errors in a topological loss function (see \cite{hu_topology-preserving_2019}, \cite{stucki_topologically_2023}, \cite{hu_structure-aware_2022}, \cite{gupta_learning_2022}, \cite{gupta_topology-aware_2023}).

The first topological segmentation loss function that uses persistent homology is the TopoNet loss proposed in \cite{hu_topology-preserving_2019}. 
It is derived from the Wasserstein distance which is applied in persistent homology to measure the similarity of persistence diagrams (see \cite{cohen-steiner_lipschitz_2010}).
Denote by $\Dgm(\mathbf I), \Dgm(\mathbf J)$ the persistence diagrams belonging to images $\mathbf I, \mathbf J$. 
In our case, $\mathbf I$ represents a prediction map and $\mathbf J$ represents a label map. 
Keep in mind that the persistence diagrams have the infinitely many points on the diagonal included by definition. 
The TopoNet loss minimizes a loss term over all bijections $\gamma: \Dgm(\mathbf I) \to \Dgm(\mathbf J)$ that respect the dimension (see \cite{stucki_topologically_2023}). 
The loss is defined as
\begin{equation}
    \ell_{\text{Topo}}(\mathbf I, \mathbf J) = \min_{\gamma} \sum_{q \in \Dgm(\mathbf I)} \lVert q - \gamma(q)\rVert_2^2
    \label{eq:toponetloss}
\end{equation}
The loss function in effect pushes the prediction persistence pair points towards their matched target points, making their birth and death values more similar to the target birth and death.
Some points are matched to points on the diagonal, and we refer to these points as \textit{unmatched points}. 
The optimal function $\gamma^\ast$ matches topological features based on their locations in the persistence diagram, but doesn't take into account the spatial relationship of features. 
The idea pursued by \cite{stucki_topologically_2023} is to replace $\gamma$ with the spatially more meaningful Betti matching $\mu_{\mathbf I, \mathbf J}$ that we introduced in \cref{sec:induced-matchings-betti-matching}. 
The Betti matching loss is defined as
\begin{gather}
   \ell_{\text{BM}}(\mathbf I, \mathbf J)
   = \sum_{q \in \text{Dgm}(\mathbf I)} 2 \lVert q - \mu_{\mathbf I, \mathbf J}(q) \rVert_2^2.
   \label{eq:betti-matching-loss}
\end{gather}
The Betti matching loss is defined with a factor of two that allows to interpret its value as the number of unmatched features in a setting where $\mathbf I$, $\mathbf J$ are binary.

The Betti matching loss decomposes into two parts, representing matched and unmatched pairs:
\begin{align*}
    \ell_{\text{BM}}(\mathbf I, \mathbf J)
     =~&~ \ell_{\scriptscriptstyle\text{BM}}^{\scriptscriptstyle\text{matched}}(\mathbf I, \mathbf J) + \ell_{\scriptscriptstyle\text{BM}}^{\scriptscriptstyle\text{unmatched}_{\mathbf I}}(\mathbf I, \mathbf J) + \ell_{\scriptscriptstyle\text{BM}}^{\scriptscriptstyle\text{unmatched}_{\mathbf J}}(\mathbf I, \mathbf J) \\[0.3cm]
    =~&~~~~
    \sum_{\substack{\scriptscriptstyle (q_b, q_d) \in \text{Dgm}(\mathbf I),\\\scriptscriptstyle\text{matched with}\\\scriptscriptstyle(p_b, p_d) \in \text{Dgm}(\mathbf J)}} 2\left((q_b-p_b)^2+(q_d-p_d)^2\right)
    \label{eq:betti-matching-loss-decomposed}
    \\[0.2cm]
    &+ \sum_{\substack{\scriptscriptstyle (q_b, q_d) \in \text{Dgm}(\mathbf I),\\\scriptscriptstyle \text{unmatched}}} (q_b - q_d)^2
    + \sum_{\substack{\scriptscriptstyle (p_b, p_d) \in \text{Dgm}(\mathbf J),\\\scriptscriptstyle \text{unmatched}}} (p_b - p_d)^2 \notag
\end{align*}%
The Betti matching loss captures that components matched across the two images should have a similar contrast range, and that unmatched components should vanish.
In the segmentation setting, $\mathbf I \in [0, 1]^{N_1 \times N_2 \times N_3}$ will be a predicted likelihood map, and $\mathbf J \in \{0, 1\}^{N_1 \times N_2 \times N_3}$ will be a binary label.
Since in this setting the foreground corresponds to voxels with a high intensity value, it is appropriate to consider \textit{superlevel sets} instead of the \textit{sublevel sets} we described previously, and use a filtration that start at the maximum and grow towards the minimum value. 
The superlevel setting can be reduced to the sublevel setting by flipping signs. 
In the superlevel setting, features are born at high intensity values and die at lower intensity values. 
Hence, all features of $\mathbf J$ will be born at 1 and die at 0. The matched features of the predicted likelihood map $\mathbf I$ will have birth and death values inbetween, which will be pushed towards those endpoints, i.e. to have sharper contrasts and less uncertainty. 
Meanwhile each unmatched feature of $\mathbf I$, with birth $b$ and death $d$, will be pushed towards $(\frac{b+d}{2}, \frac{b+d}{2})$, making the contrast between the feature and its surroundings vanish.

The Betti matching loss can be used to train a segmentation network using gradient descent. 
It admits a gradient under the mild assumption that $\mathbf I$ does not have duplicate values (and hence $\mu_{\mathbf I, \mathbf J}$ is constant in a small neighborhood around $\mathbf I$)  \cite[Equation (4)]{stucki_topologically_2023}. 
In practice, the gradient can be computed using automatic differentiation and the computation of the matching $\mu_{\mathbf I, \mathbf J}(q)$ can be treated as a black box.
It is worth considering what the gradient describes: since the birth and death of a topological feature are each tied to a single voxel, the gradient prescribes, per feature, a change in only two voxels with respect to $\mathbf I$ \footnote{with respect to the network parameters, it might influence a larger number of voxels.}.

The Betti matching loss can reinforce desired existing features, as well as make undesired features (i.e. noise) vanish. 
However, it is not designed to create new topological features, e.g. create a new connected component or connect a loop  (unless by accident). 
Furthermore, it will typically not quickly change the prediction over large regions due to it focusing on individual pairs of pixels per topological feature. 
For these reasons, it is necessary to combine $\ell_{\text{BM}}$ with a standard volumetric loss, such as the Dice loss. 
In our experiments, we use a linear combination of the Betti matching loss and the Dice loss 
\[\ell_{\text{DiceBetti}(\alpha)}(\mathbf I, \mathbf J) = \alpha \ell_{\text{BM}}(\mathbf I, \mathbf J) + \ell_{\text{Dice}}(\mathbf I, \mathbf J) \ ,\] 
the \textit{DiceBetti loss} (cf. \cite{stucki_topologically_2023}).
The working hypothesis is that the Betti matching loss can refine the results of pure volumetric training, correcting its topological mistakes and enforcing the topologically correct parts. 
In such a combination, it is also conceivable that the Betti matching loss does create new features: imagine the example of a loop, where a network trained with the Dice loss confidently predicts a large portion the loop to be present (say with voxel intensities of 0.99), but not the remaining voxels that would close the loop. 
In a typical segmentation setting, there is an input signal for those remaining voxels, and the output will still have a slightly increased likelihood for those voxels - say with intensities of 0.05, while the true background voxels are 0.01. This setting would still lead to a matched feature, with birth at 0.05 and death at 0.01: hence the Betti matching loss could pick up the weak signal and reinforce it to gain high likelihood values for the full loop. 

\section{Experimentation and Results}

\begin{figure*}[b]
    \centering
    \hspace{1.5cm}
    \begin{minipage}[t]{0.35\textwidth}
    \subfloat[ Dice-trained network]{
        \centering
        \includegraphics[width=\textwidth]{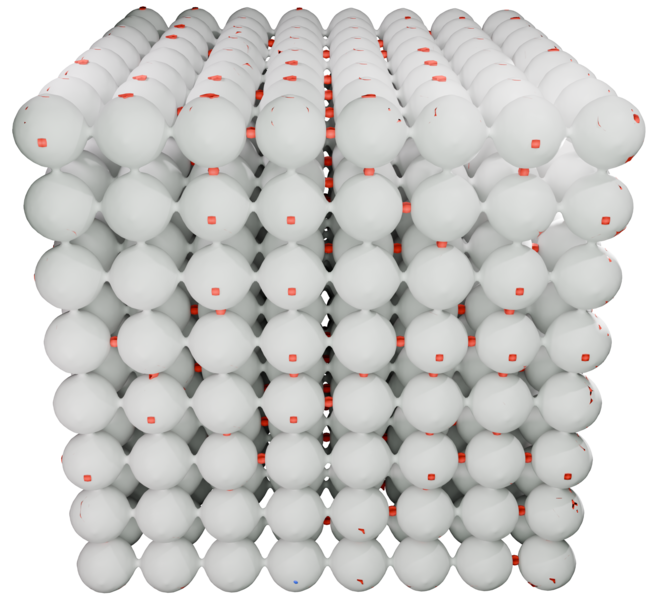}
        } 
    \end{minipage}
    \hfill
    \begin{minipage}[t]{0.35\textwidth}
    \subfloat[DiceBetti($\alpha$=0.1)-trained network]{
        \centering
        \includegraphics[width=\textwidth]{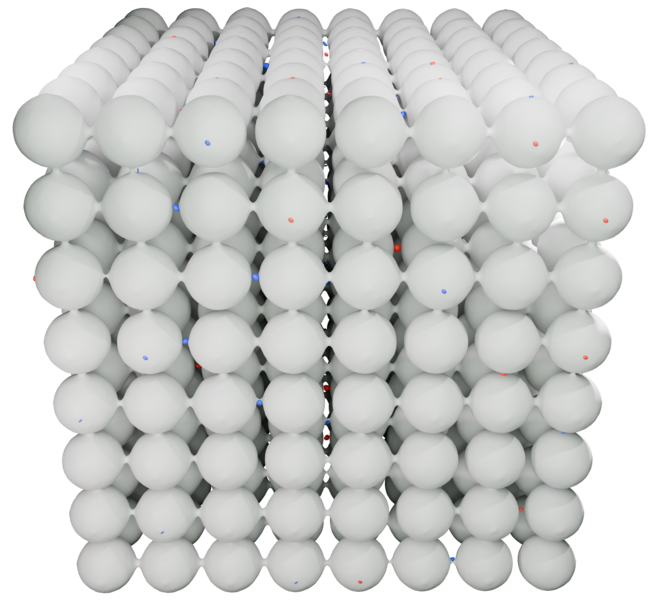}     
        }
    \end{minipage}
    \hfill
    
    \caption{Visualization of prediction quality on synthetic spheres data. False positive predictions marked in red, false negative predictions marked in blue. Gray indicates correct predictions.}
    \label{fig:synblob-confusion}
\end{figure*}

\begin{table*}[]
\caption{Results of networks trained with differently weighted DiceBetti-Losses or with baseline losses, reported on test set (VesSAP-Subset, CREMI-Boundaries, Synthetic spheres, Synthetic spheres with holes) or validation set (TopCoW-MR, TopCoW-CT, BBBC027-Downscaled), respectively. The best result per dataset and metric is highlighted in bold. $\scriptstyle {}^{\ast}$ Computed on patches}
\label{tab:main-results-table}
\centering
%\resizebox*{!}{0.95\textheight}{
\begin{tblr}{colspec={c|l|r|r|r|r|r},stretch=0}
\SetCell[r=2]{halign=c,valign=m} \textbf{Dataset} & 
\SetCell[r=2]{halign=c,valign=m} {\textbf{Training and}\\\textbf{model selection loss}} &
\SetCell[r=2]{halign=c,valign=m} {\textbf{Selected}\\\textbf{epoch}}&
\SetCell[r=2]{halign=c,valign=m} {\textbf{Dice}\\\textbf{loss} $\shortdownarrow$}
& \SetCell[r=2]{halign=c,valign=m} {\textbf{Betti matching}\\\textbf{loss} $\shortdownarrow$}  & \SetCell[r=2]{halign=c,valign=m} {\textbf{ClDice}\\\textbf{metric} $\shortuparrow$
}  & \SetCell[r=2]{halign=c,valign=m} {\textbf{TopoNet}\\\textbf{loss ($\mathbf \times 2$)} $\shortdownarrow$}
\\\\
\cline{1-7}
\SetCell[r=8]{}
{BBBC027-Downscaled} & Dice & 1900 & 0.131 & 837.500 & 0.858 & 479.750 \\
 & DiceBetti($\alpha=0.002$) & 1000 & 0.127 & 511.500 & 0.879 & 330.000 \\
 & DiceBetti($\alpha=0.01$) & 1200 & 0.137 & 380.375 & \textbf{0.886} & 147.876 \\
 & DiceBetti($\alpha=0.025$) & 700 & 0.167 & 372.500 & 0.878 & 142.500 \\
 & DiceBetti($\alpha=0.1$) & 700 & 0.230 & \textbf{364.875 }& 0.832 & 111.376 \\
 & DiceBetti($\alpha=0.5$) & 1200 & 0.534 & 437.500 & 0.570 & \textbf{59.750} \\
 & ClDice($\alpha=0.1$) & 1800 & \textbf{0.126} & 922.500 & 0.879 & 444.500 \\
 & HuTopo($\alpha=0.1$) & 2400 & 0.141 & 406.250 & 0.880 & 227.000 \\
\cline{1-7}
\SetCell[r=8]{}
{VesSAP-Subset} & Dice & 300 & \textbf{0.286} & 3551.000 & \textbf{0.783} & 2295.000 \\
 & DiceBetti($\alpha=0.002$) & 350 & 0.301 & 4425.500 & 0.714 & 3121.500 \\
 & DiceBetti($\alpha=0.01$) & 400 & 0.315 & 2242.500 & 0.723 & 830.500 \\
 & DiceBetti($\alpha=0.025$) & 400 & 0.325 & 2916.5 & 0.699 & 1385.500 \\
 & DiceBetti($\alpha=0.1$) & 450 & 0.318 & 3593.5 & 0.705 & 5127.500 \\
 & DiceBetti($\alpha=0.5$) & 400 & 0.338 & 2491.000 & 0.700 & 973.000 \\
 & ClDice($\alpha=0.1$) & 950 & 0.297 & \textbf{2078.500} & 0.750 & \textbf{672.500} \\
 & HuTopo($\alpha=0.1$) & 350 & 0.301 & 4425.500 & 0.714 & 3121.500 \\
\cline{1-7}
\SetCell[r=8]{}
{CREMI-Boundaries ${}^{\ast}$}
& Dice & 2200 & 0.166 & 44.787 & 0.873 & 29.298 \\
 & DiceBetti($\alpha=0.002$) & 400 & 0.231 & 30.380 & 0.825 & 14.046 \\
 & DiceBetti($\alpha=0.01$) & 500 & 0.220 & 25.010 & 0.834 & 10.142 \\
 & DiceBetti($\alpha=0.025$) & 400 & 0.246 & 27.569 & 0.789 & 11.482 \\
 & DiceBetti($\alpha=0.1$) & 1500 & 0.230 & \textbf{24.086} & 0.808 & \textbf{9.074} \\
 & DiceBetti($\alpha=0.5$) & 2500 & 0.258 & 24.975 & 0.818 & 9.362 \\
 & ClDice($\alpha=0.1$) & 2200 & \textbf{0.161} & 41.264 & \textbf{0.873} & 26.268 \\
 & HuTopo($\alpha=0.1$) & 300 & 0.238 & 40.194 &  0.820 & 23.118 \\
\cline{1-7}
\SetCell[r=8]{}{Synthetic spheres} & Dice & 1550 & \textbf{0.008} & 449.500 & 0.998 & 449.500 \\
 & DiceBetti($\alpha=0.002$) & 1200 & 0.016 & 163.250 & 0.999 & \textbf{42.250} \\
 & DiceBetti($\alpha=0.01$) & 2160 & 0.007 & 132.750 & \textbf{0.999} & 79.750 \\
 & DiceBetti($\alpha=0.025$) & 2190 & 0.008 & 132.000 & 0.999 & 68.000 \\
 & DiceBetti($\alpha=0.1$) & 2010 & 0.042 & \textbf{121.000} & 0.998 & 55.000 \\
 & DiceBetti($\alpha=0.5$) & 2850 & 0.178 & 169.750 & 0.984 & 139.250 \\
 & ClDice($\alpha=0.1$) & 2720 & 0.009 & 571.750 & 0.998 & 571.750 \\
 & HuTopo($\alpha=0.1$) & 2930 & 0.015 & 469.500 & 0.997 & 56.500 \\
\cline{1-7}
\SetCell[r=8]{}
{Synthetic spheres\\with holes} & Dice & 1530 & 0.042 & 485.000 & 0.993 & 199.500 \\
 & DiceBetti($\alpha=0.002$) & 2370 & 0.041 & 220.250 & 0.998 & 100.750 \\
 & DiceBetti($\alpha=0.01$) & 2170 & \textbf{0.033} & 167.750 & 0.998 & 50.750 \\
 & DiceBetti($\alpha=0.025$) & 2720 & 0.037 & \textbf{157.500} &\textbf{ 0.998} & \textbf{53.000} \\
 & DiceBetti($\alpha=0.1$) & 2990 & 0.060 & 176.000 & 0.998 & 96.000 \\
 & DiceBetti($\alpha=0.5$) & 2290 & 0.104 & 200.750 & 0.998 & 135.750 \\
 & ClDice($\alpha=0.1$) & 1950 & 0.042 & 486.000 & 0.993 & 323.500 \\
 & HuTopo($\alpha=0.1$) & 350 & 0.106 & 606.250 & 0.993 & 484.250 \\
\end{tblr}
%}
\end{table*}

\subsection{Datasets}
\label{sec:datasets}
We experimented with our developed Betti matching segmentation loss on several 3D medical imaging datasets, including vessel segmentation, cell segmentation, and two synthetic datasets designed to emphasize topological errors in training. For details on the datasets please refer to the Supplement section \ref{app:data}.

\subsection{Model Training}
\label{sec:training}
For our experiments we use a 3D-Unet \cite{cicek_3d_2016,kerfoot_left-ventricle_2019}, with hidden channel sizes $[16, 32, 64, 128, 256]$, downsampling strides of 2, two residual units, and $\text{PReLU}$ activations ($\text{PReLU}(x) = x \cdot \mathds{1}_{[0, \infty)}(x) + c\,x \cdot \mathds{1}_{(-\infty, 0)}(x)$ for a learnable $c$).
To achieve a consistent and meaningful initialization, we pre-train with the  Dice loss for 350 epochs before switching to the target loss. On the synthetic sphere datasets, we train for 3000 epochs and validate every 10 epochs, while on all other datasets, we train for 2500 epochs and validate every 50 epochs.
We train networks using the Dice-Betti matching loss $\ell_{\text{DiceBetti}(\alpha)}$ for a range of different values of the Betti matching loss weighting factor $\alpha$. For baselines, we train: 1) Dice $\ell_{\text{Dice}}$, 2) centerline Dice loss $\ell_{\text{clDice}}$ \cite{shit_cldice_2021}, an established topological loss which computes the Dice loss on \textit{centerlines} of the prediction and label, and 3) TopoNet $\ell_{\text{Topo}}$ as proposed by \cite{hu_topology-preserving_2019}, already described in equation (\ref{eq:toponetloss}).

\subsection{Model Evaluation}
\label{sec:eval}
%We evaluate models trained with the DiceBetti loss for different values of the Betti matching weighting factor $\alpha$ and compare them to models trained with the baseline losses across the seven datasets. 
Four datasets are evaluated on dedicated held-out test sets (VesSAP, CREMI-Boundaries, Synthetic spheres, and Synthetic spheres with holes). 
%The other three datasets are evaluated on the validation sets (TopCoW-MR, TopCoW-CT, BBBC027-Downsampled). 
The other BBBC027-Downsampled datasets is evaluated on the validation sets.
For each of the four training losses, we select the model at the epoch where this loss is optimized on the validation set during training. We compute the evaluation metrics (Dice loss, Betti matching loss, clDice metric, and TopoNet loss) on full volumes, except for CREMI: For CREMI we fall back to patches (covering the full volume), as the volume size makes a non-patched computation infeasible even at evaluation time, and report all metrics averaged over patches. We binarize the prediction before computing the losses, and use the non-relative setting for the Betti matching and TopoNet losses. In this binarized setting, the interpretation of the Betti matching loss is that it counts the number of unmatched features (this interpretation is possible because of the factor of two). The TopoNet loss analogously counts the number of features not matched by the Wasserstein matching, but is defined without the factor of two. To make the two losses directly comparable, we report the evaluation TopoNet loss multiplied by two.

\subsection{Results}
\label{sec:results}

All results are reported in \cref{tab:main-results-table}. We can observe several trends from the results. The clearest trend is that the DiceBetti loss generally yields a big improvement on the evaluation Betti matching loss, compared to Dice-only training. In other words, the topological optimization goal during training also significantly improves the generalization topological error. The DiceBetti trainings outperform all baselines for almost all datasets on this metric, regardless of the value of $\alpha$.
%, and the $\alpha = 0.002$ setting performs equally bad as the TopoNet loss. 
The exception is VesSAP-Subset, where clDice outperforms the DiceBetti loss.
Generally, for VesSAP-Subset, we find our method to be less effective, as in all cases, the selected model comes from an epoch shortly after turning on the Betti matching loss (after the pure-Dice loss warm-up), and in some cases the model selection epoch is 350 (i.e. just before turning on the Betti matching loss component), meaning that training with the Dice-Betti loss leads to lower performance in the validation Betti matching loss for these settings. Since we achieve significantly lower Betti matching losses on the training set, we attribute the poor performance on VesSAP-Subset to overfitting. Possibly the effect also has to do with suboptimal model selection caused by a difficult validation set, as we find test set metrics to be better than validation set metrics. Nevertheless, we significantly outperform the Dice loss on the Betti matching loss for $\alpha=0.01$ and $\alpha=0.025$. The better performance of clDice can be attributed to clDice being designed with tubular structures such as the VesSAP dataset in mind. 

Another less pronounced trend shows in a trade-off between the evaluation Dice performance and the Betti matching loss. Under the hypothesis that the Dice loss is not well aligned with topological correctness, such behavior is expected and we find that in most cases a high $\alpha$ leads to a high Dice loss. For small $\alpha$, the trend is not emphasized. 

%Surprisingly, there are a few datasets where this trend is reversed, such that a higher $\alpha$, and lower Betti matching loss, also coincide with a lower Dice loss. In particular, on TopCoW-CT, some DiceBetti-trained networks have significantly lower evaluation Dice losses. 
%This is strictly a generalization effect, as in training the Dice loss is always lowest when directly trained on. This means that on some data, a DiceBetti-trained network helps the network generalize in terms of Dice loss, and has a certain regularization effect.

\begin{figure}
    \centering
    \caption{Influence of $\alpha$ weight parameter of DiceBetti training loss on evaluation Betti matching loss and Dice loss. $\alpha=0$ corresponds to pure Dice loss training (highlighted in blue). We show results for three datasets.ul}
    \label{fig:evaluation-alpha-vs-metrics}
    
        \vspace{0.15cm}

        \includegraphics[width=0.8\linewidth]{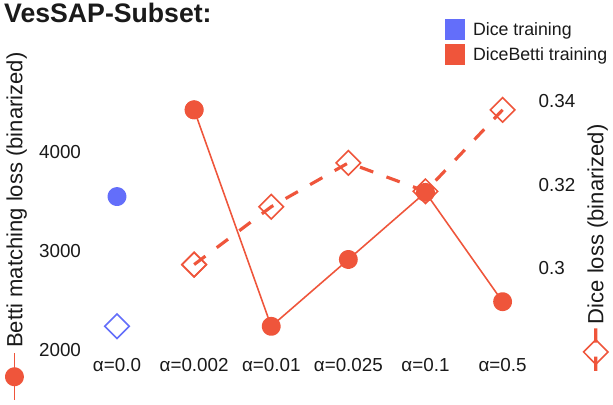}
        
        \vspace{0.2cm}

        \includegraphics[width=0.8\linewidth]{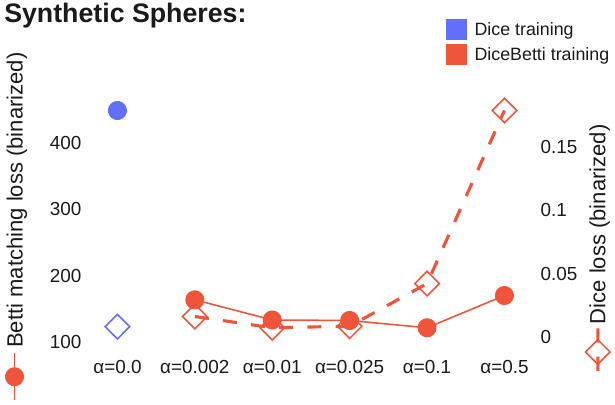}
        \vspace{0.2cm}
        
        \includegraphics[width=0.8\linewidth]{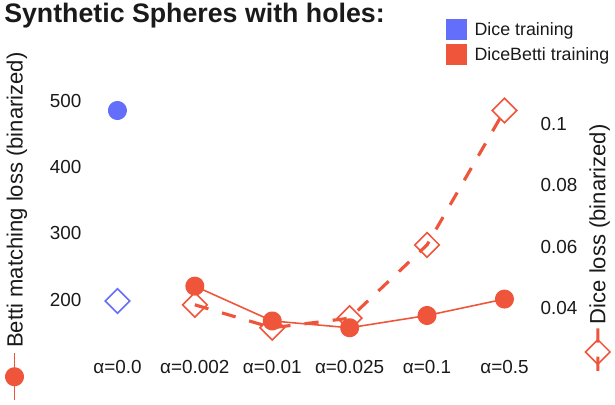}

\end{figure}

The plots in Figure \ref{fig:evaluation-alpha-vs-metrics} show the influence of the $\alpha$ parameter on the evaluation Dice and DiceBettiMatching losses. They also show a significant improvement in the evaluation of Betti matching loss when going from $\alpha=0$ to $\alpha=0.002$ on all datasets except VesSAP. They furthermore show that there is a tradeoff between the two metrics on most, but not all datasets, and that the behavior with respect to $\alpha$ is highly dataset-dependent. Importantly, what constitutes a good value of $\alpha$ depends on the number of topological features in the training patches: The Betti matching loss scales with this number, and $\alpha$ needs to be adjusted accordingly such that the Dice loss still has sufficient influence and the gradients do not explode. 
%This also explains why we can choose relatively large values of $\alpha$ for the TopCoW-Datasets, which have a small number of topological features (see \ref{subfig:topcow}). 
%To get an idea of the number of features in some of our datasets, refer to the Betti numbers in \ref{tab:dataset-betti-numbers} in the previous chapter.

The networks trained with the baseline topological losses generally perform worse with respect to the Betti matching evaluation loss than the DiceBetti trained networks (with the exception of clDice on VesSAP-Subset).
%In some cases, they still perform better than the pure Dice loss, in some they perform worse. With respect to the clDice evaluation metric, the DiceBetti training loss outperforms the baselines as well on most datasets, with a large margin in particular on the TopCoW datasets. Again, 
This appears to be a generalization effect, where the DiceBetti training helps the generalization ability. 
%Somewhat surprisingly, the TopoNet-trained network manages to outperform the DiceBetti-trained networks on the TopCoW datasets with respect to the Dice and clDice losses. Since the matching used by the TopoNet loss is necessarily more arbitrary than the Betti matching, this is likely due to chance, combined with the fact that the relatively small number of features in the TopCoW datasets makes it easier for the TopoNet loss not to mismatch features. 

With respect to the evaluation TopoNet loss, the best scores are also achieved by DiceBetti-trained networks on most datasets. We report the TopoNet loss multiplied by a factor of two such that it is directly comparable to the Betti matching loss. However, one has to keep in mind that the TopoNet loss has the flaw of not taking spatial matches into account: Both the Betti matching loss and the doubled TopoNet loss, if applied to binarized predictions (as is the case in our Table), represent the average number of unmatched features with respect to their respective matchings. However, the TopoNet loss matches features between binarized prediction and target regardless of the location, which causes it to only capture the \textit{difference in the number of features} (per dimension) between the binarized prediction and the target in the evaluation setting. Hence the doubled TopoNet loss is always lower than the Betti matching loss, and the difference becomes bigger if there is a similar number of features not matched by the Betti matching in the prediction and target. The BettiLoss($\alpha$=0.5) training for BBBC027 is an extreme example: there are on average 172 unmatched dimension 0 features in the prediction, and 214.75 in the target, giving a dimension 0 Betti matching loss of $172+214.75 = 386.75$ and a dimension 0 TopoNet loss of $\lvert 172 - 214.75 \rvert = 42.75$. Generally speaking, a low Betti matching loss also implies a low TopoNet loss, but not the other way around. In other words, the Betti matching loss measures a stronger notion of topological correctness.

We also examine the performance on the synthetic sphere datasets, which was deliberately constructed to be very challenging for the Dice training loss.
We find that, the network trained with the Dice loss fails to make a topologically correct prediction on the test set samples. Most drastically, in the dataset without holes in the sphere surfaces, the Dice loss training converges to a model where it predicts all connections between spheres to be present (this behavior was consistently observed on several runs). On the more varied \textit{synthetic spheres with holes} dataset, the Dice loss fares better with respect to connectivity, but tends to fill in the holes in the sphere surfaces. On both datasets, the Betti matching loss makes far fewer topological errors. 
%A plot of false-positive and false-negative predictions of the two losses on a sample volume is shown in \ref{fig:synblob-confusion}.

\section{Discussion}

In this work, we propose an efficient algorithm for the calculation of the Betti Matching, implemented in C++. 
Our proposed solution enables the calculation of the Betti Matching for inputs of arbitrary dimensions, being highly optimized in 1D, 2D, and 3D contexts. 
For wide applicability of the Betti Matching, it is critical that a fast implementation is available. 
To this end, the Betti-Matching-3D implementation is a contribution that enables the previously infeasible 3D-segmentation training with Betti Matching. 
The implementation-level optimizations we propose yield a considerable speedup and significant performance improvements over the state-of-the-art Cubical Ripser implementation.

\paragraph{Limiations}

While our method provides a vastly accelerated computation of barcodes and topological matchings, further optimizations are possible. 
One option is directly  utilizing GPUs to make the implementation less CPU-bound, for example in the style of Ripser++ \cite{bauer_ripser_2021}, which further parallelizes a part of the algorithm, or by delegating the expensive working boundary computations to the GPU. 
A purely CPU-based parallelization is less beneficial to our use case, since in our trainings we can trivially parallelize across multiple CPU cores by computing each instance of a batch in a separate thread.

Considering the loss, additional studies on which loss components are effective on which dataset and their fine-tuning could further improve performance; one option could be individual weighting of the loss components. 
Further, besides the matched/unmatched components, one could introduce weightings for the different dimension components of the loss ($H_0$, $H_1$, $H_2$). 
An important future research direction is to understand the behavior of the Betti matching loss in prototypical scenarios, such as growing or extinguishing a connected component, loop or cavity, which could also lead to a better understanding of the behavior under reweighting of the loss components.

\bibliographystyle{unsrt}  
\bibliography{0_main}
\newpage

%If you have an EPS/PDF photo (graphicx package needed), extra braces are needed around the contents of the optional argument to biography to prevent the LaTeX parser from getting confused when it sees the complicated $\backslash${\tt{includegraphics}} command within an optional argument. (You can create your own custom macro containing the $\backslash${\tt{includegraphics}} command to make things simpler here.)
 
\vspace{11pt}

% \end{IEEEbiography}



\vfill

%% Supplemental Material 
\appendix
% \section{Supplemental Material}

\subsection{Algorithms}
\subsubsection{Unoptimized matrix reduction algorithm for persistence pair computation}~

\label{alg:basicReductionAlgorithm}

\subsubsection{High-level structure of optimized persistence pair algorithm}~

\label{alg:computePersistencePairs_highLevel}

\subsubsection{$H_1$ persistence pair algorithm}~

\label{alg:computePersistencePairsDim1}

\paragraph{Efficient working boundary reduction}~
\begin{lstlisting}[language=python,style=PseudoPython]
function reduceWorkingBoundaryBy(workingBoundary: CubeXorQueue, column: Cube2, cache: Dict[Cube2, [Cube1]], gridComplex):
  cachedBoundary: [Cube1] = cache[column] # returns (*$\bot$*) if not present
  if cachedBoundary != (*$\bot$*):
    for facet in cachedBoundary:
      workingBoundary.push(facet)
  else:
    for facet in enumerateBoundary(column, gridComplex):
      workingBoundary.push(facet)
\end{lstlisting}%
\label{alg:reduceWorkingBoundaryBy}

\paragraph{CubeXorQueue data structure}~
\begin{lstlisting}[language=python,style=PseudoPython,escapeinside={(*}{*)}]
class CubeXorQueue:
  cubeQueue: PriorityQueue[Cube1] = []

  function push(cube: Cube1):
    cubeQueue.push(cube)

  function getPivot() -> (Cube1 | (*$\bot$*)):
    if cubeQueue.empty(): return (*$\bot$*)

    pivot: Cube1 = cubeQueue.pop()
    while (not cubeQueue.empty()) and cubeQueue.top() == pivot:
      cubeQueue.pop()
      if cubeQueue.empty():
        return (*$\bot$*)
      pivot = cubeQueue.pop()
    cubeQueue.push(pivot)
    return pivot

  # Drains the queue to a list while deduplicating modulo 2
  function toList() -> [Cube1]:
    elements: [Cube1] = []
    while not cubeQueue.empty():
      cube: Cube1 = cubeQueue.pop()
      if (not cubeQueue.empty) and cube == cubeQueue.top():
        cubeQueue.pop()
      else:
        elements.push(cube)
    return elements
\end{lstlisting}
\label{alg:cubeXorQueue}

\paragraph{Emergent pairs optimization}~
\begin{lstlisting}[language=python,style=PseudoPython,firstnumber=10]
# ... (in computePersistencePairsDim1, before inner do-while loop)
    checkEmergentPair = True
    foundEmergentPair = False
    for facet in enumerateBoundaryReverse(columnsToReduce[i], gridComplex):
      if checkEmergentPair and facet.birth == columnsToReduce[i].birth:
        j = columnIndexByPivot[facet]
        checkEmergentPair = False
        if j == (*$\bot$*):
          foundEmergentPair = True
          columnIndexByPivot[facet] = i
      workingBoundary.push(enumerator.previousFace())
    if foundEmergentPair:
      continue # continue to next column to reduce in outer loop
# ... (if no emergent pair found, continue with inner do-while loop)
\end{lstlisting}
\label{alg:emergentPairsOptimization}

\paragraph{Enumerating columns for dimension 0 with clearing}~
\begin{lstlisting}[language=python,style=PseudoPython]
function enumerateDim0Columns(columnIndexByPivot, gridComplex) -> [Cube1]:
  columnsForDim0 = [
    new Cube1(x, y, z, t)
    for x in 1..gridComplex.width for y in 1..gridComplex.height
    for z in 1..gridComplex.depth for t in 1..3
    if cube not in columnIndexByPivot
  ]
  columnsForDim0.sort()
  return columnsForDim0
\end{lstlisting}
\label{alg:enumerateDim0Columns}

\subsubsection{$H_0$ persistence pair algorithm (union find)}~

\label{alg:computePersistencePairsDim0}

\subsubsection{Computing the Betti matching}~
{}
\begin{lstlisting}[language=python,style=PseudoPython,keepspaces=True,escapeinside={(*}{*)}]
function computeBettiMatching(gridComplex1, gridComplex2, gridComplexComparison):
  inputPairs1 = computePersistencePairsDim1(gridComplex1)
  inputPairs2 = computePersistencePairsDim1(gridComplex2)
  comparisonPairs = computePersistencePairsDim1(gridComplexComparison)
  imagePairs1 = computeImagePersistencePairsDim1(gridComplex1, gridComplexComparison)
  imagePairs2 = computeImagePersistencePairsDim1(gridComplex2, gridComplexComparison)

  matchMap1 = {pair.birth: pair for pair in inputPairs1}
  matchMap2 = {pair.birth: pair for pair in inputPairs2}
  matchMapImage1 = {pair.death: pair.birth for pair in imagePairs1}
  matchMapImage1 = {pair.death: pair.birth for pair in imagePairs2}

  matches = []
  for comparisonPair in comparisonPairs:
    pair1 = matchMap1[matchMapImage1[comparisonPair.death]]
    pair2 = matchMap2[matchMapImage2[comparisonPair.death]]
    if pair1 != (*$\bot$*) and pair2 != (*$\bot$*):
      matches.append((pair1, pair2))
  return matches
\end{lstlisting}%
\label{alg:computeBettiMatching}

\subsection{Datasets}
\label{app:data}

\paragraph{VesSAP}
The VesSAP dataset \cite{todorov_machine_2020} consists of light-sheet microscopy images of mouse brains. The labels provide a segmentation of the vessel structure (see Figure X). 
The dataset admits a complicated topological structure where a label map may contain a large number of connected components with complex connectivity structure, as well as a number of cycles.
We use a set of eight volumes ($500\times500\times50$ voxels) and two input channels. Four volumes are training volumes, two are used for validation, and two are used for testing.

\paragraph{BBBC027}
The BBBC027 dataset \cite{svoboda_generation_2011} consists of synthetic colon cell microscopy images. The data exhibits a large number of small connected components that are often in close proximity and separated by thin gaps that are hard to recognize from the input image (see Figure X). 
Because of the large image size, we downsampled the 30 scans and labels in each dimension to a resolution of $32\times257\times325$.
We aim to preserve the connectivity and not falsely merge cells by downsampling; therefore, we pick the minimum value in a $4\times 4 \times 4$ region for both the image and the label. Of the 30 scans, we used 22 for training and 8 for validation.

\paragraph{CREMI (boundaries)}
CREMI is a neuroimaging dataset containing brain tissue scans of \textit{Drosophila melanogaster} \cite{funke_large_2019}. The neurons are segmented instance-wise. We consider the boundaries between neurons as our binary label, leading to a topologically complex label map with a large number of cycles and a moderate number of cavities (see Figure X).
The dataset consists of three large labeled volumes, each with different characteristics. We use volume A for training, volume B for validation, and volume C for testing (with the additional challenge that the different characteristics will make validation and test set predictions harder for the trained network).

\paragraph{Synthetic Data}
In addition, we created two new synthetic datasets with the specific goal of creating a topologically challenging dataset.

% \afterpage{
\begin{figure*}[t]
    % \centering

    % \begin{tblr}{colspec={cccccc}}
        % \SetCell[c=3]{}
        % % \begin{subfigure}[t]{0.49\textwidth}
        \subfloat[VesSAP vessel dataset sample]{
            \centering
            \includegraphics[height=0.485\textwidth]{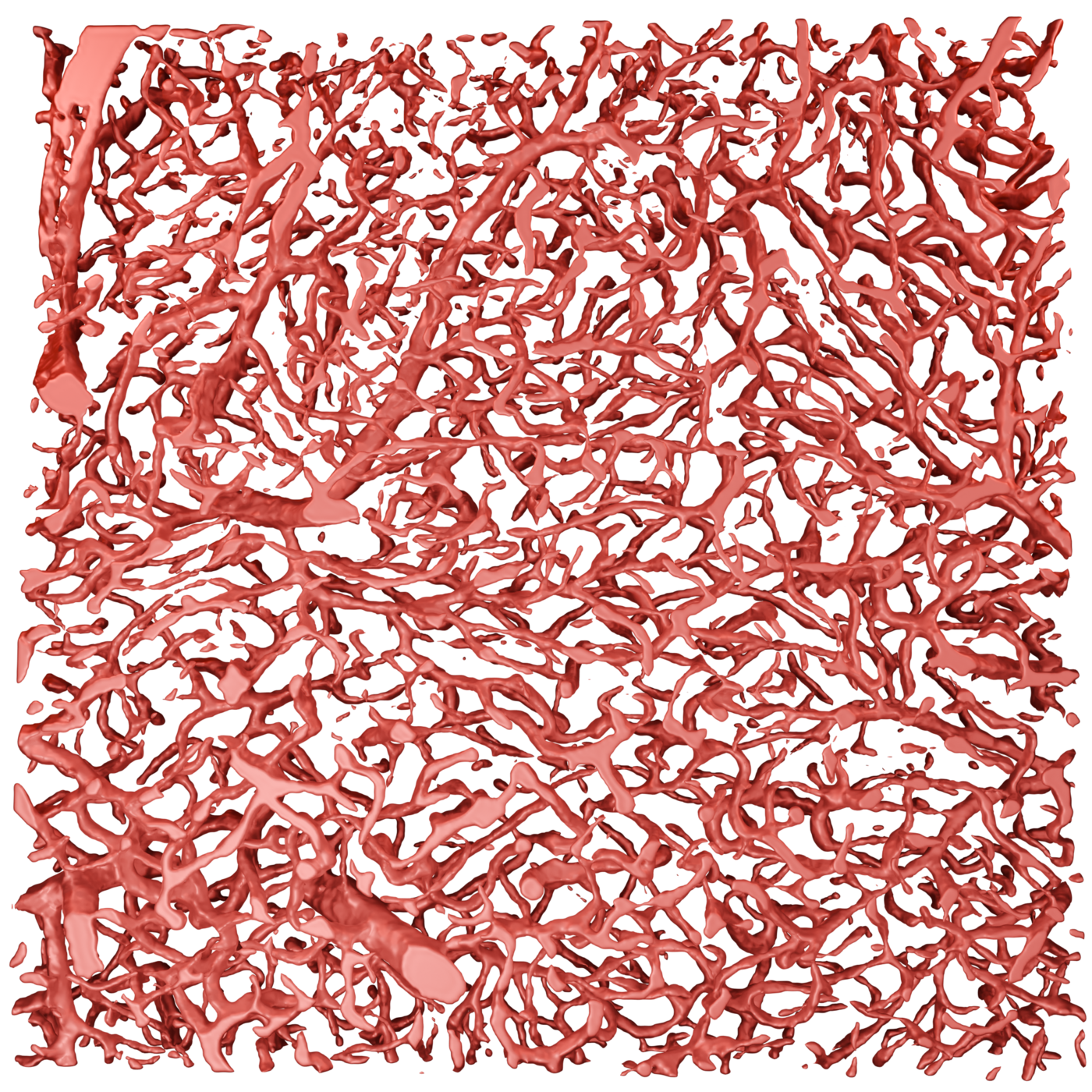}
            % \caption{VesSAP vessel dataset sample}
            % \label{subfig:vessap}
        % \end{subfigure}
        }
        % &&&
        % \SetCell[c=3]{}
        \subfloat[CREMI (boundaries) dataset (50\texttimes800\texttimes800 
        % = 2.0\texttimes3.2\texttimes3.2\textmu m
            patch from 125\texttimes1250\texttimes1250 volume)
        ]{
            \centering
            \includegraphics[height=0.45\textwidth]{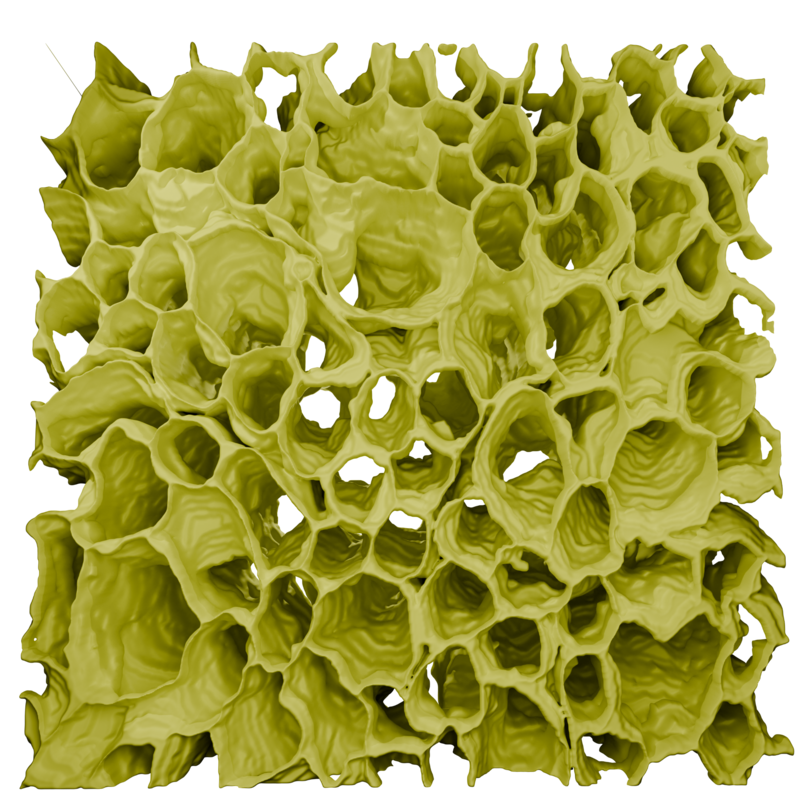}
            % \captionsetup{format=centeredMultilineCaption}
            % \caption{CREMI (boundaries) dataset (50\texttimes800\texttimes800 
            % % = 2.0\texttimes3.2\texttimes3.2\textmu m
            % patch from 125\texttimes1250\texttimes1250 volume)}
            % \label{subfig:cremi}
        }
        \\
    
    % % \hfill
    % % \begin{subfigure}[t]{0.5\textwidth}    
    
    % % \end{subfigure}
    % % \vspace{0.1cm}

    %     \SetCell[c=2]{}
    %     % \begin{subfigure}[b]{0.38\textwidth}
        \subfloat[BBBC027 cell dataset sample]{
            % \centering
            \includegraphics[width=.31\textwidth,valign=t]{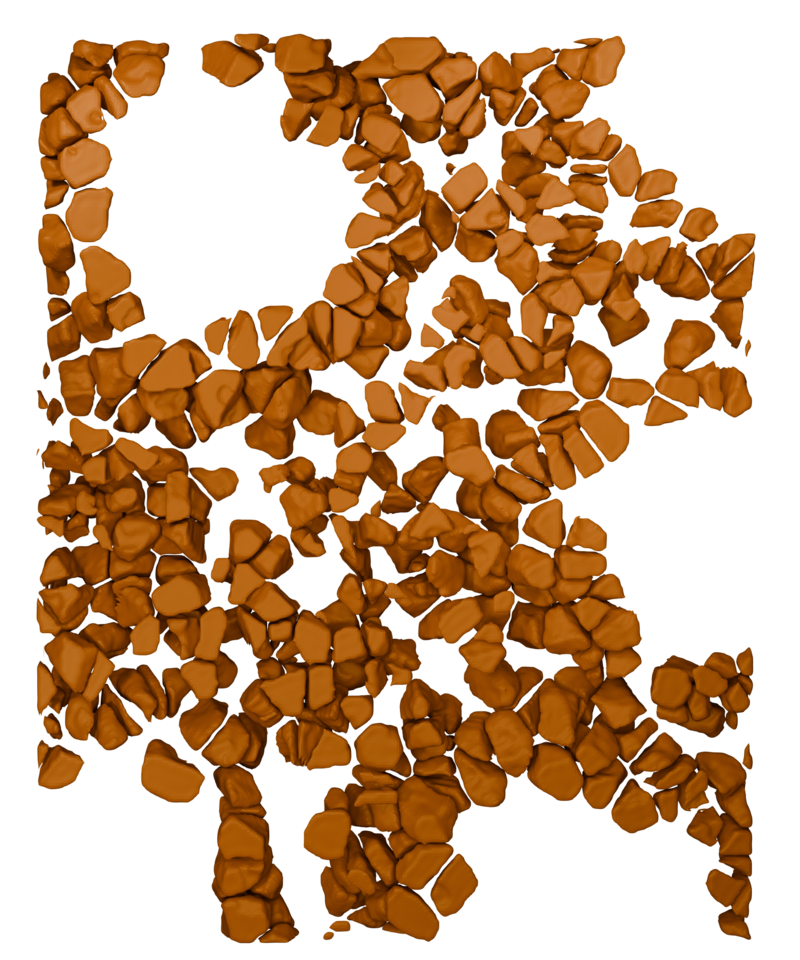}
            % \caption{BBBC027 cell dataset sample}
            % \label{subfig:bbbc027}
        }%
    %     % \end{subfigure}
    %     % \hfill
    %     % \begin{subfigure}[b]{0.3\textwidth}
    %     &&
    %     \SetCell[c=2]{}
        \subfloat[Synthetic spheres\ mesh/{label} map slice)]{
    %         % \centering
            \begin{tabular}{c}
    %         % \hfill
            \vphantom{\includegraphics[width=.31\textwidth]{figures/datasets/bbbc027-example.png}}

            \includegraphics[width=0.31\textwidth,valign=t]{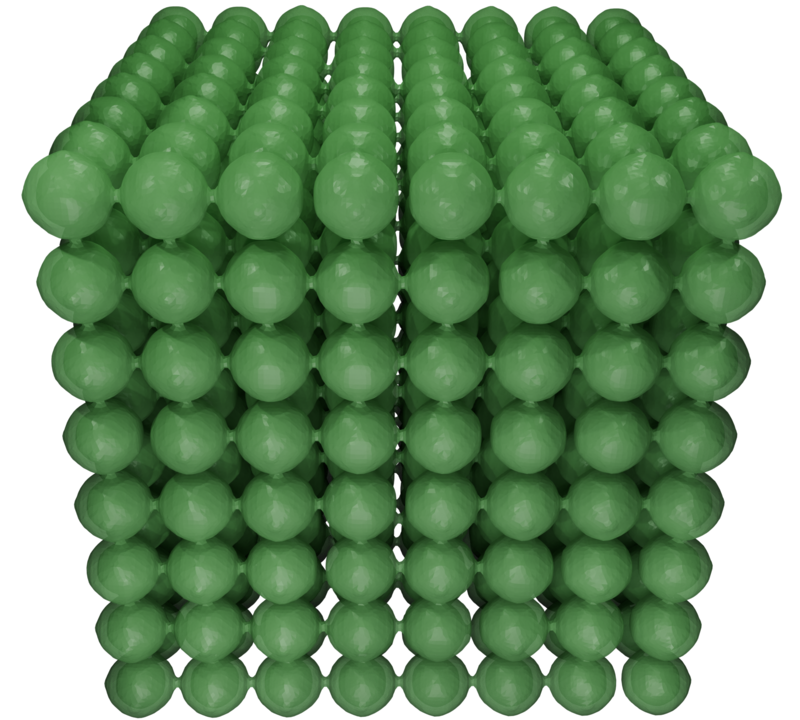}
    %         % \hfill
    \\
    %         % \vspace{0.4cm}
            \includegraphics[width=0.31\textwidth,valign=t]{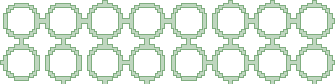}
    %         % \vspace{0.2cm}
    %         % \captionsetup{format=centeredMultilineCaption}
    %         % \caption{Synthetic spheres\\(mesh/{label} map slice)}
    %         % \label{subfig:synblob}
    %
            \end{tabular}
        }
    %     % \end{subfigure}
    %     % \hfill
    %     % \begin{subfigure}[b]{0.3\textwidth}
    %     &&
    %     \SetCell[c=2]{}
        %
        \subfloat[Synth. spheres with holes (mesh/label map slice)]{
            % \centering
            % \hfill
            \vphantom{\includegraphics[width=.31\textwidth]{figures/datasets/bbbc027-example.png}}
            \begin{tabular}[t]{c}
                \includegraphics[width=0.31\textwidth,valign=t]{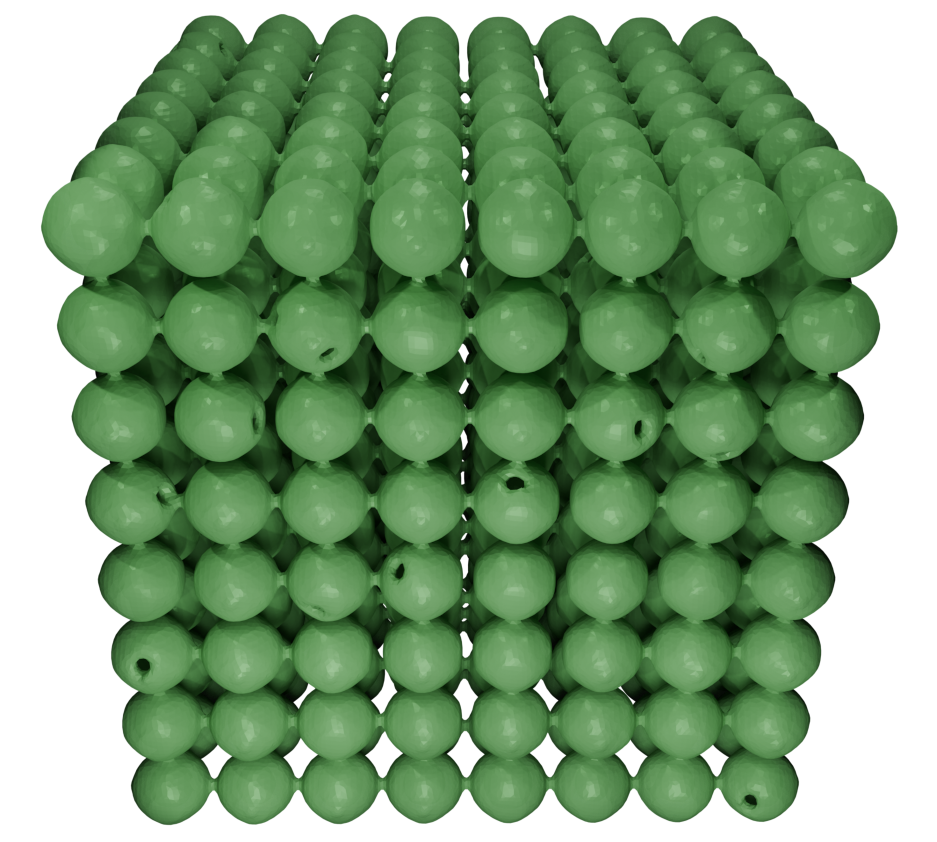}
                \\
            % \hfill
    
            % \vspace{0.4cm}
                \includegraphics[width=0.31\textwidth,valign=t]{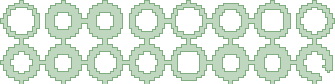}
            % \vspace{0.2cm}
            % \captionsetup{format=centeredMultilineCaption}
            % \caption{Synth. spheres with holes (mesh/label map slice)}
            % \label{subfig:synblobholes}
        % \end{subfigure}
            \end{tabular}
    }
    
    % \end{tblr}
    % \vspace{0.1cm}
    \caption{Samples of the topologically complex 3D datasets used in our experiments (smoothed meshes from label maps in Slicer3D, rendered with Blender).}.
    \label{fig:dataset-table}

\end{figure*}

% }

The first dataset, the \textit{synthetic spheres} dataset, consists of hollow spheres with a constant surface thickness that are regularly arranged in a grid. Each pair of adjacent spheres is connected by a thin bridge with a certain probability (see Figure). In the input images, the voxels making up the boundaries of the spheres have a base intensity of 1, those that make up the bridges, and the inside of the hollow spheres have a base intensity of -1, and the background has a base intensity of 0. The base intensities are overlayed with $\mathcal N(0, 1)$ Gaussian noise.
The goal is that the spheres themselves are relatively easy to segment - they are regularly spaced and uniform in appearance - but the links between them are not, as they are made up of few pixels, largely obscured by noise. Yet correctly predicting the links is crucial for getting the topology correct. The task is made harder by the fact that the correct bridges (label 1) and the sphere interiors (label 0) have the same local appearance (a base intensity of -1). The hypothesis is that Dice loss training will focus on correctly predicting the spheres and their interior, which make up a large proportion of the volume, and largely disregard the crucial bridges as they do not influence the Dice score by a lot.
We create 20 sample volumes following this process, each with a resolution of $128\times128\times128$ containing a grid of $8\times8\times8$ spheres. Of the 20 volumes, 12 are used for training, 4 for validation, and 4 for testing.

The second synthetic dataset, \textit{synthetic spheres with holes}, is a variation of the first synthetic dataset: we give spheres a more varied appearance by choosing a random thickness, and spheres may have one or more holes poked in their surface, such that they do not describe a cavity (see Figure X). The number of holes per sphere follows a geometric distribution.%, where $\sfrac{1}{4}$ of spheres have at least one hole,  $\sfrac{1}{16}$ have at least two holes, etc.
Using this process, we create another 20 volumes with the same resolution, grid size, and train/validation/test split as above.

\newcommand{\patchscalefactor}{0.75}
\newcommand{\stimes}{{\mkern-2mu\times\mkern-2mu}}

\begin{table*}[]
    \centering
    \resizebox{0.99\textwidth}{!}{
    \begin{tblr}{colspec={r|c|c|c|c|c|c||c},stretch=1}
        \SetCell[r=2]{valign=m} {\textbf{Input}\\\textbf{(prediction/label)}}
        & \SetCell[r=2]{valign=m}{\textbf{Before}\\\textbf{optimization}}
        & \SetCell[r=2]{valign=m}{\textbf{Unified}\\\textbf{computePairs(\,)}}
        & \SetCell[r=2]{valign=m}{\textbf{CubeMap}\\\textbf{optimization}}
        & \SetCell[r=2]{valign=m}{\textbf{Sorting}\\\textbf{optimizations}}
        & \SetCell[r=2]{valign=m}{\textbf{Cache boundaries}\\\textbf{as lists}}
        & \SetCell[r=2]{valign=m}{\textbf{Parallelization}}
        & \SetCell[r=2]{valign=m} {\textbf{Performance}\\\textbf{gain}} \\\\
        \cline{1-8}
        VesSAP \scalebox{\patchscalefactor}{($500\stimes500\stimes50$)} & 304.67$\pm$8.1s & 312.89$\pm$7.1s & 193.27$\pm$2.5s & 181.53$\pm$1.7s & 172.92$\pm$2.4s & 83.84$\pm$2.5s & +263.4\% \\
        SynBlobDataHoles \scalebox{\patchscalefactor}{($128\stimes128\stimes128$)} & 214.04$\pm$0.8s & 213.87$\pm$1.4s & 193.26$\pm$1.5s & 180.73$\pm$33.1s & 140.59$\pm$0.9s & 87.13$\pm$0.9s & +145.7\% \\
        BBBC027-Downscaled \scalebox{\patchscalefactor}{($32\stimes257\stimes325$)} & 49.55$\pm$0.7s & 48.98$\pm$0.5s & 30.66$\pm$0.9s & 27.66$\pm$2.9s & 24.33$\pm$0.7s & 20.02$\pm$0.9s & +147.5\% \\
        SynBlobData \scalebox{\patchscalefactor}{($128\stimes128\stimes128$)} & 45.49$\pm$0.5s & 45.67$\pm$0.8s & 32.07$\pm$0.5s & 30.36$\pm$0.6s & 25.03$\pm$0.6s & 15.96$\pm$0.7s & +185.0\% \\
        VesSAP \scalebox{\patchscalefactor}{($96\stimes96\stimes48$)} & 6.92$\pm$0.6s & 6.53$\pm$0.6s & 3.30$\pm$0.2s & 2.80$\pm$0.4s & 2.80$\pm$0.3s & 1.83$\pm$0.2s & +278.1\% \\
        SynBlobDataHoles \scalebox{\patchscalefactor}{($32\stimes32\stimes32$)} & 236$\pm$11ms & 236$\pm$10ms & 185$\pm$7ms & 166$\pm$7ms & 151$\pm$7ms & 128$\pm$12ms & +84.4\% \\
        SynBlobData \scalebox{\patchscalefactor}{($32\stimes32\stimes32$)} & 220$\pm$9ms & 220$\pm$11ms & 171$\pm$3ms & 149$\pm$6ms & 141$\pm$3ms & 116$\pm$8ms & +89.7\% \\
        BBBC027-Downscaled \scalebox{\patchscalefactor}{($32\stimes32\stimes32$)} & 199$\pm$10ms & 204$\pm$11ms & 151$\pm$7ms & 131$\pm$3ms & 127$\pm$7ms & 117$\pm$7ms & +70.1\% \\
    \end{tblr}
    }
    \caption{\small Run time performance of Betti-matching-3D at different levels of implementation optimizations. Execution times are measured on prediction-label pairs from different datasets we used in our segmentation application (see \cref{sec:datasets} for an overview of the datasets). For each dataset we benchmark on a full-size sample volume as well as a patch of this sample volume that reflects the typical patch size used in for our training on the respective dataset. The optimizations were applied to the code sequentially, and each column reflects a commit in our codebase containing all optimizations up to that column. For each input pair and each commit, the table displays the mean and standard deviation of 10 runs.}
    \label{tab:betti-matching-3d-optimizations}
\end{table*}

\begin{table*}[]
    \centering
    \resizebox{0.8\textwidth}{!}{
    \begin{tblr}{colspec={r|c|c|c||c},stretch=1}
        % \toprule
        \SetCell[r=3]{valign=h,} {\textbf{Input}\\\textbf{(prediction/label)}}
        & \SetCell[r=3]{valign=h}{\textbf{Cubical}\\\textbf{Ripser}}
        & \SetCell[r=3]{valign=h}{\textbf{Betti-matching-3D}\\\textbf{\small(barcode only,}\\\textbf{\small type/z/y/x tiebreaking)}}
        & \SetCell[r=3]{valign=h}{\textbf{Betti-matching-3D}\\\textbf{\small(barcode only,}\\\textbf{\small x/y/z/type tiebreaking)}}
        & \SetCell[r=3]{valign=h} {\textbf{Performance}\\\textbf{difference}} \\\\\\
        \cline{1-5}
        % \midrule

        Head Aneurism \scalebox{\patchscalefactor}{($512\stimes512\stimes512$)} & 492.01$\pm$35.1s & 310.50$\pm$48.5s & 269.20$\pm$15.6s & +82.8\% \\
        Bonsai \scalebox{\patchscalefactor}{($256\stimes256\stimes256$)} & 53.07$\pm$4.9s & 30.21$\pm$2.1s & 18.44$\pm$0.9s & +187.8\% \\
        VesSAP Prediction \scalebox{\patchscalefactor}{($500\stimes500\stimes50$)} & 51.36$\pm$4.7s & 41.48$\pm$3.2s & 36.50$\pm$1.3s & +40.7\% \\
        VesSAP Label \scalebox{\patchscalefactor}{($500\stimes500\stimes50$)} & 26.95$\pm$1.6s & 41.81$\pm$1.5s & 59.63$\pm$1.8s & $-$54.8\% \\
        VesSAP Prediction \scalebox{\patchscalefactor}{($96\stimes96\stimes48$)} & 1.07$\pm$0.1s & 1.01$\pm$0.0s & 957$\pm$33ms & +11.8\% \\
        VesSAP Prediction, Sublevel \scalebox{\patchscalefactor}{($96\stimes96\stimes48$)} & 1.00$\pm$0.0s & 978$\pm$49ms & 908$\pm$65ms & +10.1\% \\
        VesSAP Label, Sublevel \scalebox{\patchscalefactor}{($96\stimes96\stimes48$)} & 772$\pm$66ms & 683$\pm$38ms & 468$\pm$33ms & +65.0\% \\
        VesSAP Label \scalebox{\patchscalefactor}{($96\stimes96\stimes48$)} & 750$\pm$65ms & 744$\pm$24ms & 545$\pm$13ms & +37.6\% \\
        % \bottomrule
    \end{tblr}
    }
    \caption{\small Run time performance comparison of Betti-matching-3D with the previous state-of-the-art implementation Cubical Ripser on which Betti-matching-3D is based. We benchmark by running each implementation 10 times and reporting the mean run time and its standard deviation. We called both Cubical Ripser and Betti-matching-3D via their Python bindings, where we used the official \code{cripser} package distributed on PyPI. We also compare to a version of Betti-matching-3D that uses Cubical-Ripser-style tiebreaking (reversing the criteria in the lexicographic order), which performs generally worse but better on one specific input. The \textit{Head Aneurism} \cite{volvisorg_rotational_nodate} and \textit{Bonsai} \cite{volvisorg_ct_nodate} volumes are used because they are also benchmark volumes in the Cubical Ripser paper \cite{kaji_cubical_2020}.}
    \label{tab:betti-matching-3d-vs-cripser}
\end{table*}

\newcommand{\rot}[1]{\rotatebox{90}{\shortstack[l]{#1}}}

\end{document}